\documentclass[twoside, 10pt]{article}

\usepackage[hmargin= 0.75in, height=9.4in]{geometry}
\usepackage{amssymb,amsthm, amsfonts,amsmath,mathrsfs, amsxtra, url, hyperref}

\usepackage{fancyhdr}
\fancypagestyle{plain}{
 \fancyhf{}

}

\renewcommand{\theequation}{\thesection.\arabic{equation}}
 \numberwithin{equation}{section}

\newtheorem {thm}{Theorem}[section]
\newtheorem {prop}{Proposition}[section]
\newtheorem {lemm}{Lemma}[section]

\newtheorem {rem}{Remark}[section]
\newtheorem{assum}{Assumption}[section]
\newtheorem {eg}{Example}[section]

\makeatletter
\newenvironment{highlightequation}{%
  \def\tagform@##1{\maketag@@@{(\ignorespaces##1\unskip\@@italiccorr*)}}%
  \ignorespaces
}{%
  \def\tagform@##1{\maketag@@@{(\ignorespaces##1\unskip\@@italiccorr)}}%
  \ignorespacesafterend
}
\makeatother

\def\ba{\begin{array}}
\def\ea{\end{array}}
\def\bea{\begin{eqnarray}}
\def\eea{\end{eqnarray}}
\def\beas{\begin{eqnarray*}}
\def\eeas{\end{eqnarray*}}
\def\bi{\begin{itemize}}
\def\ei{\end{itemize}}
\def\bc{\begin{cases}}
\def\ec{\end{cases}}

\def\bhe{\begin{highlightequation}  }
\def\ehe{\end{highlightequation}  }

\def\a{\alpha}
\def\ga{\gamma}
\def\d{\delta}
\def\e{\varepsilon}
\def\z{\zeta}
\def\k{\kappa}
\def\l{\lambda}
\def\vr{\varrho}
\def\si{\sigma}
\def\vs{\varsigma}
\def\t{\tau}

\def\o{\omega}
\def\vf{\varphi}

\def\D{\Delta}
\def\G{\Gamma}

\def\O{\Omega}

\def\U{\Upsilon}

\def\bF{{\bf F}}

\def\bz{{\bf 0}}

\def\bd{{\bf d}}

\def\cA{{\cal A}}

\def\cD{{\cal D}}
\def\cE{{\cal E}}
\def\cF{{\cal F}}
\def\cG{{\cal G}}

\def\cL{{\cal L}}

\def\cN{{\cal N}}

\def\cP{{\cal P}}

\def\cS{{\cal S}}
\def\cT{{\cal T}}
\def\cU{{\cal U}}

\def\cW{{\cal W}}
\def\cX{{\cal X}}
\def\cY{{\cal Y}}

\def\hC{\mathbb{C}}

\def\hE{\mathbb{E}}

\def\hN{\mathbb{N}}

\def\hP{\mathbb{P}}
\def\hQ{\mathbb{Q}}
\def\hR{\mathbb{R}}
\def\hS{\mathbb{S}}

\def\sB{\mathscr{B}}

\def\sE{\mathscr{E}}

\def\sL{\mathscr{L}}

\def\sN{\mathscr{N}}

\def\sP{\mathscr{P}}

\def\sU{\mathscr{U}}

\def\sX{\mathscr{X}}
\def\sY{\mathscr{Y}}

\def\fF{\mathfrak{F}}

\def\fM{\mathfrak{M}}

\def\fP{\mathfrak{P}}

\def\fp{\mathfrak{p}}

\def\fm{\mathfrak{m}}

\def\fk{\mathfrak{k}}

\def\ft{\mathfrak{t}}

\def\ti{\n \times \n}
\def\oti{\n \otimes \n}
\def\df{\n := \n}
\def\ls{\n \le \n}
\def\gs{\n \ge \n}
\def\={\n = \n}
\def\+{\n + \n}
\def\-{\n - \n}
\def\ins{\n \in \n}

\def\sb{\n \subset \n}
\def\>{\n > \n}
\def\<{\n < \n}

\def\({\textnormal{(}}
\def\){\textnormal{)}}
\def\[{[\n[}
\def\]{]\n]}

\def\no{\noindent}

\def\ss{\smallskip}

\def\q{\quad}
\def\qq{\qquad}

\def\n{\negthinspace}
\def\dn{\n \n}
\def\tn{\n \n \n}

\def\ol{\overline}
\def\ul{\underline}
\def\ua{\mathop{\uparrow}}
\def\da{\mathop{\downarrow}}

\def\wt{\widetilde}
\def\wh{\widehat}
\def\fra{\mathfrak{a}}

\def\pas{{\hbox{$\hP-$a.s.}}}

\def\hb{\hbox}
\def\dis{\displaystyle}
\def\cd{\cdot}
\def\cds{\cdots}

\def\fa{\,\forall \,}

\def\es{\emptyset}

\def\b1{{\bf 1}}
\def\qed{\hfill $\Box$ \medskip}

\def\esssup{\mathop{\rm esssup}}
\def\liminf{\mathop{\ul{\rm lim}}}
\def\limsup{\mathop{\ol{\rm lim}}}

\newcommand{\lsup}[1]{ \underset{#1}{\limsup}}
\newcommand{\linf}[1]{ \underset{#1}{\liminf}}
\newcommand{\lmt}[1]{ \underset{#1}{\lim}}
\newcommand{\lmtu}[1]{ \underset{#1}{\lim} \n \ua \,}
\newcommand{\lmtd}[1]{ \underset{#1}{\lim} \n \da \,}

\begin{document}

 \title{\bf  On the Robust Dynkin Game
  \thanks{We are grateful to   Jianfeng Zhang for insightful comments.}
 }

\author{
  Erhan Bayraktar\thanks{ \noindent Department of
  Mathematics, University of Michigan, Ann Arbor, MI 48109; email:
{\tt erhan@umich.edu}.}  \thanks{E. Bayraktar is supported in part by the National Science Foundation under DMS-1613170,
 and by the Susan M. Smith Professorship.  Any opinions, findings, and conclusions
 or recommendations expressed in this material are
those of the authors and do not necessarily reflect the views of the National Science Foundation.} $\,\,$,
$~~$Song Yao\thanks{
\noindent Department of
  Mathematics, University of Pittsburgh, Pittsburgh, PA 15260; email: {\tt songyao@pitt.edu}. } \thanks{S. Yao is supported in part by the National Science Foundation under DMS-1613208}}
\date{}

\maketitle

 \begin{abstract}

  We analyze a robust version of the Dynkin game over  a set $\cP$ of mutually singular probabilities.
  We first prove that   conservative player's   lower and upper value   coincide (Let us denote the value by $V$).
  Such a result   connects
  the robust Dynkin game with second-order doubly reflected backward stochastic differential equations.
  Also, we show that  the   value process $V$   is a submartingale
  under an appropriately  defined   nonlinear expectation $\ul{\sE}  $
  up to the first time  $\tau_*$   when $V$ meets the lower payoff process $L$.
  If the probability set $\cP$ is  weakly compact,
  one can  even  find an optimal triplet  $ (\hP_*, \tau_*, \ga_*) $  for the value $V_0$.

     The mutual singularity of  probabilities in $\cP$ causes major technical   difficulties.
  To deal with them, we use some new methods including
  two   approximations with respect to the set   of  stopping times.

 \smallskip   {\bf Keywords:}\;  robust Dynkin game, nonlinear expectation,
 dynamic programming principle, controls in weak formulation,
 weak stability under pasting,     martingale approach,
 path-dependent stochastic differential equations with controls, optimal triplet,
 optimal stopping with random maturity.

\end{abstract}

 \smallskip


  \section{Introduction}

  We analyze  a continuous-time {\it robust}  Dynkin game   with respect to  a non-dominated set $\cP$
  of mutually singular   probabilities on the canonical space $\O$  of continuous paths.  In this game, Player 1,
  who negatively/conservatively thinks that the {\it Nature} is also against her,
  will receive the following payment from Player 2
  if the two players choose $\tau \in \cT$ and $\ga \in \cT$ respectively to quit the game:
  \beas
   R(\tau,\ga) : = \int_0^{\tau \land \ga} \n g_s ds +
     \b1_{\{\tau \le \ga\}}  L_\tau + \b1_{\{ \ga < \tau \}}  U_\ga .
  \eeas
  Here $\cT$ denotes the set of all stopping times with respect to the natural filtration $\bF$ of the canonical
 process $B$, and the running payoff $g$, the terminal payoff $L \le U$ are
 $\bF-$adapted processes uniformly continuous in sense of \eqref{eq:aa211}.

 As   probabilities in $\cP$ are mutually singular, one can not define  the conditional expectation of
 the nonlinear expectation $ \underset{\hP \in \cP  }{\inf} \hE_\hP [\cd] $,
 and thus Player 1's  lower  value process $\ul{V}$ and upper  value process $\ol{V}$,   in  essential extremum  sense.
 Instead, we use shifted processes and regular conditional probability distributions
 (see Section \ref{subsec:preliminary} for details)  to define
  \beas 
  \q     \ul{V}_t(\o) \n :=  \n   \underset{\t  \in \cT^t }{\sup} \,   \underset{\ga  \in \cT^t }{\inf} \,
   \underset{\hP \in \cP(t,\o) }{\inf} \, \hE_\hP  \big[  R^{t,\o} (\t,\ga)    \big] , \q
   \ol{V}_t (\o)  \n  :=  \n     \underset{\hP \in \cP(t,\o) }{\inf} \,  \underset{\ga  \in \cT^t }{\inf} \,
  \underset{\t  \in \cT^t }{\sup} \,   \hE_\hP  \big[  R^{t,\o} (\t,\ga)    \big]   , \q (t,\o) \in [0,T] \n \times \n \O .
 \eeas
 Here $\cT^t$ denotes the set of all stopping times with respect to the natural filtration $\bF^t$
 of the shifted canonical process $B^t$
  on the shifted canonical space  $\O^t $,
  $\cP (t,\o)$ is a path-dependent  probability set which includes all regular conditional probability distributions
 stemming from   $\cP$ (see (P2)),
 and $   R^{t,\o} (\tau,\ga) \n  := \n  \int_t^{\tau \land \ga} \n g^{t,\o}_s ds  \n + \n
     \b1_{\{\tau \le \ga\}}  L^{t,\o}_\tau  \n + \n  \b1_{\{ \ga < \tau \}}  U^{t,\o}_\ga $.

  In Theorem \ref{thm_RDG}, we demonstrate that Player 1's lower and upper value processes    coincide
 and thus she has a value process
 $V_t (\o)   \n =  \n  \ul{V}_t(\o)   \n =  \n  \ol{V}_t(\o)$, $(t,\o)  \n \in \n  [0,T]  \n \times \n  \O  $
  in the robust  Dynkin game.
 We also see in Theorem \ref{thm_RDG}  that   the first time $\tau_*$ when  $V$ meets $L$
 is an optimal stopping time for Player 1, i.e.
 \bea \label{eq:RDOSRT}
 V_0    \n =  \n   \underset{\ga \in \cT }{\inf} \,  \underset{\hP \in \cP}{\inf} \,  \hE_\hP  \big[ R(\tau_*,\ga)   \big]  ,
 \eea
 and that processes $V_t   \n +  \n  \int_0^t g_s ds $, $t   \n \in  \n  [0,T]$ is a submartingale
 under  the pathwise-defined nonlinear expectation
 $  \ul{\sE}_t [\xi] (\o)   \n :=  \n    \underset{\hP \in \cP(t,\o) }{\inf} \hE_\hP [\xi^{t,\o}]  $, $(t,\o)   \n \in  \n  [0,T]   \n \times  \n  \O$
 up to time $\tau_*$.

Since a Dynkin game is actually a coupling of two  optimal stopping problems, the {\it martingale}  approach
introduced by Snell \cite{Snell_1952} to solve the  optimal stopping problem was later extended to Dynkin games,
see e.g. \cite{Neveu_1975, Bismut_DG_1979, ALM_1982, Lepeltier_Maingueneau_1984, Morimoto_1984}.
 In the current paper,
 we will adopt a generalized martingale method with respect to the nonlinear expectations
 $  \ul{\sE} \n = \n \{\ul{\sE}_t\}_{t \in [0,T]} $.
 The mutual singularity of probabilities in $\cP$ gives rise to some major  technical hurdles:
 First,  no dominating probability in $\cP$ means that we do not have  a dominated convergence theorem
 for the  nonlinear expectations  $  \ul{\sE}  $.
 Because of this,   one can not follow the classic  approach
 for Dynkin games to obtain the  $\ul{\sE}-$martingale property of $V_\cd \n + \n \int_0^\cd g_s ds $.
 Second,   we do not have a measurable selection theorem for stopping strategies,
 which complicates the proof of the  dynamic programming principle.

      Our martingale approach starts with a dynamic programming principle (DPP) for process $\ol{V}$.
      The ``subsolution" part of DPP (Proposition \ref{prop_DPP})
      relies on a ``weak stability under pasting" assumption (P3) on the   probability class
   $\{ \cP(t,\o)\}_{(t,\o) \in [0,T] \times \O}$,
   which  allows  us to    construct   approximating measures
   by pasting together  local $\e-$optimal probabilistic models.
  We show in Section~\ref{sec:example} that (P3),
  along with our   other assumptions  on the  probability class, are satisfied
  in the case of some path-dependent SDEs with controls,
  which represents a large class of models on simultaneous drift and volatility uncertainty.
  We demonstrate that
  the  ``supersolution" part of the DPP (Proposition \ref{prop_DPP2})  by employing a countable dense subset $\G$
   of $\cT^t$ to construct  a suitable approximation.
   This dynamic programming result implies the continuity of process $\ol{V}$ (Proposition \ref{prop_conti_V}),
   which plays a crucial role in the approximation scheme (to be described in the following paragraph)
   for proving Theorem \ref{thm_RDG}.

 The key to Theorem \ref{thm_RDG} is the $\ul{\sE}-$submartingality of
 process $ \big\{ \ol{V}_t \n + \n \int_0^t g_s ds \big\}_{t \in [0,T]} $ up to $\tau_*$.
 Inspired by Nutz and Zhang \cite{NZ_2012}'s idea on using stopping times with finitely many values for approximation,
 we define an approximating sequence of value processes $V^n$'s to $\ol{V}$ by
 \beas
    V^n_t (\o)     :  =      \underset{\hP \in \cP (t,\o)  }{\inf} \,  \underset{\ga  \in \cT^t }{\inf} \,
  \underset{\t  \in \cT^t(n) }{\sup} \,   \hE_\hP  \big[  R^{t,\o} (\t,\ga)    \big]
  \le \ol{V}_t (\o) , \q (t,\o) \n \in \n [0,T] \n \times \n \O ,
  \eeas
  where $ \cT^t(n)$ collects all $ \cT^t-$stopping times taking values in   $\big\{ t \n \vee \n (i2^{-n}T)\big\}^{2^n}_{i=0}$.
  By (P3),  Proposition \ref{prop_DPP} still holds for $V^n$, which leads to that for any $\d \n > \n 0$ and $k \n \ge \n n$,
   the process    $ \big\{ V^n_t \n + \n \int_0^t g_s ds \big\}_{t \in [0,T]} $ is an  $\ul{\sE}-$submartingale
   over the grid $\{i2^{-k}T \}^{2^k}_{i=0} $ up to the first time $\nu^{n,\d}$ when $V^n$    meets $L \n + \n \d$
   (see \eqref{eq:dc021}).    Letting $k \n \to \n \infty$, $n \n \to \n \infty$ and then $\e \n \to \n 0$, we
   can deduce from $\lmtu{n \to \infty} V^n \n = \n \ol{V}$ (Proposition \ref{prop_grid})
     and the continuity of $\ol{V}$ that
     the  process $ \big\{ \ol{V}_t \n + \n \int_0^t g_s ds \big\}_{t \in [0,T]} $
     is an $\ul{\sE}-$submartingale up to $\tau_*$. Theorem \ref{thm_RDG} then easily follows.
     It is worth pointing out that our argument does not require  the payoff processes to be bounded.

     At the cost of some additional conditions such as the weak compactness of $\cP$
     and the stronger pasting condition of \cite{STZ_2011b}
     (all of which are satisfied for {\it controls of weak formulation}, see Example \ref{eg_weak_formulation}),
     we can apply the main result of \cite{RDOSRT}
     to find in Theorem \ref{thm_triplet} a pair $(\hP_*, \ga_*)   \n \in  \n  \cP   \n \times  \n  \cT$ such that
         \bea   \label{eq:optimal_triplet}
      V_0   =       \hE_{\hP_*}  \big[ R (\tau_*,\ga_*)   \big]  .
       \eea

  \no {\bf Relevant Literature.}
Since its introduction by   \cite{Dynkin_1967},
Dynkin games have been analyzed  in discrete and continuous-time models for decades.
Bensoussan and Friedman  \cite{Friedman_1973, Bensoussan_Friedman_1974, Bensoussan_Friedman_1977}
first analyzed the games in  the setting of  Markov diffusion processes by means of
variational inequalities and free boundary problems. Bayraktar and S\^{i}rbu in \cite{MR3162260} had a fresh  look at this problem using the Stochastic Perron's method (a verification approach without smoothness).
For a more general class of reward processes  martingale approach was developed under Mokobodzki's condition
 (see e.g. \cite{Neveu_1975, Bismut_DG_1977, Bismut_DG_1979, ALM_1982}) and
   certain regularity assumption on  payoff processes (see e.g. \cite{Lepeltier_Maingueneau_1984,KQC_2014}).

Cvitani\'c and Karatzas \cite{Cvitanic_Karatzas_1996} connected Dynkin games to
  backward stochastic differential equations (BSDEs) with two reflecting barriers $L$ and $U$.
  Along with the growth of  the BSDEs theory, Dynkin games have attracted much   attention
  in the probabilistic  framework with Brownian filtration, see e.g.
   \cite{Hamadene_Lepeltier_Wu_1999, Hamad_Lepeltier_2000, Hamadene_Hassani_2005,  Hamadene_2006, Xu_2007,
   Hamadene_Hdhiri_2006, HRZ_2009, Buckdahn_Li_3, Lp_DRBSDE,  DRBSDEUI}. Among these works,
    \cite{Hamadene_Hassani_2005, Hamadene_Hdhiri_2006, HRZ_2009, Buckdahn_Li_3, Lp_DRBSDE,  DRBSDEUI}
    only require  ``$L \n < \n U$" rather than Mokobodski's condition     via a penalization method.

In Mathematical Finance, the theory of Dynkin games can be applied to
 pricing and hedging game options (or Israeli options) and their derivatives, see \cite{Kifer_2000, Ma_Cvitanic_2001, Kallsen_Kuhn_2004, Hamadene_2006, Ekstrom_2006, Dolinsky_2014} and the references in the survey paper  \cite{Kifer_2013_survey}.
 Also,  \cite{Ekstrom_2006, Alvarez_2008} analyzed the sensitivity of the Dynkin game value with respect to
 changes in the volatility of the underlying.
 There is plentiful research on Dynkin games in many other areas: for examples,
 \cite{Hamadene_Lepeltier_Wu_1999, Hamad_Lepeltier_2000, Hamadene_2006,Hamadene_Hdhiri_2006, HRZ_2009}
 added stochastic {\it controls} into the Dynkin games to study mixed zero-sum stochastic differential games
 of control and stopping;
 \cite{Taksar_1985, Karatzas_Wang_2001,Fukushima_Taksar_2002, Boetius_2005}
 and      \cite{Stettner_1982,Cosso_2013}  studied some Dynkin games through the associated
  singular control problems and impulse control problems respectively;
  \cite{Yasuda_1985, RSV_2001, Touzi_Vieille_2002, Laraki_Solan_2005}  considered
the Dynkin games in which the players can choose randomized stopping times;
 and   \cite{Bensoussan_Friedman_1977, Ohtsubo_1987, Nagai_1987, Cattiaux_Lepeltier_1990,
Hamadene_Zhang_2010, Hamadene_Hassani_2014a, Hamadene_Hassani_2014b} analyzed non-zero sum Dynkin games.

  However, there are only a few works on  Dynkin games under model uncertainty:
  Hamadene and Hdhiri \cite{Hamadene_Hdhiri_2006} and Yin \cite{Yin_JL_2012}    studied the Dynkin games over
      a set of equivalent probabilities,  
      which represents  drift uncertainty (or {\it Knightian uncertainty}).
 When the probability set contains mutually singular probabilities (or equivalently, both drift and volatility of the underlying can
 be ``manipulated" against Player 1),  Dolinsky \cite{Dolinsky_2014} derived dual expressions for the super–replication prices
 of game options in the discrete time,  and Matoussi et al. \cite{MPP_2014} related the Dynkin games under $G-$expectations
 (introduced by Peng \cite{Peng_G_2007b}) to second-order doubly reflected BSDEs.

 In this paper we substantially benefit from the martingale techniques developed for robust optimal stopping problems by
 \cite{Kara_Zam_2008, OS_CRM} (which analyzed the problem when $\mathcal{P}$ is dominated), \cite{ETZ_2012} ($\mathcal{P}$ is non-dominated but the Nature and the stopper cooperate) and \cite{NZ_2012, ROSVU} (in which $\mathcal{P}$ is non-dominated and the Nature and stopper are adversaries.) Especially the results of \cite{RDOSRT} are crucial for determining a saddle point. (The latter results also recently proved to be useful for defining the viscosity solutions of fully non-linear degenerate path dependent PDEs in \cite{EkrenZhang1016}).

      The rest of the paper is organized as follows:
   In Section~\ref{subsec:preliminary}, we will introduce some notation and preliminary results
  such as the regular conditional probability distribution.
  In Section \ref{sec:WSUP},  we set-up the stage for our main
result by imposing some assumptions on the reward process and the classes of mutually singular probabilities. Then
Section  \ref{sec:DPP} derives properties of Player 1's upper value processes and approximating value processes
 such as path regularity and dynamic  programming principles.
They play essential roles in  deriving our main result on the robust Dynkin games
stated in Section  \ref{sec:RDG}. In Section \ref{sec:example},
we give an example of path-dependent SDEs with controls that satisfies all our  assumptions.
In Section  \ref{sec:optim_tri},   we discuss the optimal triplet for  Player 1's  value  under additional conditions.
  Section \ref{sec:proofs} contains proofs of our results while
the demonstration of some auxiliary statements with starred labels (in the corresponding equation numbers)
in these proofs  are deferred to the Appendix.
We also include in the appendix a technical lemma necessary for the proof of Theorem \ref{thm_RDG}.

\subsection{Notation and Preliminaries} \label{subsec:preliminary}


 Throughout this paper, we fix $d  \n  \in   \n  \hN$.
 Let 
  $  \cS^{>0}_d $ stand for  all $\hR^{d \times d}-$valued positively definite matrices and
    denote by $ \sB(\cS^{>0}_d)$ 
   the Borel $\si-$field of $\cS^{>0}_d$ 
   under the relative Euclidean topology. We also fix
  a time horizon $T \n \in \n  (0,\infty)  $  and  let $t  \n \in \n  [0 , T ] $.

 We set $\O^t   \n := \n   \big\{\o    \n  \in   \n   \hC \big([t,T]; \hR^d \big)  \n  :   \o(t)   \n  =  \n   0 \big\}$
 as the   canonical space    over    period    $[t,T]$
         and denote its null path by
     $ \bz^t  \n := \n  \{\o (s)   \n = \n  0 , \fa s  \n \in \n  [t,T] \}$.
     For any $ s  \n \in \n  [ t , T ]  $,
     $ \|\o\|_{t,s}  \n := \n  \underset{r \in [t,s]}{\sup} |\o(r)|  $, $ \fa \o  \n \in \n  \O^t $
     defines  a semi-norm   on   $\O^t   $.
     In particular, $\|\cd\|_{t,T}$ is the {\it uniform} norm on $\O^t$.

     The canonical process  $ B^t $  of  $\O^t$
   is a  $d-$dimensional standard Brownian motion  under  the   Wiener measure $\hP^t_0$ of $\big(\O^t,  \cF^t_T\big)$.
    Let      $\bF^t \n  =   \n    \{ \cF^t_s     \}_{s \in [t,T]}$,
    with $\cF^t_s   \n  :=   \n  \si \big(B^t_r; \, r   \n  \in  \n   [t,s]\big)$,
     be the  natural filtration of $ B^t $ and denote its $\hP^t_0-$augmentation by
    $\ol{\bF}^t  \n = \n   \big\{ \ol{\cF}^t_{\n s}    \big\}_{s \in [t,T]}$,
    where $\ol{\cF}^t_{\n s}  \n : = \n  \si \big(\cF^t_s  \n \cup \n  \ol{\sN}^t \big)$ and
     $ \ol{\sN}^t  \n := \n  \big\{ \cN  \n \subset \n  \O^t \n : \cN  \n \subset \n  A \hb{ for some  }
    A  \n \in \n  \cF^t_T  \hb{ with } \hP^t_0 (A )  \n = \n 0  \big\} $.
    The expectation on $\big(\O^t, \ol{\cF}^t_T, \hP^t_0\big)$ will be simply denoted by $\hE_t$.
    Also, we
    let $\sP^t  $  be  the   $\bF^t-$progressively  measurable sigma$-$field of $ [t,T]  \n \times \n  \O^t$
    and let $\cT^t$ \big(resp. $\ol{\cT}^t$\big)  collect all $\bF^t $\big(resp. $\ol{\bF}^t $\big)$-$stopping times.

    Given $s  \n \in \n  [t, T]$,
    we set $\cT^t_s  \n := \n  \{\tau  \n \in \n  \cT^t \n : \tau (\o)  \n \ge \n  s , \fa \o  \n \in \n  \O^t \}$,
    $\ol{\cT}^t_s  \n := \n  \{\tau  \n \in \n  \ol{\cT}^t \n : \tau (\o)  \n \ge \n  s , \fa \o  \n \in \n  \O^t \}$
    and define   the  {\it truncation}  mapping $\Pi^t_s$  from
 $ \O^t $ to $  \O^s $   by
 $ \big(\Pi^t_s (\o)\big)(r) \n := \n \o (r) \n - \n  \o (s) $,
 $  \fa (r,\o)  \n \in \n  [s, T] \n \times \n  \O^t $.
  By Lemma A.1 of \cite{ROSVU},
  \bea  \label{eq:xax001}
     \tau  (\Pi^t_s)   \n  \in \n  \cT^t_s  , \q   \fa \tau  \n  \in \n  \cT^s .
   \eea
  For any $\d  \n > \n  0     $ and $ \o  \n \in \n  \O^t$,
  \bea   \label{eq:bb237}
 O^s_{\d} (\o)  := \big\{\o' \in \O^t:  \| \o'  -   \o  \|_{t,s} < \d \big\}
 \hb{ is an $  \cF^t_s-$measurable open set of }\O^t ,
    \eea
    and $\ol{O}^s_{\d} (\o) \n := \n  \big\{\o'  \n \in \n  \O^t:  \| \o'   \n - \n    \o  \|_{t,s}  \n \le \n  \d \big\}$
   is an $  \cF^t_s-$measurable closed set of $\O^t$
    \big(see e.g.   (2.1) of \cite{ROSVU}\big).
  In particular, we will simply denote $O^T_{\d} (\o)$ and $\ol{O}^T_{\d} (\o)$ by
 $O_{\d} (\o)$ and $\ol{O}_{\d} (\o)$ respectively.

 For any $n  \n \in \n  \hN$ and $s  \n \in \n  [t,T]$, let
  $\cT^t(n)$ denote all $\bF^t-$stopping times taking values in
  $ \{ t^n_i  \}^{2^n}_{i=0}$ with
  \bea \label{eq:uxu185}
  t^n_i \df t  \n \vee \n  (i 2^{-n}   T ) , \q i= 0, \cds , 2^n ,
  \eea
  and set $ \cT^t_s (n)  \n := \n   \{\tau  \n \in \n  \cT^t (n) \n : \tau (\o)  \n \ge \n  s , \fa \o  \n \in \n  \O^t \}$.
  In particular, 
  we literally set $\cT^t (\infty)  \n : = \n  \cT^t $   and $\cT^t_s (\infty)  \n : = \n  \cT^t_s $.

    Let $\fP_t  $ collect all probabilities   on   $\big(\O^t,  \cF^t_T  \big)   $.
    For any $\hP \n \in \n \fP_t  $, we   consider the following     spaces   about $\hP$:

   \no 1)  For any sub sigma-field  $\cG$  of $\cF^t_T$,  let
     $L^1 (\cG , \hP) $ be  the space of all  real-valued,
$  \cG-$measurable random variables $\xi$ with $\|\xi\|_{L^1(\cG ,\hP)}
   :=   \hE_\hP  \big[ |\xi|  \big]  < \infty$.

 \no 2)        Let   $\hS      (\bF^t,   \hP)$
 be the space of all  real$-$valued, ${\bF^t}-$adapted   processes $\{X_s\}_{s \in [t,T]}$
       with all   continuous paths   and satisfying
       $\hE_\hP  [ X_*  ] \n  < \n \infty $, where $X_*  \n := \n \|X\|_{t,T} \= \underset{s \in [t,T]}{\sup}|X_s| $.

      We will drop the superscript $t$ from the above notations if it is $0$.
  For example,   $  (\O,\cF)    \n =  \n  (\O^0,\cF^0) $.

  We say that a process $X$ is   bounded by some $C  > 0$
  if $ |X_t(\o)| \le C   $ for any $ (t,\o) \in [0,T] \times \O $. Also,
  a real-valued  process $X$ is said to be   uniformly   continuous
  on $  [0,T] \times \O$   with respect to some modulus of continuity function $\rho   $ if
       \bea  \label{eq:aa211}
        | X_{t_1}( \o_1) \n - \n   X_{t_2}( \o_2) |
   \n \le \n  \rho  \Big(  \bd_\infty \big((t_1, \o_1), (t_2, \o_2)\big)  \Big)  , \q
   \fa   (t_1, \o_1), (t_2, \o_2) \in [0, T] \times \O   ,
   \eea
 where $\bd_\infty \big((t_1, \o_1), (t_2, \o_2)\big)
 \n := \n  |t_1 \n - \n t_2|   \n + \n  \|\o_1 (\cd  \n \land \n  t_1)  \n - \n  \o_2(\cd  \n \land \n  t_2) \|_{0,T} $.
 For any $t  \n \in \n  [0,T]$, taking $t_1  \n = \n t_2  \n = \n t$ in \eqref{eq:aa211} shows that
  $ \big| X_t( \o_1) \n - \n   X_t( \o_2)\big|  \n \le \n  \rho \big( \|\o_1  \n - \n  \o_2\|_{0,t} \big) $,
  $  \o_1 ,  \o_2  \n \in \n  \O$, which implies the $\cF_t-$measurability of  $X_t$. So
 \beas 
 \hb{  $X$ is   indeed an $\bF-$adapted process with all continuous paths. }
 \eeas

 Moreover, let $\fM$ denote all modulus of continuity functions $\rho$ such that
 for some $ C  \n > \n  0    $ and
 $ 0  \n < \n  p_1  \n \le \n    p_2  $,
 \bea \label{eq:gb017}
 \rho (x )  \n \le \n    C  (x^{p_1}  \n \vee \n  x^{p_2} )      ,
   \q \fa x  \n \in \n    [0, \infty) .
 \eea
 In this paper,    we will  use the convention $ \inf \es := \infty $.

 \subsection{Shifted Processes and Regular Conditional Probability Distributions}

 In this subsection, we fix $ 0 \n \le \n  t   \n \le \n   s   \n \le \n   T$.
 The concatenation $\o  \n \otimes_s \n   \wt{\o}$ of  an  $\o  \n \in \n  \O^t$
 and an $ \wt{\o}  \n \in \n  \O^s$ at time $s$:
 \beas  
 \big(\o \otimes_s  \wt{\o}\big)(r)  :=   \o(r) \, \b1_{\{r \in [t,s)\}}
  + \big(\o(s) + \wt{\o}(r) \big) \, \b1_{\{r \in [s,T]\}} , \q \fa  r \in [t,T]
 \eeas
 defines another path in $\O^t$.
 Set $\o  \n \otimes_s \n  \es   \n = \n \es $ and $  \o  \n \otimes_s \n  \wt{A}  \n := \n
 \big\{ \o   \n \otimes_s \n  \wt{\o} \n : \wt{\o}  \n \in \n   \wt{A} \big\}$
 for  any non-empty  subset   $\wt{A} $ of $  \O^s $.


\begin{lemm}  \label{lem_element}
  If  $ A \in \cF^t_s$, then  $  \o \otimes_s \O^s   \subset A  $   for any $\o \in A  $.

  \end{lemm}

    For any    $ \cF^t_s-$measurable random variable $\eta$,
  since $   \{ \o'  \n \in \n  \O^t  \n : \eta (\o') \n
  = \n  \eta (\o) \} \n \in \n \cF^t_s$,
    Lemma \ref{lem_element} implies that
 \bea   \label{eq:bb421}
  \o  \n \otimes_s \n  \O^s \subset \{\o'  \n \in \n  \O^t  \n : \eta (\o') \n = \n  \eta (\o) \}
 \q \hb{i.e.,} \q
    \eta(  \o \otimes_s \wt{\o} )  \n = \n  \eta(\o), \q
  \fa \wt{\o}  \n \in \n  \O^s .
  \eea
  To wit, the value $ \eta(\o)  $  depends only on $\o|_{[t,s]}$.

   Let $\o  \n \in \n  \O^t$.  For  any   $A  \n \subset \n  \O^t$ we set $A^{s, \o}  \n := \n
   \{ \wt{\o}  \n \in \n  \O^s \n : \o  \n \otimes_s \n  \wt{\o}  \n \in \n  A  \} $
   as the  projection of $A$  on $\O^s $ along $\o$. In particular, $\es^{s,\o}  \n = \n  \es$.
 Given a random variable $\xi$  on $\O^t$,  define  the {\it shift} $\xi^{s,\o}$ of $\xi$ along $\o|_{[t,s]}$
 by  $ \xi^{s,\o}(\wt{\o})  \n := \n  \xi ( \o  \n \otimes_s \n   \wt{\o} ) $,  $  \fa  \wt{\o}  \n \in \n  \O^s $.
 Correspondingly, for a process $X  \n = \n  \{X_r\}_{r \in [t,T]}$ on $\O^t$,
     its {\it shifted} process $X^{s,\o}$ is
      \beas
       X^{s,\o}  (r, \wt{\o}) := (X_r)^{s,\o}(\wt{\o}) = X_r ( \o \otimes_s \wt{\o}) , 
      \q  \fa   ( r, \wt{\o} )  \in [s,T] \times \O^s   .
      \eeas

   Shifted random variables and shifted processes ``inherit" the measurability of   original ones:

 \begin{prop}  \label{prop_shift0}
   Let $ 0 \n \le \n  t   \n \le \n   s   \n \le \n   T$ and     $\o \in \O^t$.

 \no \(1\)  If a  real-valued random variable $\xi $ on $\O^t$ is   $\cF^t_r-$measurable for some $r \in [s,T]$,
   then   $\xi^{s,\o} $ is $  \cF^s_r-$measurable.

   \no \(2\)  For any $n \in \hN \cup \{\infty\}$ and $\tau  \n \in \n  \cT^t (n) $,
  if $\t(\o  \n \otimes_s \n  \O^s )  \n \subset \n  [r,T]$
  for some $r   \n \in \n  [s,T]$,  then $\t^{s,\o}  \n \in \n  \cT^s_r  (n)  $.

  \no \(3\)  Given  $\tau  \n \in \n  \cT^t$,  if $\tau (\o)  \n \le \n  s$, then
 $ \tau (\o  \n \otimes_s \n  \O^s )  \n \equiv \n  \tau (\o) $; if $\tau (\o)  \n \ge \n  s$
 \(resp.  $ \n > \n  s$\), then
 $ \tau (\o  \n \otimes_s \n  \wt{\o})  \n \ge \n  s $ \(resp.  $ \n > \n  s$\), $\fa \wt{\o}  \n \in \n  \O^s$
 and thus $\tau^{s,\o}  \n \in \n  \cT^s$.

  \no \(4\)  If a real-valued   process $  \{X_r \}_{r \in [t, T]}$
  is $\bF^t-$adapted \(resp. $\bF^t-$progressively measurable\),
 then  $  X^{s,\o} $ is $\bF^s-$adapted \(resp. $\bF^s-$progressively measurable\).

 \end{prop}

  Let  $\hP \n \in \n  \fP_t $. In light of the
   {\it regular conditional probability distributions} (see e.g. \cite{Stroock_Varadhan}),
   we can follow Section 2.2 of \cite{ROSVU} to introduce
   a family of {\it shifted} probabilities $\{\hP^{s,\o}  \}_{\o \in \O^t } \n \subset \n \fP_s $,
   under which the corresponding shifted random variables
   and shifted processes inherit the $\hP$ integrability of original ones:

   \begin{prop}  \label{prop_shift1}
    \(1\)  It holds for $\hP^t_0-$a.s.  $ \o  \n \in \n  \O^t$  that  $\big(\hP^t_0\big)^{s, \o } = \hP^s_0$.

  \no   \(2\)   If   $\xi \n \in \n  L^1 \big( \cF^t_T,\hP \big)$ for some   $\hP  \n \in \n  \fP_t$,
  then it holds for $\hP -$a.s.     $\o  \n \in \n  \O^t$  that
   $\xi^{s,\o}  \n \in \n  L^1 \big( \cF^s_T ,  \hP^{s,\o}   \big) $ and
  \bea   \label{eq:f475}
 \hE_{\hP^{s,\o}}  \big[ \xi^{s,\o} \big]= \hE_\hP \big[\xi\big| \cF^t_s\big](\o) \in  \hR  .
    \eea
  \(3\)     If   $X  \n \in \n  \hS \big( \bF^t ,\hP \big)$   for some  $\hP  \n \in \n  \fP_t$,
  then it holds for $\hP -$a.s.     $\o  \n \in \n  \O^t$  that
     $X^{s,\o}    \in  $ $  \hS \big( \bF^s ,  \hP^{s,\o}   \big) $.

 \end{prop}

  As a consequence of \eqref{eq:f475},
  a shifted $ \hP^t_0- $null set   also  has zero measure.  

   \begin{lemm} \label{lem_null_sets}
  For any      $  \cN  \in \ol{\sN}^t        $,
  it holds  for  $\hP^t_0-$a.s.  $\o  \in \O^t  $ that $\cN^{s,\o} \in \ol{\sN}^s$.
  \end{lemm}

 This subsection was presented in   \cite{ROSVU} with more details and proofs.
In the next three sections, we will gradually provide the technical set-up and preparation for our main results
(Theorem \ref{thm_RDG} and Theorem \ref{thm_triplet}) on the robust Dynkin game.

 \section{Weak Stability under Pasting}

  \label{sec:WSUP}



   To study the robust Dynkin game, we need some regularity conditions on the payoff processes.

  \no   {\bf Standing assumptions on payoff processes $(g, L ,U)$.}

   \no  \(\textbf{A}\) $g$, $L$ and $U$ are three real-valued  processes that are
   uniformly   continuous  on $  [0,T] \times \O$
   with respect to   the same   modulus  of continuity function $\rho_0$ and satisfy
      $L_t (\o) \le U_t (\o)$,    $\fa (t,\o) \in [0,T] \times \O$.

   For any $(t,\o) \in [0,T] \times \O$ and $s, s' \in [t, T]$, we technically  define
$ R(t,s,s',\o) : = \int_t^{s \land s'} \n g_r (\o) dr  +
     \b1_{\{s \le s'\}}  L_s (\o)  + \b1_{\{ s' < s \}}  U_{s'} (\o) $. By \eqref{eq:aa211},
 \bea  \label{eq:da025}
   |R(t,s,s',\o_1) \n - \n  R(t,s,s',\o_2) |
  & \tn  \le & \tn   \int_t^{s \land s'} \n | g_r (\o_1) \n - \n  g_r (\o_2)| dr
    \n + \n   \b1_{\{s \le s'\}} | L_s (\o_1 )  \n - \n   L_s (\o_2 ) |
      \n + \n  \b1_{\{ s' < s \}} | U_{s'} (\o_1)  \n - \n  U_{s'} (\o_2) | \nonumber \\
  & \tn  \le & \tn   (1 \n + \n  s \n \land  \n  s'  \n - \n  t )
  \rho_0 \big( \|\o_1  \n - \n  \o_2 \|_{0, s \land s'}  \big) ,
  \q \fa \o_1, \o_2 \in \O .
 \eea

 Let the robust Dynkin game start  from time $t  \n \in \n  [0,T]$ when the history has
 been evolving along path $\o|_{[0,t]}$ for some $\o  \n \in \n  \O$.
 Player 1 and   2   make their own choices on the exiting time of the game.
 If Player 1 selects $\tau \in \cT^t $ and Player 2 selects $ \ga  \in \cT^t   $,
 the game ceases at  $\tau  \n \land \n  \ga$.
 Then Player 1 will receive from  her opponent
 an accumulated  reward $\int_t^{\tau \land \ga} \n g^{t,\o}_s ds $
 and a terminal    payoff $L^{t,\o}_\tau$ (resp. $U^{t,\o}_\ga$)
 if $\tau  \n \le \n  \ga$ (resp. $\ga  \n < \n  \tau$).
 Here negative $\int_t^{\tau \land \ga} \n g^{t,\o}_s ds $, $L^{t,\o}_\tau$ or   $U^{t,\o}_\ga$
 means a payment from Player 1 to Player 2.
 So Player 1's total wealth    at time $\tau  \n \land \n  \ga$ is
  \beas
   R^{t,\o} (\tau,\ga) := \int_t^{\tau \land \ga} \n g^{t,\o}_s ds +
     \b1_{\{\tau \le \ga\}}  L^{t,\o}_\tau + \b1_{\{ \ga < \tau \}}  U^{t,\o}_\ga
   = \int_t^{\tau \land \ga} \n g^{t,\o}_s ds   +   \b1_{\{\tau \le \ga\}}  L^{t,\o}_{\tau \land \ga}
   + \b1_{\{ \ga <  \tau \}}  U^{t,\o}_{\tau \land \ga}    .
 \eeas
 Since    Proposition  \ref{prop_shift0} (4)    shows  that
 $g^{t,\o}$, $L^{t,\o}$ and  $U^{t,\o}$ are $\bF^t-$adapted processes with all continuous paths,
 \bea \label{eq:ef031}
     R^{t,\o} (\tau,\ga)  \n \in \n \cF^t_{\tau \land \ga} , \q  \fa  \tau, \ga  \n \in \n  \cT^t  .
     \eea
 Also,  it is clear that
 \bea  \label{eq:da021}
    \big( R^{t,\o} (\tau,\ga)\big) (\wt{\o}) =
 R \big(t, \tau(\wt{\o}), \ga(\wt{\o}), \o \otimes_t \wt{\o} \big)  , \q \fa \wt{\o} \in \O^t .
 \eea

 Next, we define   $ \Psi_t    := ( - L_t ) \vee   U_t \vee 0     $, $t \in [0,T]  $.      By   \eqref{eq:aa211},
 one can deduce that
  \bea
      \big|  \Psi_t(\o_1 ) -  \Psi_t(\o_2) \big|      \le
     \rho_0 \big( \|\o_1 \n - \n \o_2 \|_{0,t} \big)  ,
     \q  \fa t \in [0,T] , ~   \fa \o_1 , \o_2 \in \O     ;     \label{eq:ab025}
   \eea
 (For the reader's convenience we provided a proof in Section~\ref{subsect:proof_ROSVU}.)

It is clear that
  \bea  \label{eq:ab015}
 \big|  R^{t,\o} (\tau,\ga)  \big| \le  \int_t^{\tau \land \ga} |g^{t,\o}_s| ds  +   \Psi^{t,\o}_{\tau \land \ga} ,
 \q \fa (t,\o) \in [0,T] \times \O , ~ \fa  \tau, \ga  \n \in \n  \cT^t .
 \eea

     The following  result shows that the  integrability of   shifted payoff processes
is independent of the given path  history.

   \begin{lemm} \label{lem_Phi_integr}
   Assume  \(A\).
   For any   $t \n \in \n  [0,T]$ and  $\hP  \n \in \n  \fP_t $,
      if   $   \Psi^{t,\o}   \n \in \n   \hS  (\bF^t,  \hP )$ and
   $\hE_\hP    \int_t^T \n  |g^{t,\o}_s| ds   \n < \n  \infty$
  for some $\o  \n \in \n  \O $, then   $    \Psi^{t,\o'}    \n \in \n   \hS   (\bF^t,  \hP )$ and
   $\hE_\hP   \int_t^T \n  |g^{t,\o'}_s| ds   \n < \n  \infty$
    for all $\o'  \n \in \n  \O $.

  \end{lemm}

 We will   concentrate on those probabilities $\hP$ in $\fP_t$
 under which   shifted payoff processes are  integrable:

\begin{assum}  \label{assum_fP_Phi}
    For any  $ t \n \in   \n  [0,T]  $,    $ \wh{\fP} _t   \n  :=   \n
   \Big\{\hP  \n   \in  \n   \fP_t   \n  :  \Psi^{t, \bz}   \n   \in  \n    \hS   (\bF^t,\hP)   \hb{ and }
   \hE_\hP    \int_t^T \n  |g^{t,\bz}_s| ds    \n < \n  \infty   \Big\}$ is not empty.

  \end{assum}

  \begin{rem} \label{rem_fP_LU}
  \(1\)    If  $\Psi \in \hS (\bF,\hP_0)$
  and $\hE_{\hP_0}     \int_0^T \n  |g_s| ds   \n < \n  \infty$,
  then   $\hP^t_0 \in \wh{\fP}_t  $ for any $ t \n \in   \n  [0,T]   $.
  \(2\)  As we will show in Proposition \ref{prop_P1P2P3_Ass}, when the   modulus  of continuity $\rho_0$
  in  \(A\) has polynomial growth,
  the laws  of  solutions to the controlled SDEs \eqref{FSDE1} over period $[t,T]$ belong to $\wh{\fP}_t $.
  \end{rem}

 Under  (A) and Assumption \ref{assum_fP_Phi}, one can deduce from Lemma \ref{lem_Phi_integr}   that
   for any   $ t \n \in   \n  [0,T] $ and   $\hP  \n \in \n  \wh{\fP}_t$,
   \bea  \label{eq:xxx111}
   \Psi^{t,\o} \n \in \n  \hS  \big( \bF^t ,  \hP    \big)
   \q \hb{and} \q
   \hE_\hP   \int_t^T \n  |g^{t,\o}_s| ds     < \n  \infty , \q \fa \o \in \O       .
   \eea

 Next, we  need the    probability class  to be  adapted and weakly stable under pasting in the following sense: \\

    \no   {\bf Standing assumptions on the probability class.}

    \no \(\textbf{P1}\) For any $t \in [0,T]$, we   consider a family $\{\cP(t,\o)  \}_{ \o \in \O}$
 of  subsets of $\wh{\fP}_t$ such  that
  \bea \label{eq:uxu111}
    \cP  (t,\o_1) \n = \n  \cP  (t,\o_2) \; \hb{ if } \;  \o_1 |_{[0,t]}  \n   = \n  \o_2 |_{[0,t]}  .
    \eea

      Assume  further   that the probability class
 $\{\cP(t,\o)\}_{(t,\o) \in [0,T] \times \O}$ satisfy the following two conditions
 for some   modulus  of continuity function $\wh{\rho}_0 $:
 for any  $ 0 \le  t <s \le T     $,
 $ \o  \in \O$ and $\hP  \n \in \n  \cP (t,  \o )$:



   \no \(\textbf{P2}\)  
   There exists an   extension $(\O^t,\cF',\hP')$ of $(\O^t,\cF^t_T,\hP)$
   \big(i.e. $\cF^t_T \n \subset \n  \cF' $ and $\hP'|_{\cF^t_T}  \n = \n  \hP$\big) and
   $\O' \in \cF'$ with $\hP' (\O') = 1$ such that
    $\hP^{s,   \wt{\o}} $ belongs to $  \cP (s,  \o  \otimes_t \wt{\o} ) $ for any $\wt{\o} \in \O'$.

  \no \(\textbf{P3}\) \({\it weak stability under pasting}\)
  For any  $ \d \n \in \n  \hQ_+   $ and $\l \n  \in \n \hN$,
 let $\{\cA_j\}^\l_{j=0}$ be a $\cF^t_s-$partition of $\O^t$ such that for $j \n = \n 1,\cds \n , \l$,
  $\cA_j  \n \subset \n  O^s_{\d_j} (\wt{\o}_j)$ for some
  $\d_j  \n \in \n  \big((0,\d]  \n \cap \n  \hQ\big)  \n \cup \n  \{\d\}$ and $\wt{\o}_j  \n \in \n  \O^t $.
  Then for any
    $   \hP_j   \n   \in   \n    \cP(s, \o   \n   \otimes_t   \n     \wt{\o}_j)$,
       $j  \n = \n 1,\cds  \n ,\l$,
 there exists a $\wh{\hP}   \n \in \n  \cP(t,\o) $  such that

   \no (\,i) $\wh{\hP}  (A \cap \cA_0 ) \n =  \n  \hP   (A \cap \cA_0 ) $, $ \fa A \in \cF^t_T$;

   \no  (ii) For any $j \n = \n 1,\cds  \n ,\l $ and $A \in \cF^t_s$,
   $\wh{\hP}  (A \cap \cA_j ) =  \hP   (A \cap \cA_j ) $;

   \no  (iii) For any $n \n \in \n  \hN  \n \cup \n  \{\infty\} $  and $  \wp  \n \in \n  \cT^s$,
  there exist  $\wp^n_j  \n \in \n  \cT^t_s$,   $j \n = \n 1,\cds  \n ,\l $  such that for any $A \n \in \n  \cF^t_s$
  and $\tau \in \cT^t_s (n)$
  \bea    \label{eq:xxx617}
  \sum^\l_{j=1}   \hE_{\wh{\hP} }   \big[ \b1_{A \cap \cA_j} R^{t,\o}   (\tau, \wp^n_j  )  \big]
  \n  \le \n  \sum^\l_{j=1}  \hE_{ \hP  } \bigg[ \b1_{\{\wt{\o} \in A \cap  \cA_j\}}
  \bigg(   \underset{\vs \in \cT^s (n) }{\sup}
   \hE_{\hP_j}   \big[ R^{s,\o \otimes_t \wt{\o}}  (\vs, \wp)  \big] \n + \n \int_t^s \n g^{t,\o}_r (\wt{\o}) dr
    \bigg) \bigg]  \n + \n  \wh{\rho}_0 (\d)   .
   \eea

 \begin{rem} \label{rem_P3}
 \(1\) By \eqref{eq:uxu111}, one can regard $\cP(t,\o)$ as a path-dependent subset of  $\fP_t$.
  In particular, $\cP  \n := \n  \cP (0,\bz) \n = \n  \cP (0,\o)$, $\fa \o \n \in \n  \O $.

  \no \(2\)  Both sides of  \eqref{eq:xxx617} are finite as we will show in Section \ref{sec:proofs}.
 In particular,  the expectations on the right-hand-side
 are well-defined since    the mapping $ \wt{\o} \to \n  \underset{\vs \in \cT^s (n) }{\sup}
   \hE_{\wt{\hP}}   \big[  R^{s,\o \otimes_t \wt{\o}}  (\vs, \wp)     \big]$ is continuous under   norm
   $\|~\|_{t,T}$ for any
   $n \n \in \n  \hN    \cup    \{\infty\} $,
   $\wt{\hP} \ins \wh{\fP}_s $   and   $ \wp \n \in \n  \cT^s $.

  \no \(3\) Analogous to \(P2\) assumed in \cite{ROSVU},
 the condition \(P3\)   can be regarded as a weak form of stability under pasting
 since it  is implied by   the  ``stability under finite pasting"
 \big(see e.g. \(4.18\) of \cite{STZ_2011b}\big):  for any  $0  \n \le \n  t  \n < \n  s  \n \le \n  T$,
 $ \o   \n \in \n  \O$, $\hP  \n \in \n  \cP (t,  \o )$,
    $ \d  \n \in \n  \hQ_+   $ and $\l  \n \in \n  \hN$,
    let $\{\cA_j\}^\l_{j=0}$ be a $\cF^t_s-$partition of $\O^t$ such that for $j=1,\cds \n , \l$,
  $\cA_j \subset O^s_{\d_j} (\wt{\o}_j)$ for some
  $\d_j  \n \in \n  \big((0,\d]  \n \cap \n  \hQ\big)  \n \cup \n  \{\d\}$ and
  $\wt{\o}_j \in \O^t $. Then for any
    $   \hP_j   \n   \in   \n    \cP(s, \o    \otimes_t      \wt{\o}_j)$,
       $j  \n = \n 1,\cds  \n ,\l$,
   the probability defined by
 \bea   \label{eq:xxx131c}
    \wh{\hP} (A)   \n =  \n   \hP ( A \cap    \cA_0  \big)
     +   \sum^\l_{j=1} \hE_\hP   \n  \left[   \b1_{\{\wt{\o} \in \cA_j\}}  \hP_j \big( A^{s,\wt{\o}} \big) \right]
          , \q \fa  A \in \cF^t_T
   \eea
 is in $\cP(t,\o) $.

\end{rem}

 As pointed out in Remark 3.6 of \cite{Nutz_2012a} (see also Remark 3.4 of \cite{ROSVU}),
 \eqref{eq:xxx131c} is not suitable for
 the example of path-dependent SDEs with controls
 \(see Section \ref{sec:example}\). Thus we assume the weak pasting condition \(P3\), which turns out
 to be sufficient for our approximation scheme in proving the main results.

 \section{The Dynamic Programming Principle}
  \label{sec:DPP}

    Consider the robust Dynkin game  with payoff processes $(g,L,U)$
  and over the   probability class $\{\cP(t,\o)\}_{(t,\o) \in [0,T] \times \O}$
   as described in Section \ref{sec:WSUP}.
  If Player 1  conservatively thinks that {\it Nature} is also against her, then
  for any $(t,\o) \in [0,T] \times \O$,
       \beas
     \ul{V}_t(\o) :=   \underset{\t  \in \cT^t }{\sup} \,   \underset{\ga  \in \cT^t }{\inf} \,
   \underset{\hP \in \cP(t,\o) }{\inf} \, \hE_\hP  \big[  R^{t,\o} (\t,\ga)    \big]
     \q \hb{and} \q
  \ol{V}_t (\o)  :=    \underset{\hP \in \cP(t,\o) }{\inf} \,  \underset{\ga  \in \cT^t }{\inf} \,
  \underset{\t  \in \cT^t }{\sup} \,   \hE_\hP  \big[  R^{t,\o} (\t,\ga)    \big]
    \eeas
  define   the {\it lower }  value and {\it upper } value of
  Player 1 at time $t$ given the historical path $\o|_{[0,t]}$.

     As  we will see in Theorem \ref{thm_RDG} that  $\ul{V} $   coincides with    $  \ol{V}     $
  as Player 1's value process $V$, whose sum with $\int_0^\cd g_s ds $
  is an $\ul{\sE}-$submartingale up to  the first time $\tau_*$ when $V$ meets $L$.
     For this purpose,   we  derive in this section some basic properties of   $\ol{V} $
     and its approximating values including  dynamic programming principles.
  Let  \(A\),    \(P1\)$-$\(P3\) and  Assumption \ref{assum_fP_Phi}   hold throughout the section.

  For any $(t,\o) \in [0,T] \times \O$, following  \cite{NZ_2012}'s idea, we technically    define approximating value processes of $\ol{V}$ by
 \bea   \label{eq:da371}
  V^n_t (\o)  :  =   \underset{\hP \in \cP (t,\o)  }{\inf} \,  \underset{\ga  \in \cT^t }{\inf} \,
  \underset{\t  \in \cT^t(n) }{\sup} \,   \hE_\hP  \big[  R^{t,\o} (\t,\ga)    \big]
  \le     \underset{\hP \in \cP (t,\o)  }{\inf} \,  \underset{\ga  \in \cT^t }{\inf} \,
  \underset{\t  \in \cT^t  }{\sup} \,   \hE_\hP  \big[  R^{t,\o} (\t,\ga)    \big] = \ol{V}_t (\o) , \q  \fa n \in \hN ,
  \eea
  and  set in particular $ V^\infty_t (\o) :=  \ol{V}_t (\o)  $.

   Let $n \in \hN \cup \{\infty\}$. It is clear that
 \bea  \label{eq:ek011}
  V^n  (T, \o  )
 =  \underset{\hP \in \cP (T,\o  )  }{\inf} \,  \underset{\ga  \in \cT^T }{\inf} \,
  \underset{\t  \in \cT^T(n) }{\sup} \,   \hE_\hP  \big[  R^{T,\o  } (\t,\ga)    \big]
 = \underset{\hP \in \cP (T,\o  )  }{\inf} \,  \hE_\hP  \big[ R^{T,\o  } (T,T) \big] =  L_T (\o) , \q \fa \o \in \O .
 \eea
 And we can show that
 \bea \label{eq:uxu170}
 - \Psi_t (\o) \ls L_t (\o) \ls   V^n_t (\o)  \ls U_t (\o) \ls \Psi_t (\o) , \q \fa (t,\o) \in [0,T] \times \O .
 \eea
For the reader's convenience we provide a proof in Section~\ref{subsect:proof_ROSVU}.

We need the following assumption on $V^n$'s to discuss the dynamic programming principles they satisfy.

     \begin{assum} \label{assum_V_conti}
     There exists a   modulus  of continuity function $\rho_1  \ge \rho_0    $
     such that  for any $n \in \hN \cup \{\infty\}$
      \bea \label{eq:aa213}
         \big|   V^n_t (\o_1) - V^n_t (\o_2) \big| \le   \rho_1 \big( \|\o_1  - \o_2 \|_{0,t} \big)  ,
     \q    \fa t \in [0,T] , ~ \fa     \o_1 , \o_2 \in \O .
   \eea
    \end{assum}

    \begin{rem}  \label{rem_V_conti}
     If $\cP(t,\o)$ does not depend on $\o$ for all $t \in [0,T]$, then
     Assumption \ref{assum_V_conti} holds automatically.
    \end{rem}

 \begin{rem}     \label{rem_V_adapted}
    Assumption \ref{assum_V_conti} implies that   $V^n$ is   $\bF-$adapted for any   $n \in \hN \cup \{\infty\}$.
 \end{rem}

   We first present the sub-solution side of dynamic programming principle
  for $V^n$\,'s:

 \begin{prop}   \label{prop_DPP}
 For any $n \in \hN \cup \{\infty\}$, $0 \le  t \le s \le T $ and $  \o \in \O $,
 \bea   \label{eq:bb013}
 V^n_t (\o) \le \underset{\hP \in \cP(t,\o)}{\inf} \, \underset{\ga \in \cT^t  }{\inf} \,
  \underset{\t \in \cT^t (n) }{\sup} \,  \hE_\hP \bigg[ \b1_{\{\tau \land \ga < s\}}   R^{t,\o}  (\tau, \ga)
   + \b1_{\{\tau \land \ga \ge s\}}  \Big( \big( V^n_s \big)^{t,\o}  +
   \int_t^s g^{t,\o}_r dr  \Big)      \bigg]   .
 \eea
 \end{prop}

Conversely, we only need to show  the super-solution side of dynamic programming principle for $V^\infty = \ol{V} $.

 \begin{prop}   \label{prop_DPP2}
 For any   $0 \le  t \le s \le T $ and $  \o \in \O $,
 \beas 
 \ol{V}_t (\o) \ge \underset{\hP \in \cP(t,\o)}{\inf} \, \underset{\ga \in \cT^t  }{\inf} \,
  \underset{\t \in \cT^t   }{\sup} \,  \hE_\hP \bigg[ \b1_{\{\tau \land \ga < s\}}   R^{t,\o}  (\tau, \ga)
   + \b1_{\{\tau \land \ga \ge s\}}  \Big(   \ol{V}^{t,\o}_s   +
   \int_t^s g^{t,\o}_r dr  \Big)      \bigg]   .
 \eeas
 \end{prop}

 As a consequence of Propositions \ref{prop_DPP} and \ref{prop_DPP2},
 the upper value process $\ol{V}$ of Player 1 satisfies a true dynamic programming principle.

   We rely on another condition to further show the convergence of $V^n$ to $\ol{V}$ and their path regularities
 in the next two propositions.

     \begin{assum}  \label{assum_V_conti_2}
 For any $\a > 0$,
       there exists a modulus  of continuity function   $ \rho_\a $ such that  for any $t \in [0,T)$
   \bea    \label{eq:aa213b}
     \underset{\o \in O^t_\a (\bz) }{\sup} \;      \underset{\hP \in \cP(t,\o)}{\sup} \;
       \underset{\z \in \cT^t}{\sup}  \hE_\hP
 \bigg[  \rho_1 \Big(    \d +      \underset{r \in [\z ,  (\z + \d) \land T ]    }{\sup}
  \big|  B^{t}_r - B^{t}_\z  \big|  \Big) \bigg] \le   \rho_\a  (\d) , \q \fa \d \in  ( 0 , T  ] .
   \eea

   \end{assum}

    \begin{prop} \label{prop_grid}
 Let $n \in \hN$, $t \in [0,T]$ and $\a > 0$. It holds for any $ \o \in O^t_\a (\bz)$ that
 \bea    \label{eq:da365}
 \ol{V}_t (\o) \le  V^n_t (\o) + \rho_\a (2^{-n})    + 2^{-n} \big( |g_t(\o)| + \rho_\a (T \n - \n t) \big) .
 \eea
 \end{prop}

\begin{prop}  \label{prop_conti_V}
 \(1\)  For any $n \in \hN \cup \{\infty\}$, all paths of process $V^n$ are both   left-upper-semicontinuous and right-lower-semicontinuous. In particular, the process $\ol{V}$ has all continuous paths.

 \no  \(2\) For any $(t,\o) \in [0,T] \times \O$ and $\hP \in \cP(t,\o)$,  $\ol{V}^{t,\o} \in \hS  ( \bF^t,   \hP)$.

\end{prop}

 \section{Main Result}

 \label{sec:RDG}

 In this section, we state   our first main result  on robust Dynkin games.
 Let  \(A\),         \(P1\)$-$\(P3\)  and
 Assumptions \ref{assum_fP_Phi}, \ref{assum_V_conti}, \ref{assum_V_conti_2}  hold throughout the section.

       Given $t \n \in \n  [0,T]$,    set
    $   \sL_t  \n :=  \n    \{\hb{random variable $\xi$ on $\O$} \n :
     \xi^{t,\o}  \n \in \n  L^1(\cF^t_T, \hP),~ \fa \o  \n \in \n  \O ,\;
   \fa \hP  \n \in \n  \cP(t,\o)  \}$. Clearly, $\sL_t$ is closed under linear combination: i.e.
   for any $\xi_1, \xi_2 \in \sL_t$ and $\a_1, \a_2 \in \hR$, $\a_1 \xi_1 + \a_2 \xi_2 \in \sL_t$.
   Then we  define on $\sL_t$ a nonlinear  expectation:
    \beas
       \ul{\sE}_t [\xi] (\o) :=   \underset{\hP \in \cP(t,\o) }{\inf} \hE_\hP [\xi^{t,\o}]   ,  \q
    \fa    \o \in \O , ~ \fa   \xi \in   \sL_t  .
    \eeas
    For any $n \in \hN \cup \{\infty\}$ and $\tau   \n \in \n  \cT$,
 \bea \label{eq:uxu173}
   \hb{ both  $   V^n_\tau $ and $ \int_0^\tau g_r dr  $ belong to $ \sL_t$. }
  \eea
  (We demonstrate this claim in Section~\ref{pfsec4}.)

     Similar to the classic Dynkin game,  we will show that
   $\ol{V} $ coincides with $\ul{V} $ as the value process $V$ of Player 1 in the robust Dynkin game and that
   $V$ plus $\int_0^\cd g_s ds $ is a  submartingale  with respect to the nonlinear expectation $\ul{\sE}  $.

         \begin{thm} \label{thm_RDG}
  Let   \(A\),    \(P1\)$-$\(P3\)  and Assumptions \ref{assum_fP_Phi}, \ref{assum_V_conti}, \ref{assum_V_conti_2} hold.

  \no \(1\) For any $ (t,\o)  \n \in \n  [0,T]  \n \times \n  \O $,
  \bea \label{eq:value_equal}
  V_t(\o)  \n :=  \n  \ul{V}_t (\o)  \n = \n  \ol{V}_t (\o)
  \eea
   in   the robust Dynkin game starting from time $t$ given
  the historical  path $\o \big|_{[0,t]}$. Moreover,
    \bea \label{eq:xax015}
   V_t (\o) =  \underset{\ga \in \cT^t}{\inf} \,   \underset{\hP \in \cP(t,\o)}{\inf} \,
  \hE_\hP \big[ R^{t,\o}   \big( \tau^*_{(t,\o)}  ,   \ga \big)   \big]    , \;
  \hb{   where $\tau^*_{(t,\o)}    \n : = \n  \inf\big\{s  \n \in \n  [t,T] \n :  V^{t,\o}_s
  \n = \n    L^{t,\o}_s    \big\}   \n \in \n  \cT^t   $. }
 \eea

   \no \(2\) The $\bF-$adapted process with all continuous paths $ \U_t  :=
   V_t + \int_0^t g_r dr $, $ t \in [0,T]$ is an $\ul{\sE}-$submartingale
    up to time  $\tau_*  \n : = \n \tau^*_{(0,\bz)}
    \n = \n  \inf\big\{t  \n \in \n  [0,T] \n :  V_t
  \n = \n    L_t    \big\}    \n \in \n  \cT   $
   in sense that for any $\zeta  \n \in \n  \cT$
  \bea    \label{eq:cc761}
     \U_{\tau_*   \land \zeta   \land  t}       ( \o)
   \n  \le  \n   \ul{\sE}_t  \big[ \U_{\tau_*  \land \zeta }   \big] ( \o)  ,
   \q \fa (t,\o)  \n \in \n  [0,T]  \n \times \n  \O    .
  \eea

 \end{thm}

 \if{0}
\begin{rem} \label{rem_WCRM}

    Define a worst-case  risk measure $\vf (\xi) \n :=  \n   \underset{\hP \in \cP}{\sup} \,  \hE_\hP  [ - \xi   ] $
      for any bounded financial position $\xi$.   Theorem \ref{thm_RDG} implies that
      $   V_0 \n = \n   \underset{\t  \in \cT }{\sup} \,   \underset{\ga  \in \cT }{\inf} \,
   \underset{\hP \in \cP }{\inf} \, \hE_\hP  \big[  R (\t,\ga)    \big]
   =   \underset{\ga  \in \cT }{\inf} \,   \underset{\t  \in \cT }{\sup} \,
    \underset{\hP \in \cP }{\inf} \,   \hE_\hP  \big[  R (\t,\ga)    \big]
     =    \underset{\ga  \in \cT }{\inf} \,  \underset{\hP \in \cP }{\inf} \,
  \underset{\t  \in \cT }{\sup} \,   \hE_\hP  \big[  R (\t,\ga)    \big]  $.
  Multiplying  $-1$ to each side yields that
  \bea  \label{eq:xax017}
  \underset{\tau  \in \cT  }{\inf} \;  \underset{\ga  \in \cT  }{\sup} \; \vf    \big(  R (\tau, \ga)   \big)
  =  \underset{\ga  \in \cT  }{\sup} \; \underset{\tau  \in \cT  }{\inf} \;  \vf    \big(  R (\tau, \ga)   \big)  .
 \eea
 If $ R (\tau, \ga) $ is the payment from Player 1 to Player 2, then  \eqref{eq:xax017}
 is exactly Player 2's value in the corresponding Dynkin game under the worst-case  risk measure $\vf$.

\end{rem}
\fi

 \section{Examples: Controlled Path-dependent  SDEs}

 \label{sec:example}

In this section,  we  provide an example of the probability class $\{\cP(t,\o)\}_{(t,\o) \in [0,T] \times \O}$
in case of path-dependent stochastic differential equations with controls.

   Let $\k \n > \n 0$  
  and     let  $b  \n : [0,T]  \n \times \n  \O   \n \times \n  \hR^{d \times d} \to \hR^d $ be a
 $\sP    \n  \otimes  \n   \sB(\hR^{d \times d})\big/\sB(\hR^d)-$measurable   function such that
 \beas 
 |b(t,\o ,u) \n - \n b(t,\o',u)|  \n \le \n  \k \|\o  \n - \n \o' \|_{0,t}
 \q \hb{and} \q |b(t,\bz,u)|  \n \le \n  \k (1 \n + \n |u|)  ,
 \q \fa \o , \o'  \n \in \n  \O,  ~   (t,u)  \n \in \n  [0,T]   \n   \times  \n  \hR^{d \times d} .
 \eeas

 Fix $t \in [0,T]$. We   let   $\sU_t$ collect  all
       $ \cS^{>0}_d-$valued, $\bF^t-$progressively measurable   processes $\{\mu_s\}_{s \in [t,T]}$
       such that $|\mu_s| \le \k $, $ds \times d\hP^t_0-$a.s.
       Let $\o \n \in \n \O$,    $ b^{t,\o} (r,\wt{\o},u) \n := \n  b(r,\o \otimes_t \wt{\o},u)$,
  $  (r,\wt{\o},u)\in [t,T] \times \O^t \times \hR^{d \times d}$
  is clearly a $\sP^t \n \otimes \n \sB(\hR^{d \times d})\big/\sB(\hR^d)-$measurable  function that satisfies
     \beas
 \q |b^{t,\o}(r,\wt{\o} ,u) \n - \n b^{t,\o}(r,\wt{\o}',u)|  \n \le \n  \k \|\wt{\o}  \n - \n \wt{\o}' \|_{t,r}
 ~ \hb{and} ~ |b^{t,\o}(r, \bz^t,u)|  \n \le \n  \k \big( 1 \n + \n  \|\o\|_{0,t}  \n + \n  |u| \big) ,
 ~ \fa \wt{\o} , \wt{\o}'  \n \in \n  \O^t ,  \;   (r,u)  \n \in \n  [t,T]   \n   \times  \n  \hR^{d \times d} .
 \eeas

  Given    $\mu \in \sU_t$,
  a slight extension of Theorem V.12.1 of \cite{Rogers_Williams_2} shows that
  the following     SDE   on the probability space $\big( \O^t, \cF^t_T, \hP^t_0 \big)$:
    \bea   \label{FSDE1}
     X_s   =   \int_t^s b^{t,\o} (r,   X, \mu_r  )  dr + \int_t^s  \mu_r  \, dB^t_r ,  \q s \in [t, T]  ,
     \eea
           admits a unique solution $  X^{t,\o,\mu}$,   which is an   $\ol{\bF}^t-$adapted continuous  process
           satisfying $E_t \big[ \big( X^{t,\o,\mu}_* \big)^p \big] \< \infty$ for any $p \gs 1$
           (or see the complete ArXiv version of \cite{ROSVU} for its proof).

      Note that the SDE \eqref{FSDE1} depends on  $\o \big|_{[0,t]}$ via the generator $b^{t,\o}$.
    Without loss of generality, we   assume that all paths of   $X^{t,\o,\mu}$ are continuous and starting from $0$.
   \big(Otherwise, by setting $\cN  \n  :=  \n \{\o  \n \in \n  \O^t \n : X^{t,\o,\mu}_t(\o)
    \n \ne  \n  \bz \hb{ or the path $  X^{t,\o,\mu}_\cd (\o)$ is not continuous}\}  \n \in  \n \ol{\sN}^t $,
    one can take
        $ \wt{X}^{t,\o,\mu}_s \n := \n  \b1_{   \cN^c   }    X^{t,\o,\mu}_s   $, $s  \n \in \n  [t,T]$.
    It is  an $\ol{\bF}^t-$adapted process that  satisfies \eqref{FSDE1}
    and  whose paths are all continuous  and  starting from $0$.\big)

     Applying the Burkholder-Davis-Gundy inequality, Gronwall's inequality
    and using the Lipschitz continuity of $b$
     in $\o-$variable, one can easily derive the following estimates for $  X^{t,\o,\mu}$:   for any $p \ge 1$
    \bea
  &&      \hE_t \bigg[ \underset{r \in [t, s]}{\sup} \big|X^{t,\o,\mu}_r  \n - \n X^{t,\o'\n,\,\mu}_r \big|^p \bigg]
     \n \le  \n
     C_p    \|\o \n - \n \o'\|^p_{0,t} \, (s \n - \n t)^p , \q
     \fa \o'  \n \in \n  \O , \; \fa s \in [t,T] , \label{eq:xxx151} \\
    \;  \hb{ and } &&   \hE_t  \bigg[  \underset{r  \in [\z, ( \z + \d) \land T]}{\sup}
         \big|  X^{t,\o,\mu}_r  -  X^{t,\o,\mu}_\z \big|^p \bigg]
          \n \le \n   \vf_p (\|\o\|_{0,t}) \, \d^{\, p/2} , \q
          \hb{for any $\ol{\bF}^t-$stopping time $\z$   and   $\d  \n > \n  0$} , \qq \label{eq:xxx153}
    \eea
    where  $C_p$ is a   constant depending on $p,  \k, T$ and
     $  \vf_p \n  : \hR_+  \n  \to \n   \hR_+  $ is a continuous function   depending on $ p,   \k, T $
     \big(see the complete ArXiv version of \cite{ROSVU} for the proofs of \eqref{eq:xxx151}   and   \eqref{eq:xxx153}\big).

   For any $s \in [t,T]$,  we see from \cite{ROSVU} that
  $    \cF^t_s     \subset \cG^{X^{t,\o,\mu}}_s  \n := \n
     \Big\{A  \n \subset \n  \O^t \n :
     \big(X^{t,\o,\mu}\big)^{-1}(A)  \n \in  \n   \ol{\cF}^t_{ s}  \Big\} $,
   i.e.,
   \bea   \label{eq:xxx439}
    \big(X^{t,\o,\mu}\big)^{-1}(A) \in   \ol{\cF}^t_s , \q \fa    A \in \cF^t_s .
   \eea
   Namely,  $  X^{t,\o,\mu}$
   is $ \ol{\cF}^t_s  \big/   \cF^t_s  -$measurable as a mapping from $\O^t$ to $\O^t$.
   \if{0}
   To see this, let us pick up  an arbitrary
   $   \cE  \in   \sB(\hR^d)$.        The $\ol{\bF}^t-$adaptness of  $ X^{t,\o,\mu} $    shows that
   for any $r \in [t,s]$
   \bea   \label{eq:xx193}
      \big( X^{t,\o,\mu}\big)^{-1}\Big( \big(B^t_r\big)^{-1}(\cE)\Big)
   = \big\{\wt{\o} \in \O^t:  X^{t,\o,\mu}   (\wt{\o}) \in \big(B^t_r\big)^{-1}(\cE) \big\}
   = \big\{\wt{\o} \in \O^t:  X^{t,\o,\mu}_r (\wt{\o})   \in \cE \big\}  \in \ol{\cF}^t_s .
   \eea
   Thus $\big(B^t_r\big)^{-1}(\cE) \n  \in \n  \cG^{X^{t,\o,\mu}}_s  \n := \n
     \Big\{A  \n \subset \n  \O^t \n :
     \big(X^{t,\o,\mu}\big)^{-1}(A)  \n \in  \n   \ol{\cF}^t_{ s}  \Big\}$,
   which is  clearly       a $\si-$field  of $\O^t$.
   It then follows that     $    \cF^t_s     \subset \cG^{X^{t,\o,\mu}}_s  $,
   i.e.,
   \bea   \label{eq:xxx439}
    \big(X^{t,\o,\mu}\big)^{-1}(A) \in   \ol{\cF}^t_s , \q \fa    A \in \cF^t_s ,
   \eea
    proving the measurability of the mapping $ X^{t,\o,\mu} $.
    \fi
  Define the law of $ X^{t,\o,\mu} $ under $\hP^t_0 $ by
       \beas  
   \fp^{t,\o,\mu} (A)   := \hP^t_0   \circ    \big( X^{t,\o,\mu} \big)^{-1} (A), \q \fa A \in \cG^{X^{t,\o,\mu}}_T ,
       \eeas
  and  denote by $\hP^{t,\o,\mu}  $ the restriction of $ \fp^{t,\o,\mu} $ on    $\big(\O^t,  \cF^t_T\big)  $.

  Now,  let us set $\cP(t,\o) \n :=  \n \big\{ \hP^{t,\o,\mu} \n : \mu  \n \in \n  \sU_t   \big\}
 \n \subset \n \fP_t $.

    \begin{prop}   \label{prop_P1P2P3_Ass}
     Let $\vr_0$ be a  modulus  of continuity function such that for some $\varpi  \n \ge \n 1 $,
 $\vr_0(\d)  \n  \le  \n   \k ( 1  \n + \n  \d^\varpi ) $,    $    \fa         \d  \n > \n  0$.
  Assume  that $g$,  $L$, $U$  satisfy \(A\) with respect to   $\vr_0$
   and that   $ \int_0^T  \n  | g_t(\bz)| dt  \n < \n  \infty $.
 Then for any   $(t,\o) \n \in \n  [0,T]  \n \times \n  \O$,
  we have $\cP(t,\o)  \n \subset \n \wh{\fP}_t$.
  And the probability class  $\{\cP(t,\o)\}_{(t,\o) \in [0,T] \times \O }$
   satisfies  \(P1\)$-$\(P3\),   Assumption \ref{assum_V_conti}$-$\ref{assum_V_conti_2}.

  \end{prop}

\begin{rem}\label{rem:matou}
\(1\) When $b \n \equiv \n 0$, Proposition \ref{prop_P1P2P3_Ass} and
the result \eqref{eq:value_equal} verify  Assumption 5.7 of   \cite{MPP_2014}
\(particularly for $t  \n = \n  0$\).   Then we know from   Theorem 5.8 therein   that
in case of controlled path-dependent SDEs with null drift,
     Player 1's value   $V$ is closely related to the  solution of a
   second-order doubly reflected backward stochastic differential equation.

 \no \(2\) Similar to \cite{ROSVU},
    the reason we consider the law of $X^{t,\o,\mu}$ under $\hP^t_0$   over   $ \cG^{X^{t,\o,\mu}}_T $
    \big(the largest $\si-$field  to induce  $ \hP^t_0  $ under the mapping $       X^{t,\o,\mu}  $\big)
     rather than $\cF^t_T$ lies in the fact that the proof of   Proposition \ref{prop_P1P2P3_Ass}
     relies heavily on the   inverse mapping $W^{t,\o,\mu}$  of  $X^{t,\o,\mu}$.
     According to the proofs of Proposition 6.2 and 6.3 in \cite{ROSVU},
     since  $W^{t,\o,\mu}$   is an $\bF^t-$progessively measurable processes
     that has only $ \fp^{t,\o,\mu} -$a.s. continuous paths,
      it holds for $\fp^{t,\o,\mu}-$a.s.  $  \wt{\o}  \in \O^t $ that
       the shifted probability $ \big(\hP^{t,\o,\mu}\big)^{s,\wt{\o}} $ is  the law of the solution to the
 shifted SDE \big(and thus   $ \big(\hP^{t,\o,\mu}\big)^{s,\wt{\o}} \n \in \n \cP( s, \o \otimes_t  \wt{\o} ) $\big).
   This explains why
  our assumption \(P2\) needs an extension $(\O^t,\cF',\hP')$ of the probability space $(\O^t,\cF^t_T,\hP)$.

\end{rem}

\section{The Optimal Triplet}

 \label{sec:optim_tri}

   In this section, we identify an optimal triplet for Player 1's value in
   the robust Dynkin game under the following additional conditions
   on  the payoff processes and the probability class.

  \no {\bf (A$'$)} Let $g \equiv 0$ and let
    $L$,  $U$ be two real-valued    processes bounded by some $M_0  \n > \n  0$
  such that  they are   uniformly   continuous   on $ [0,T]  \n \times \n  \O$
  with respect to  the same    $\rho_0 \n \in \n  \fM$, that
  $L_t (\o)  \n \le \n  U_t (\o)$,    $\fa (t,\o)  \n \in \n  [0,T)  \n \times \n  \O$,
  and that   $L_T (\o)  \n = \n  U_T (\o) $, $\fa \o  \n \in \n  \O$.

 \if{0}


 When the robust Dynkin game starts at   $t  \n \in \n  [0,T]$ given a historical path $\o|_{[0,t]}$,
 $\o  \n \in \n  \O$, Player 1's
     total wealth    is
   \beas
   R^{t,\o} (\tau,\ga) :=
     \b1_{\{\tau \le \ga\}}  L^{t,\o}_\tau + \b1_{\{ \ga < \tau \}}  U^{t,\o}_\ga
   \eeas
 if she chooses $\tau  \n \in \n  \cT^t$ and her opponent chooses $\ga  \n \in \n  \cT^t$ to quit the game.

 \fi

   Also,  let  a family $\{\cP_t  \}_{ t \in [0,T]}$ of subsets $ \cP_t $ of $
   \wh{\fP}_t = \fP_t $, $t \in [0,T]$  satisfy:

    \no \(\textbf{H1}\)    $ \cP : = \cP_0      $ is a  weakly compact subset   of $\fP_0$.

   \no ({\bf H2})  For any $\rho \in \fM$,
       there exists another   $ \ol{\rho}  $ of $\fM$ such that
   \beas
  \underset{(\hP,\z) \in \cP_t  \times \cT^t}{\sup}  \hE_\hP
 \bigg[  \rho  \Big(    \d +      \underset{ r \in [\z, (\z+\d) \land T]    }{\sup}
  \big|  B^{t}_r - B^{t}_\z  \big|  \Big) \bigg]  \le   \ol{\rho}  (\d) , \q \fa t \in [0,T) , ~
   \fa \d \in  ( 0 , \infty)   .
   \eeas
 In particular, we require $\ol{\rho}_0$ to satisfy \eqref{eq:gb017}
 with some $\ol{C}>0$ and $1<\ol{p}_1 \le \ol{p}_2  $.

   \no \(\textbf{H3}\)  For any  $ 0 \le  t <s \le T     $,
 $ \o  \in \O$ and $\hP  \n \in \n  \cP_t$,
   there exists an   extension $(\O^t,\cF',\hP')$ of $(\O^t,\cF^t_T,\hP)$
   \big(i.e. $\cF^t_T \n \subset \n  \cF' $ and $ \hP' |_{\cF^t_T}  \n = \n  \hP$\big) and
   $\O' \in \cF'$ with $\hP' (\O') = 1$ such that
    $\hP^{s,   \wt{\o}} $ belongs to $  \cP_s  $ for any $\wt{\o} \in \O'$.


   \no ({\bf H4}) Moreover, let the finite stability under pasting stated in Remark \ref{rem_P3} (3) hold.
  \if{0}
   For any  $0  \n \le \n  t  \n < \n  s  \n \le \n  T$,
 $ \o   \n \in \n  \O$, $\hP  \n \in \n  \cP_t$,
    $ \d  \n \in \n  \hQ_+   $ and $\l  \n \in \n  \hN$,
    let $\{\cA_j\}^\l_{j=0}$ be a $\cF^t_s-$partition of $\O^t$ such that for $j=1,\cds \n , \l$,
  $\cA_j \subset O^s_{\d_j} (\wt{\o}_j)$ for some $\d_j  \n \in \n  \big((0,\d]  \n \cap \n  \hQ\big) \cup \{\d\}$ and $\wt{\o}_j \in \O^t $.
Then for any
    $  \{ \hP_j  \}^\l_{j=1} \n   \subset   \n    \cP_s $,
   the probability defined by
 $   \wh{\hP} (A)   \n =  \n   \hP ( A \cap    \cA_0  \big)
     +   \sum^\l_{j=1} \hE_\hP   \n  \left[   \b1_{\{\wt{\o} \in \cA_j\}}  \hP_j \big( A^{s,\wt{\o}} \big) \right] $,
          $ \fa  A \in \cF^t_T $    is in $\cP_t$.
          \fi

      The next example shows that {\it controls of weak formulation }
(i.e. $\cP$ contains all semimartingale measures under which $B$ has uniformly bounded drift and diffusion coefficients)
satisfies (H1)$-$(H4).

\begin{eg} \label{eg_weak_formulation}
   Given $\ell > 0$, let $\{\cP^\ell_t\}_{t \in [0,T]}$
 be the family of semimartingale measures considered in \cite{ETZ_2014} such that
 $\cP^\ell_t$   collects all continuous semimartingale measures   on $(\O^t,\cF^t_T)$
 whose drift and diffusion characteristics are bounded by $\ell$ and $\sqrt{2\ell}$ respectively.
 According to Lemma 2.3 therein,  $\{\cP^\ell_t\}_{t \in [0,T]}$ satisfies
 \(H1\), \(H3\) and \(H4\).
 Also,   one can  deduce from the Burkholder-Davis-Gundy inequality that  $\{\cP^\ell_t\}_{t \in [0,T]}$ satisfies
 \(H2\), see the proof of \cite[Example 3.3]{RDOSRT}  for details.
 \end{eg}

    Remark \ref{rem_P3} (3) and a revisit of   Remark \ref{rem_V_conti}'s proof
 show  that the path-independent probability class
 $\{\cP_t\}_{t \in [0,T]}$ satisfies (P1)$-$(P3) and Assumption \ref{assum_V_conti} with $\rho_1 = \rho_0$,
  while  Assumption \ref{assum_V_conti_2} is clearly implied by (H2)
  with $\rho_\a \n \equiv \n  \ol{\rho}_0$, $\fa \a  \n > \n  0$.
  So   Theorem \ref{thm_RDG} still holds
  for the robust Dynkin game over $\{\cP_t\}_{t \in [0,T]}$.
   In addition, (H1) enables us to apply the result of \cite{RDOSRT} to solve \eqref{eq:optimal_triplet}.

         \begin{thm} \label{thm_triplet}
  Under Assumptions  \(A$'$\) and     \(H1\)$-$\(H4\),
 \if{0}
   \no \(1\) For any $ (t,\o)  \n \in \n  [0,T]  \n \times \n  \O $,
  Player 1 has a value
        \beas
      V_t(\o) : =   \underset{\t  \in \cT^t }{\sup} \,   \underset{\ga  \in \cT^t }{\inf} \,
   \underset{\hP \in \cP_t }{\inf} \, \hE_\hP  \big[  R^{t,\o} (\t,\ga)    \big]
   =    \underset{\hP \in \cP_t }{\inf} \,  \underset{\ga  \in \cT^t }{\inf} \,
  \underset{\t  \in \cT^t }{\sup} \,   \hE_\hP  \big[  R^{t,\o} (\t,\ga)    \big] \in \big[ L_t (\o) ,  U_t (\o) \big]
    \eeas
   in    the robust Dynkin game starting from time $t$ given    the historical  path $\o \big|_{[0,t]}$.

   \no \(2\) Player 1's value process $V$ is an   $\bF-$adapted process with all continuous paths such that
  for any $ (t,\o)  \n \in \n  [0,T]  \n \times \n  \O $
          \bea \label{eq:xax014}
    V_t(\o)   =     \underset{(\hP, \ga) \in \cP_t \times \cT^t}{\inf} \,
  \hE_\hP \big[ R^{t,\o}   \big( \tau^*_{(t,\o)}  ,   \ga \big)   \big] ,
    \eea
   where $\tau^*_{(t,\o)}    \n : = \n  \inf\big\{s  \n \in \n  [t,T] \n :  V^{t,\o}_s
  \n = \n    L^{t,\o}_s    \big\}   \n \in \n  \cT^t   $.

   \no \(3\) Process $V$ is an $\ul{\sE}-$submartingale
    up to time  $\tau_*  \n : = \n \tau^*_{(0,\bz)}
    \n = \n  \inf\big\{t  \n \in \n  [0,T] \n :  V_t
  \n = \n    L_t    \big\}    \n \in \n  \cT   $
   in sense that for any $\zeta  \n \in \n  \cT$
  \beas
     V_{\tau_*   \land \zeta   \land  t}       ( \o)
   \n  \le  \n   \ul{\sE}_t  \big[ V_{\tau_*  \land \zeta }   \big] ( \o)  ,
   \q \fa (t,\o)  \n \in \n  [0,T]  \n \times \n  \O    .
  \eeas
 \fi
 there exists a pair $(\hP_*,\ga_*) \n \in \n \cP \n \times \n \cT$ such that
 $ V_0  \n = \n  \hE_{\hP_*} \n  \big[ R (\tau_*,\ga_*)\big] $.

 \end{thm}

 \if{0}

With Theorem \ref{thm_triplet}, we can continue our discussion in Remark \ref{rem_WCRM}.

\begin{rem}

    In light of Theorem \ref{thm_triplet},
      \beas
  \hspace{-2mm}   V_0 \n = \n   \underset{\t  \in \cT }{\sup} \,   \underset{\ga  \in \cT }{\inf} \,
   \underset{\hP \in \cP }{\inf} \, \hE_\hP  \big[  R (\t,\ga)    \big]
    \n =  \n   \underset{\ga  \in \cT }{\inf} \,   \underset{\t  \in \cT }{\sup} \,
    \underset{\hP \in \cP }{\inf} \,   \hE_\hP  \big[  R (\t,\ga)    \big]
  \n = \n   \hE_{\hP_*} \n  \big[ R (\tau_*,\ga_*)\big]
  \n  \ge  \n   \underset{\hP \in \cP }{\inf} \, \hE_\hP \n  \big[ R (\tau_*,\ga_*)\big]
    \n  \ge  \n    \underset{(\hP, \ga) \in \cP \times \cT}{\inf} \hE_{\hP } \n  \big[ R (\tau_*,\ga )\big]
    \n = \n  V_0 .
      \eeas
      In particular, one has
   $       \underset{\t  \in \cT }{\sup} \,   \underset{\ga  \in \cT }{\inf} \,
   \underset{\hP \in \cP }{\inf} \, \hE_\hP  \big[  R (\t,\ga)    \big]
   =   \underset{\ga  \in \cT }{\inf} \,   \underset{\t  \in \cT }{\sup} \,
    \underset{\hP \in \cP }{\inf} \,   \hE_\hP  \big[  R (\t,\ga)    \big]
   =   \underset{\hP \in \cP }{\inf} \, \hE_\hP \n  \big[ R (\tau_*,\ga_*)\big]  $.
 Multiplying  $-1$ to each side shows that
 $  \underset{\tau  \in \cT  }{\inf} \;  \underset{\ga  \in \cT  }{\sup} \; \vf    \big(  R (\tau, \ga)   \big)
  =  \underset{\ga  \in \cT  }{\sup} \; \underset{\tau  \in \cT  }{\inf} \;  \vf    \big(  R (\tau, \ga)   \big)
  = \vf \big( R (\tau_*, \ga^*) \big)  $
 is  Player 2's value in the corresponding Dynkin game under the worst-case  risk measure $\vf$,
  of which     $(\tau_*,\ga^*)$ is a saddle point.

\end{rem}

\fi

\begin{rem}
  Theorem \ref{thm_RDG} \(1\) and Theorem \ref{thm_triplet} imply that
\beas
V_0  \n = \n  \hE_{\hP_*} \n  \big[ R (\tau_*,\ga_*)\big] \n \ge \n \underset{\hP \in \cP}{\inf}
\hE_\hP  \big[ R (\tau_*,\ga_*)\big] \n \ge \n \underset{\ga \in \cT}{\inf} \, \underset{\hP \in \cP}{\inf}
\hE_\hP  \big[ R (\tau_*,\ga )\big] \n = \n V_0 ,
\eeas
which shows that $ V_0  \n = \n    \underset{\hP \in \cP}{\inf}
\hE_\hP  \big[ R (\tau_*,\ga_*)\big] \n = \n \ul{\sE}_0 \big[ R (\tau_*,\ga_*)\big] $.
Hence, we   see that the pair $ (\tau_*,\ga_*) $ is robust with respect to $\hP \n \in \n \cP$, or
$ (\tau_*,\ga_*) $ is a saddle point of the Dynkin game under the nonlinear expectation $\ul{\sE}_0$.
\end{rem}

  \section{Proofs}
  \label{sec:proofs}

    \subsection{Proofs of technical results in Sections \ref{subsec:preliminary}, \ref{sec:WSUP} and \ref{sec:DPP}}

     \label{subsect:proof_ROSVU}

      \no {\bf Proof of Proposition \ref{prop_shift0} (2):}
 Let $n  \n \in \n  \hN$ and $\t  \n \in \n  \cT^t (n) $. Assume that
 $\t(\o  \n \otimes_s \n  \O^s )  \n \subset \n  [r,T]$ for some $r   \n \in \n     [s,T]$.
 For any $ i  \n = \n  0, \cds, 2^n$ such that $t^n_i \= t  \n \vee \n  (i 2^{-n} T)  \n \ge \n  r$,
 since $r  \n \ge \n  s  \n \ge \n  t$, one has $\wt{r}  \n := \n  t  \n \vee \n  (i 2^{-n} T)
  \n = \n    \big( t  \n \vee \n  (i 2^{-n} T) \big) \n \vee \n  s  \n = \n s  \n \vee \n  (i 2^{-n} T) $.
 Setting $A  \n := \n  \{ \o'  \n \in \n  \O^t \n : \t (\o')  \n \le \n  \wt{r}  \}  \n \in \n  \cF^t_{\wt{r} } $,
 we can deduce from Lemma 2.2 of \cite{ROSVU}   that
  \beas
  \{ \wt{\o} \in \O^s : \t^{s,\o} (\wt{\o}) \le \wt{r} \, \}
  =   \{ \wt{\o} \in \O^s : \t  ( \o \otimes_s \wt{\o}) \le \wt{r} \, \}
  =   \{ \wt{\o} \in \O^s :      \o \otimes_s \wt{\o}  \in A  \}
  =  A^{s,\o}  \in \cF^s_{ \wt{r} }   .
  \eeas
  So  $\t^{s,\o}$ is an $\bF^s-$stopping time  valued in
  $ \{t  \n \vee \n  (i 2^{-n} T) \n \in \n [r,T] \n : i \n = \n  0, \cds \n , 2^n \}
  \n \subset \n \{s  \n \vee \n  (i 2^{-n} T) \n \in \n [r,T] \n : i \n = \n  0, \cds \n , 2^n \} $,
  i.e. $\t^{s,\o} \ins \cT^s_r (n) $.

  For the case of   $n \= \infty$, see Corollary 2.1 of \cite{ROSVU}.   \qed

     \no {\bf Proof of \eqref{eq:ab025}:} Let $t \in [0,T] $ and $ \o_1 , \o_2 \in \O$. We see from
 \eqref{eq:aa211} that
 \beas
   -L_t(\o_1)  & \tn \le & \tn  -L_t(\o_2) + |L_t(\o_1)-L_t(\o_2)| \le  \Psi_t(\o_2)
  + \rho_0 \big( \|\o_1 \n - \n \o_2 \|_{0,t} \big) ,   \\
 \hb{and} \q  U_t (\o_1)  & \tn \le & \tn   U_t(\o_2) + |U_t(\o_1)-U_t(\o_2)| \le  \Psi_t(\o_2)
 + \rho_0 \big( \|\o_1 \n - \n \o_2 \|_{0,t} \big) .
 \eeas
 It follows that
 $\Psi_t (\o_1)   \= \big( \- L_t (\o_1) \big) \vee   U_t (\o_1) \vee 0 \ls
  \Psi_t(\o_2)  \+ \rho_0 \big( \|\o_1 \n - \n \o_2 \|_{0,t} \big) $. Then exchanging the roles of
  $\o_1$ and $\o_2$ proves \eqref{eq:ab025}. \qed

        \no {\bf Proof of Lemma \ref{lem_Phi_integr}:}
   Let  $t \n \in \n  [0,T]$ and   $\hP  \n \in \n  \fP_t $.
    Suppose  that  $   \Psi^{t,\o}   \n \in \n   \hS   (\bF^t, \hP )$ and
   $\hE_\hP \n \int_t^T |g^{t,\o}_s| ds < \infty$   for some $\o  \n \in \n  \O $.
   Let $\o' \in \O$.      For any $(s,\wt{\o}) \n \in \n [t,T] \times \O^t$,    \eqref{eq:aa211} implies that
   \bea    \label{eq:bb373}
         \big| g^{t,\o'}_s(   \wt{\o}) \n  -  \n  g^{t,\o}_s(   \wt{\o} ) \big|
 \n = \n \big| g_s( \o' \n \otimes_t \n  \wt{\o})  \n - \n   g_s( \o  \n \otimes_t \n  \wt{\o} ) \big|
 \n  \le    \n \rho_0 \big( \|\o'  \n \otimes_t \n  \wt{\o} \n - \n \o  \n \otimes_t \n  \wt{\o} \|_{0,s} \big)
 \n = \n  \rho_0 \big( \|\o' \n - \n \o\|_{0,t} \big)    ,
   \eea
   so $\hE_\hP \int_t^T \n  |g^{t,\o'}_s | ds
      \n \le  \n  \hE_\hP  \int_t^T \n  |g^{t,\o}_s | ds
       \n + \n  (T \n - \n t) \rho_0 \big( \|\o' \n - \n \o\|_{0,t} \big) \n < \n \infty$\,.

   Proposition \ref{prop_shift0} (4) shows    that both $L^{t,\o'}$ and $U^{t,\o'}$ are $ \bF^t-$adapted processes
   with all continuous  paths,
   so is the process $\Psi^{t,\o'}_s = \big( \n - \n L^{t,\o'}_s \big)\vee   U^{t,\o'}_s \vee 0 $, $s \in [t,T]$.
  Similar to \eqref{eq:bb373},    we see from   \eqref{eq:ab025}      that
   \beas
            \big| \Psi^{t,\o'}_s(   \wt{\o}) \n  -  \n  \Psi^{t,\o}_s(   \wt{\o} ) \big|
 \n  \le    \n 
 \rho_0 \big( \|\o' \n - \n \o\|_{0,t} \big)  , \q  \fa (s,\wt{\o}) \n \in \n [t,T] \times \O^t  .
 \eeas
   It follows that   $ \hE_\hP \Big[ \Psi^{t,\o'}_*   \Big] \n = \n
  \hE_\hP \bigg[ \underset{s \in [t,T]}{\sup} \big|  \Psi^{t,\o'}_s  \big| \bigg]
       \n \le \n   \hE_\hP \bigg[ \underset{s \in [t,T]}{\sup} |  \Psi^{t,\o}_s  | \bigg]
       \n + \n    \rho_0 \big( \|\o' \n - \n \o\|_{0,t} \big)
        \n = \n   \hE_\hP \big[ \Psi^{t,\o}_*   \big]
        \n + \n    \rho_0 \big( \|\o' \n - \n \o\|_{0,t} \big) \n < \n \infty $.
        Therefore,  $   \Psi^{t,\o'}   \n \in \n   \hS  \big( \bF^t ,   \hP   \big)$.    \qed

  \no {\bf Proof of Remark \ref{rem_fP_LU} (1):}
 Let $t \n \in   \n  [0,T]$.
   Proposition \ref{prop_shift1}     implies that  for $\hP_0-$a.s. $\o \in \O$,
   $\Psi^{t,\o} \in  \hS  \big(\bF^t,(\hP_0)^{t,\o}\big)= \hS  \big(\bF^t,\hP^t_0 \big) $ and
   \beas
   \hE_{\hP^t_0}   \n  \int_t^T \n  |g^{t,\o}_s | ds
   \n = \n  \hE_{(\hP_0)^{t,\o}} \dn  \left[   \bigg( \n  \int_t^T \n  |g_s | ds   \bigg)^{t,\o} \right]
    \n \le \n  \hE_{(\hP_0)^{t,\o}} \dn  \left[   \bigg( \n   \int_0^T \n  |g_s | ds   \bigg)^{t,\o} \right]
    \n = \n  \hE_{\hP_0}  \dn \left[      \int_0^T \n  |g_s | ds   \bigg| \cF_t \right]  \n  (\o)
     \n < \n  \infty .
   \eeas
    It then follows from
    Lemma \ref{lem_Phi_integr} that   $\Psi^{t,\bz} \n \in \n   \hS  \big(\bF^t, \hP^t_0 \big)$
    and $ \hE_{\hP^t_0}     \n  \int_t^T \n  |g^{t,\bz}_s | ds   \n < \n \infty $.
    Hence,   $ \hP^t_0 \in \wh{\fP}_t $. \qed

  \no {\bf Proof of Remark \ref{rem_P3}:}
  {\bf 2)} Fix  $t \n \in \n  [0,T]$ and let
  $ \o_1, \o_2  \n \in \n  \O  $, $\tau  , \ga  \n \in \n  \cT^t$.
 By \eqref{eq:da021} and  \eqref{eq:da025},
       \bea
   \big| \big( R^{t,\o_1} (\tau, \ga) \big) (\wt{\o})
 \n  - \n \big( R^{t,\o_2} (\tau, \ga) \big) (\wt{\o}) \big|
  & \tn \dn =  & \tn \dn   \big| R \big(t, \tau (\wt{\o}), \ga (\wt{\o}) , \o_1  \n \otimes_t \n  \wt{\o} \big)
  - R \big(t, \tau (\wt{\o}), \ga (\wt{\o}) , \o_2  \n \otimes_t \n  \wt{\o} \big) \big|  \nonumber  \\
    & \tn  \dn    \le  & \tn  \dn     (1 \n + \n  T  )
    \rho_0 \big( \|\o_1 \n \otimes_t \n  \wt{\o} \n - \n \o_2 \n \otimes_t \n  \wt{\o} \|_{0,T} \big)
   \n  =  \n  (1 \n + \n  T  )  \rho_0 \big( \|\o_1 \n - \n \o_2 \|_{0,t} \big) ,
   \q \fa \wt{\o} \in \O^t. \q \qq  \label{eq:ab021}
  \eea

 Now, let $\o \n \in \n  \O$,  $s  \n \in \n  [t,T]$,
 $n \in \hN \cup \{\infty\}$, $\wt{\hP} \ins \wh{\fP}_s$ and $ \wp  \n \in \n  \cT^s$.
 Given   $ \wt{\o}_1,\wt{\o}_2  \n \in \n  \O^t $  and $\vs \n \in \n  \cT^s(n) $,
 similar to \eqref{eq:ab021},
 \bea \label{eq:xxx867}
 \big| R^{s,\o \otimes_t \wt{\o}_1 }  (\vs,\wp) \n  - \n  R^{s,\o \otimes_t \wt{\o}_2 }  (\vs,\wp) \big|
 \le  (1 \n + \n  T  )  \rho_0 \big( \| \o \otimes_t \wt{\o}_1 \n - \n \o \otimes_t \wt{\o}_2 \|_{0,s} \big)
 =  (1 \n + \n  T  )  \rho_0 \big( \|   \wt{\o}_1 \n - \n   \wt{\o}_2 \|_{t,s} \big) .
 \eea
 It follows that
    $  \hE_{\wt{\hP}} \big[ R^{s,\o \otimes_t \wt{\o}_1}  (\vs,\wp) \big] \ls
   \hE_{\wt{\hP}} \big[ R^{s,\o \otimes_t \wt{\o}_2}  (\vs,\wp) \big]
    \+  (1 \n + \n  T  )  \rho_0 \big(\|   \wt{\o}_1  \n  -  \n  \wt{\o}_2 \|_{t,s} \big) $.
   Taking supremum over $\vs \n \in \n \cT^s (n) $ yields that
    $  \underset{\vs \in \cT^s(n) }{\sup} \hE_{\wt{\hP}} \big[ R^{s,\o \otimes_t \wt{\o}_1}  (\vs,\wp) \big] \ls
  \underset{\vs \in \cT^s(n) }{\sup} \hE_{\wt{\hP}} \big[ R^{s,\o \otimes_t \wt{\o}_2}  (\vs,\wp) \big]
  \+  (1 \n + \n  T  )  \rho_0 \big(\|   \wt{\o}_1 -   \wt{\o}_2 \|_{t,T} \big)$.
  Exchanging the roles of $\wt{\o}_1$ and $\wt{\o}_2$ shows that
   the mapping $ \wt{\o} \to \underset{\vs \in \cT^s(n) }{\sup}
   \hE_{\wt{\hP}}   \big[ R^{s,\o \otimes_t \wt{\o}}  (\vs,\wp) \big]$ is continuous
   under   norm    $\|~\|_{t,T}$ and thus $\cF^t_T-$measurable.

       Next, let us show that both sides of  \eqref{eq:xxx617} are finite:
          Let $A \n \in \n  \cF^t_s$, $\tau \ins \cT^t_s (n) $ and $j \= 1, \cds \n , \l $.
          By    \eqref{eq:ab015} and \eqref{eq:xxx111},
   \beas
      \Big| \hE_{\wh{\hP} } \big[ \b1_{A \cap \cA_j} R^{t,\o  }  \big(\tau,\wp^n_j \big) \big] \Big|
    \ls \hE_{\wh{\hP} } \Big[  \big|  R^{t,\o  }  \big(\tau,\wp^n_j \big)  \big| \Big]
   \le \hE_{\wh{\hP} } \bigg[  \int_t^{\tau  \land \wp^n_j} |g^{t,\o}_s| ds  +   \Psi^{t,\o}_{\tau  \land \wp^n_j}  \bigg]
    \ls \hE_{\wh{\hP} } \bigg[  \int_t^T |g^{t,\o}_s| ds +   \Psi^{t,\o}_*  \bigg] < \infty .
   \eeas

    On the other hand, given $\wt{\o} \n \in \n  A \cap \cA_j$ and $\vs  \n \in \n  \cT^s (n) $,  taking
       $(\wt{\o}_1, \wt{\o}_2)  \n = \n  (\wt{\o},\wt{\o}_j)$   in \eqref{eq:xxx867},
    we can deduce from   \eqref{eq:ab015} and \eqref{eq:xxx111} again that
    \beas
    \hspace{-6mm}   \Big| \hE_{\hP_j}   \big[ R^{s,\o   \otimes_t    \wt{\o}}  (\vs, \wp )  \big] \Big|
         \ls       \hE_{\hP_j}   \Big[  \big|  R^{s,\o    \otimes_t    \wt{\o}_j} (\vs  , \wp )  \big| \Big]
       \dn + \n ( 1 \n + \n T )
      \rho_0 \big(\|   \wt{\o}   \n - \n    \wt{\o}_j \|_{t,s} \big)
      \n \le  \n  \hE_{\hP_j} \n   \bigg[ \int_s^T \n \big| g^{s,\o    \otimes_t    \wt{\o}_j}_r \big| dr \dn + \n
      \Psi^{s,\o    \otimes_t    \wt{\o}_j}_* \bigg]  \n + \n ( 1 \n + \n T )    \rho_0  (\d  )
      \df \a_j \< \infty .
    \eeas
        It then follows 
        that
    \beas
    \hE_{ \hP  } \bigg[ \b1_{\{\wt{\o} \in A \cap  \cA_j\}} \bigg(   \underset{\vs \in \cT^s (n) }{\sup}
    \hE_{\hP_j}   \big[ R^{s,\o   \otimes_t    \wt{\o}}  (\vs, \wp )  \big]
    \n + \n \int_t^s g^{t,\o}_r (\wt{\o}) dr    \bigg) \bigg]
       \n \le \n \hE_\hP \bigg[ \b1_{   A \cap  \cA_j } \n \int_t^T \big| g^{t,\o}_r \big| dr \bigg]
     \n + \n  \a_j   \hP(A  \n \cap \n  \cA_j)  \n < \n  \infty ,
    \eeas
    as well as that
      \beas
   \hE_{ \hP  } \bigg[ \b1_{\{\wt{\o} \in A \cap  \cA_j\}} \bigg(   \underset{\vs \in \cT^s  (n)  }{\sup}
    \hE_{\hP_j}   \big[ R^{s,\o   \otimes_t    \wt{\o}}  (\vs, \wp )  \big]
    \n + \n \int_t^s g^{t,\o}_r (\wt{\o}) dr    \bigg) \bigg]
     \n \ge \n - \hE_\hP \bigg[ \b1_{   A \cap  \cA_j } \n \int_t^T \big| g^{t,\o}_r \big| dr \bigg]
       \- \a_j   \hP(A  \n \cap \n  \cA_j)  \n > \n  - \infty \, .
    \eeas
 Summing both up over $j \ins \{1,\cds \n , \l\}$ shows that
 the right-hand-side of  \eqref{eq:xxx617} is finite.

  \no {\bf 3)} The proof of Remark 3.3 (2) in \cite{ROSVU} has shown that
 the probability $\wh{\hP}$ defined in \eqref{eq:xxx131c} satisfies (P3) (i) and  (ii):
 $ \wh{\hP}  (A \cap \cA_0 ) \n =  \n  \hP   (A \cap \cA_0 )  $, $ \fa   A  \n \in \n  \cF^t_T $, and $
 \wh{\hP}  (A \cap \cA_j )  \n = \n   \hP  (A \cap \cA_j )$, $ \fa  j \n = \n 1,\cds  \n ,\l $,
   $ \fa    A  \n \in  \n  \cF^t_s $.
 \if{0}
 Given $  A \n \in \n  \cF^t_T$, for any $j  \n = \n  1,  \cds \n , \l $ and $\wt{\o}  \n \in \n  \cA_j $,
    since $\cA_j  \n \in \n  \cF^t_s$,
    Lemma \ref{lem_element}  shows that $( \cA_j  )^{s , \wt{\o}} = \O^s $,
    which implies that  $(A \cap \cA_0 )^{s , \wt{\o}} = \es$. Then one can easily calculate  that
     $\wh{\hP}  (A \cap \cA_0 ) =  \hP   (A \cap \cA_0 ) $.

       Next, let $j = 1,  \cds \n , \l $ and  $A \in \cF^t_s$.
  We see from Lemma \ref{lem_element}  again that
  \bea  \label{eq:xxx614}
 \hb{ if $\wt{\o} \in A \cap \cA_j$ (resp. $\notin A \cap \cA_j$),
 then $( A \cap \cA_j )^{s , \wt{\o}}  = \O^s  $ (resp. $= \es $). }
 \eea
    It follows that
\beas
\wh{\hP}  (A \cap \cA_j ) =  \sum^\l_{j'=1} \hE_\hP   \n  \left[   \b1_{\{\wt{\o} \in \cA_{j'}\}}
  \hP_{j'}  \big( (A \cap \cA_j )^{s,\wt{\o}}  \big)    \right]
 = \sum^\l_{j'=1} \hE_\hP   \n  \left[ \b1_{\{\wt{\o} \in A \cap \cA_j \}}  \b1_{\{\wt{\o} \in \cA_{j'} \}}
  \hP_{j'}  \big(\O^s\big)  \right]   = \hP(A \cap \cA_j) .
\eeas
 \fi
 To see   $\wh{\hP}$  satisfying \eqref{eq:xxx617},
 let us fix $ n \in \hN \cup \{\infty\}$   and $\wp \in \cT^s$.
 We set $\wp^n_j : = \wp (\Pi^t_s)$, $j \n = \n  1,  \cds \n , \l $, which are of $\cT^t_s$ by \eqref{eq:xax001}.

 Let $A  \n \in \n  \cF^t_s$   and $\tau  \n \in \n  \cT^t_s (n) $.
 Given $ \wt{\o}  \in \O^t $,
 Proposition \ref{prop_shift0} (2) shows that  $\t^{s,\wt{\o}} \in \cT^s (n)$.
 Since the $\bF-$adaptness of $g$ and \eqref{eq:bb421} imply  that
 \bea \label{eq:xax141}
    g_r (\o \otimes_t \O^t) = g_r (\o) , \q  \fa  r \in [0,t]  \q \hb{and} \q
  g_r \big( (\o \otimes_t \wt{\o} ) \otimes_s \O^s \big) = g_r (\o \otimes_t \wt{\o}) , \q \fa  r \in [0, s] ,
  \eea
 we see from \eqref{eq:da021}    that for any $\wh{\o} \in \O^s$
 \bea
 &&  \hspace{-1cm}  \big(  R^{t,\o}  (\tau , \wp^n_j)   \big)^{s,\wt{\o}} (\wh{\o})
     =  \big(  R^{t,\o} \big( \tau, \wp^n_j \big)   \big)  (\wt{\o} \otimes_s \wh{\o}   )
     = R \big(t,  \tau (\wt{\o} \otimes_s \wh{\o} ) , \wp \big( \Pi^t_s (\wt{\o} \otimes_s \wh{\o} ) ) \big)
     ,  \o \otimes_t ( \wt{\o}   \otimes_s \wh{\o} ) \big)       \nonumber    \\
 &&     = \n  R \big(s,  \tau^{s,\wt{\o}} (  \wh{\o} ) , \wp  (   \wh{\o} )
     , (\o \otimes_t \wt{\o} ) \otimes_s \wh{\o} \big)
     \n  +  \n  \int_t^s  \n  g_r \big( (\o \otimes_t \wt{\o} ) \otimes_s \wh{\o} \big) dr
      \n = \n  \big( R^{s, \o \otimes_t \wt{\o} } (\tau^{s,\wt{\o}}, \wp) \big) (\wh{\o})
       \n + \n   \int_t^s  \n  g_r ( \o \otimes_t \wt{\o} ) dr  . \qq     \label{eq:ab027}
 \eea
 By  Lemma \ref{lem_element},
   $( A \cap \cA_j )^{s , \wt{\o}}  = \O^s  $ (resp. $= \es $)
 \, if \, $\wt{\o} \in A \cap \cA_j$ (resp. $\notin A \cap \cA_j$).
 Then \eqref{eq:ab027} leads to  that
   \beas
  &&  \hspace{-0.7cm}  \hE_{\wh{\hP} }  \big[  \b1_{ A \cap \cA_j}   R^{t,\o} (\tau, \wp^n_j)       \big]
    \n   =     \n   \sum^\l_{j' = 1}    \hE_{ \hP  }  \bigg[
  \b1_{ \{ \wt{\o} \in  \cA_{j'}  \} } \hE_{ \hP_{j'}}
  \Big[ \big( \b1_{   A \cap   \cA_j}   R^{t,\o} (\tau, \wp^n_j)  \big)^{s,\wt{\o}} \Big]   \bigg]   \\
 && = \n \sum^\l_{j' = 1}    \hE_{ \hP  }  \bigg[
  \b1_{ \{ \wt{\o} \in  A \cap \cA_j  \} } \b1_{ \{ \wt{\o} \in  \cA_{j'}  \} } \hE_{ \hP_{j'}}
  \Big[ \big(    R^{t,\o} (\tau, \wp^n_j)  \big)^{s,\wt{\o}} \Big]   \bigg]
  \n   =     \n  \hE_{ \hP  }  \Bigg[  \b1_{ \{ \wt{\o} \in  A \cap  \cA_j  \} }
 \bigg( \hE_{\hP_j} \bigg[   R^{s, \o   \otimes_t    \wt{\o}} \big(\tau^{s,\wt{\o}}, \wp \big) \bigg]
  \n +  \n   \int_t^s g^{t,\o}_r  (     \wt{\o}      ) dr   \bigg)  \Bigg] \q \\
 && \le  \hE_{ \hP  }  \bigg[  \b1_{ \{ \wt{\o} \in  A \cap  \cA_j  \} }
 \bigg(    \underset{\vs \in \cT^s (n) }{\sup}
 \hE_{\hP_j} \big[  R^{s, \o   \otimes_t    \wt{\o}}  (\vs, \wp  )  \big]
 \n + \n \int_t^s g^{t,\o}_r  (     \wt{\o}      ) dr  \bigg)  \bigg] .
     \eeas
 Taking summation  over $j \ins \{1,\cds \n  ,\l\}$ yields \eqref{eq:xxx617}. \qed

  \no {\bf Proof of \eqref{eq:uxu170}:}
  Let $(t,\o) \in [0,T] \times \O$.
 Since the $\cF_t-$measurability of $L_t$, $U_t$ and \eqref{eq:bb421} show  that
 \bea \label{eq:uxu205}
 L^{t,\o}_t (\wt{\o}) = L_t(\o \otimes_t \wt{\o}) = L_t(\o) \q \hb{and} \q
 U^{t,\o}_t (\wt{\o}) = U_t(\o \otimes_t \wt{\o}) = U_t(\o) , \q \fa \wt{\o} \in \O^t  .
 \eea
   it holds for any $\tau \n \in \n  \cT^t (n) $ that
   $   R^{t,\o}  (\tau  ,t ) \n = \n    \b1_{\{\tau  = t\}} L^{t,\o}_\tau  \n + \n  \b1_{\{ t < \tau   \}} U^{t,\o}_t
       \n = \n    \b1_{\{\tau  = t\}} L^{t,\o}_t   \n + \n  \b1_{\{ t < \tau   \}} U^{t,\o}_t
        \n  \le  \n  U^{t,\o}_t  \n  = \n   U_t   ( \o ) $. So
     \beas 
  V^n_t (\o)  \le     \underset{\hP \in \cP(t,\o) }{\inf} \,
  \underset{\t  \in \cT^t(n) }{\sup} \,   \hE_\hP  \big[  R^{t,\o} (\t,t)    \big]
  \le   \underset{\hP \in \cP(t,\o) }{\inf} \,  \hE_\hP  \big[  U_t (\o)   \big] = U_t (\o) \le \Psi_t (\o).
    \eeas

    On the other hand, since $t \in \cT^t(n)$ and
    since  $ R^{t,\o}  (t,\ga )  \n  =  \n   \b1_{\{t   \le \ga \}} L^{t,\o}_t
       \n + \n  \b1_{\{ \ga < t   \}}  U^{t,\o}_\ga
       \n = \n    L^{t,\o}_t
        \n =  \n    L_t (\o) $ for any $\ga  \n \in \n  \cT^t$,
      \beas 
  \hspace{2.7cm}
  V^n_t (\o)  \ge   \underset{\hP \in \cP(t,\o) }{\inf} \,   \underset{\ga  \in \cT^t }{\inf} \,
      \hE_\hP  \big[  R^{t,\o} (t,\ga)    \big]
       =    \underset{\hP \in \cP(t,\o) }{\inf} \, \hE_\hP \big[  L_t (\o) \big]  =  L_t (\o) \ge - \Psi_t (\o)   .
  \hspace{2.7cm} \hb{\qed}
    \eeas

             \no {\bf Proof of Remark \ref{rem_V_conti}:}
       Fix $n  \n \in \n  \hN \cup \{\infty\}$.
   Let $t \n \in \n  [0,T]$,    $\o_1, \o_2  \n \in \n  \O$,  $\hP  \n \in \n  \cP_t  $ and
      $\tau, \ga  \n \in \n  \cT^t$. By \eqref{eq:ab021},
     $    \hE_\hP \big[  R^{t,\o_1} (\tau,\ga)    \big]
   \n \le  \n  \hE_\hP \big[ R^{t,\o_2} (\tau,\ga)      \big]
    \n +  \n  (1 \n + \n  T  ) \rho_0 \big(  \|\o_1  \n - \n  \o_2\|_{0,t}  \big)    $.
     Taking supremum over $\t \n \in \n  \cT^t (n)$, taking infimum over $\ga  \n \in \n  \cT^t$
      and then taking infimum over $\hP \n \in \n  \cP_t$ yield  that
   $  V^n_t(\o_1)  \n \le \n  V^n_t(\o_2)  \n  + \n
  (1 \n + \n  T  ) \rho_0 \big(  \|\o_1  \n - \n  \o_2\|_{0,t}  \big)    $.
   Exchanging the roles of $\o_1$ and $\o_2$, we obtain  \eqref{eq:aa213} with
   $\rho_1  \n = \n   (1 \n + \n  T  ) \rho_0 $ for each $n \in \hN \cup \{\infty\}$. \qed

  \subsection{Proofs of the Dynamic Programming Principles}

      \no {\bf Proof of Proposition \ref{prop_DPP}:}
   Fix $n \in \hN \cup \{\infty\}$,  $0 \n \le  \n  t  \n \le \n  s  \n \le \n  T $ and   $  \o  \n \in \n  \O $.

 \no {\bf 1)} When $t \n = \n s$, since $ V^n $ is $\bF-$adapted by
 Remark \ref{rem_V_adapted},
 an analogy to \eqref{eq:uxu205}
 shows that $ \big( V^n \big)^{t,\o}_t (\wt{\o})   \n  = \n  V^n     (t, \o \otimes_t \wt{\o})
    \n = \n  V^n_t      (  \o )  $, $\fa \wt{\o}  \n \in \n  \O^t$.
 Then
 \beas
 \underset{\hP \in \cP(t,\o)}{\inf} \, \underset{\ga \in \cT^t}{\inf} \,
  \underset{\t \in \cT^t(n)}{\sup} \,  \hE_\hP \Big[ \b1_{\{\tau \land \ga < t\}}   R^{t,\o}  (\tau, \ga)
   + \b1_{\{\tau \land \ga \ge t\}}   \big( V^n_t \big)^{t,\o}        \Big]
   = \underset{\hP \in \cP(t,\o)}{\inf} \, \hE_\hP \big[   V^n_t (  \o )  \big]
   = V^n_t (  \o ) .
 \eeas

 \no {\bf 2)} {\it  To demonstrate \eqref{eq:bb013} for case $t \n < \n s$,
 we shall paste the local approximating $\hP-$minimizers of $ ( V^n_s )^{t,\o}   $
  according to \(P3\)   and then make some estimations. }


 \no {\bf 2a)} Under  norm $\|\cd\|_{t,T} $,   since  $\O^t$ is  a separable complete metric space,
 there exists     a  countable  dense subset $\big\{ \wh{\o}^t_j \big\}_{j \in \hN}$  of $\O^t$.
  Fix    $\e  \n > \n 0$ and let $\d \n \in \n \hQ_+$  satisfy
 $  \dis  \rho_1 (\d)     \n  \vee  \n    \wh{\rho}_0 (\d)     \n  \vee  \n    \big( (1 \n + \n T) \rho_0 (\d) \big)
     \n < \n   \e / 5   $.
 Let $j \n \in \n \hN$.
 By \eqref{eq:bb237},   $\cA_j  \n  := \n     O^s_\d (\wh{\o}^t_j)
   \big\backslash \big( \underset{j'<j}{\cup} O^s_\d (\wh{\o}^t_{j'}) \big)     \\ \in\n  \cF^t_{s} $.
 We can find a  $\hP_j     \n  \in \n  \cP (s, \o \otimes_t \wh{\o}^t_j ) $
   and  a $\ga_j     \n \in \n  \cT^s$
   such that
   \bea   \label{eq:bb417}
   V^n_s (\o \otimes_t \wh{\o}^t_j )
   \ge  \underset{\ga \in \cT^s}{\inf} \,\underset{\t  \in \cT^s (n) }{\sup} \,  \hE_{  \hP_j  }
    \big[  R^{s,\o \otimes_t \wh{\o}^t_j} (\t,\ga)  \big] - \frac15 \e
   \ge \underset{\t  \in \cT^s (n)   }{\sup} \,  \hE_{  \hP_j  }
   \big[   R^{s,\o \otimes_t \wh{\o}^t_j} (\t,\ga_j)  \big] - \frac25 \e .
   \eea

 Given $\wt{\o}  \in O^s_{\d } ( \wh{\o}^t_j )   $, an analogy    to   \eqref{eq:xxx867} shows that
 for any $\tau \in \cT^s (n)$
 \beas
 \q  \big|   R^{s,\o \otimes_t \wt{\o}}   (\t,\ga_j)
  \n - \n     R^{s,\o \otimes_t \wh{\o}^t_j }   (\t,\ga_j)   \big|
  \n \le \n 
   (1 \n + \n T) \rho_0 \big(  \| \wt{\o}  \n - \n  \wh{\o}^t_j \|_{t,s} \big)
   \n \le \n (1 \n + \n T) \rho_0 (\d)   \n \le \n   \frac15   \e    ,
 \eeas
 so  $   \hE_{  \hP_j }  \n \big[ R^{s,\o \otimes_t \wt{\o}}   (\t,\ga_j)  \big]
         \n  \le   \n    \hE_{  \hP_j }  \n  \big[ R^{s,\o \otimes_t \wh{\o}^t_j} (\t,\ga_j)  \big]
         \n  +  \n     \e / 5     $.
       Taking supremum over $ \t  \n  \in \n  \cT^s (n) $,
       we see from    \eqref{eq:bb417} and \eqref{eq:aa213}  that
    \bea
    \hspace{-5mm}  \underset{\t  \in \cT^s (n) }{\sup}
    \hE_{  \hP_j }  \n  \big[  R^{s,\o    \otimes_t    \wt{\o}}   (\t,\ga_j)  \big]
  & \tn \tn \le & \tn \tn   \underset{\t  \in \cT^s (n) }{\sup}
    \hE_{  \hP_j }  \n   \big[  R^{s,\o    \otimes_t    \wh{\o}^t_j} (\t,\ga_j) \big]
       \dn   + \n  \frac15   \e
    \n  \le  \n   V^n_s     \big(\o  \n \otimes_t \n  \wh{\o}^t_j \big)  \dn + \n  \frac35   \e
    \n   \le  \n    V^n_s   ( \o  \n \otimes_t  \n  \wt{\o})
    \n + \n    \rho_1  \big(  \|   \o  \n \otimes_t  \n  \wt{\o}
     \dn -  \n   \o  \n \otimes_t  \n   \wh{\o}^t_j \|_{0,s}  \big) \dn + \n  \frac35  \e    \nonumber  \\
   & \tn \tn   =  & \tn \dn   ( V^n_s )^{t,\o}(   \wt{\o})
    \n + \n    \rho_1  \big(  \|   \wt{\o}  \n -  \n    \wh{\o}^t_j \|_{t,s}  \big) \n + \n  \frac35  \e
     \n \le \n  ( V^n_s )^{t,\o}(   \wt{\o}) \n + \n \rho_1 (\d) \n + \n  \frac35   \e
    \n \le \n  ( V^n_s )^{t,\o}(   \wt{\o})  \n + \n  \frac45   \e   .
    \label{eq:aa103}
    \eea

   Next,    fix  $\hP \n \in \n  \cP (t,\o) $,   $\l  \n \in \n  \hN  $ and
    let $ \wh{\hP}_\l $ be the probability of $\cP(t,\o)$  in (P3)
    for $\big\{(\cA_j, \d_j, \wt{\o}_j, \hP_j) \big\}^\l_{j=1} \=
    \big\{(\cA_j, \d , \wh{\o}^t_j, \hP_j) \big\}^\l_{j=1} $ and
   $\cA_0  \n := \n  \Big(  \underset{ j=1}{\overset{\l}{\cup}} \cA_j  \Big)^c   \n \in \n  \cF^t_{s} $.
   Then we have
    \bea    \label{eq:ff011}
     \hE_{ \wh{\hP}_\l} [\xi]   \n = \n  \hE_\hP [\xi]  , ~
     \fa \xi  \n \in \n  L^1 \big( \cF^t_s, \wh{\hP}_\l \big)  \n \cap \n  L^1 \big( \cF^t_s, \hP \big)
     \q \hb{and} \q     \hE_{ \wh{\hP}_\l} [\b1_{\cA_0}\xi]   \n = \n  \hE_\hP [\b1_{\cA_0}\xi]  ,
     ~ \fa \xi  \n \in \n  L^1 \big( \cF^t_T, \wh{\hP}_\l \big) \n  \cap \n  L^1 \big( \cF^t_T, \hP \big)  .  \q
         \eea
  Also, in light of \eqref{eq:xxx617} and \eqref{eq:aa103},
   there exist  $\wp^n_j  \n \in \n  \cT^t_s$, $ j \n = \n 1 , \cds \n , \l $, such that
  for any $A \n \in \n  \cF^t_s$ and $\tau \ins \cT^t_s (n)$
  \bea    \label{eq:xxx617c}
  \sum^\l_{j=1}  \hE_{ \wh{\hP}_\l}   \big[ \b1_{A \cap \cA_j} R^{t,\o}   (\tau, \wp^n_j  )  \big]
  & \tn  \le & \tn   \sum^\l_{j=1}
  \hE_{ \hP  } \bigg[ \b1_{\{\wt{\o} \in A \cap  \cA_j\}} \bigg(   \underset{\vs \in \cT^s (n) }{\sup}
   \hE_{\hP_j}   \big[ R^{s,\o \otimes_t \wt{\o}}  (\vs, \ga_j)  \big] \n + \n \int_t^s \n g^{t,\o}_r (\wt{\o}) dr
    \bigg) \bigg]  \n + \n  \wh{\rho}_0 (\d) \nonumber \\
  & \tn       \le   & \tn   \hE_{\hP}  \bigg[ \b1_{  A \cap \cA_0^c }
     \Big( ( V^n_s )^{t,\o}     \n + \n  \int_t^s g^{t,\o}_r   dr  \Big)          \bigg] +   \e    .
   \eea

   \no {\bf 2b)} Now,   let $  \ga \n \in \n  \cT^t $  and  $ \tau \n \in \n  \cT^t(n) $.
  Applying \eqref{eq:xxx617c}
   with $A = \{\tau  \n \land \n  \ga   \n \ge \n   s \}  \n \in \n  \cF^t_s $, one can show that
         \bhe \bea
    \sum^\l_{j=1}  \hE_{\wh{\hP}_\l}  \big[ \b1_{\{\t  \land \ga  \ge  s\} \cap \cA_j}
         R^{t,\o}  (\tau   , \wp^n_j )   \big]
     \le   \hE_{\hP}  \bigg[ \b1_{  \{\t \land \ga   \ge s\} \cap \cA_0^c }
     \Big( ( V^n_s )^{t,\o}     \n + \n  \int_t^s g^{t,\o}_r   dr  \Big)          \bigg] +   \e   . \qq \label{eq:da041}
     \eea \ehe
  We glue $\ga$ with $\{\wp^n_j\}^\l_{j=1}$ to form a new $\bF^t-$stopping time
  \bhe \bea \label{eq:uxu151}
     \wh{\ga}_\l  \n := \n  \b1_{\{  \ga    < s\}} \ga  \n + \n  \b1_{\{  \ga   \ge   s\}}
         \bigg( \b1_{\cA_0} \ga  \n + \n  \sum^\l_{j=1} \b1_{\cA_j} \wp^n_j \bigg) .
    \eea \ehe

  Since $  \wh{\ga}_\l \gs s \> \tau $ on $\{  \ga   \gs   s\} \n \cap \n \{ \tau     \< s\} $,
      \eqref{eq:ef031} shows  that
  \beas 
  \b1_{\{ \tau \land \ga    < s\}}  R^{t,\o} ( \t, \wh{\ga}_\l )
  \n = \n \b1_{\{   \ga    < s\}} R^{t,\o} ( \t, \ga) \n + \n  \b1_{\{   \ga    \ge s\} \cap \{ \tau     < s\}  }
  \Big( \int_t^{\tau  } \n g^{t,\o}_s ds   +      L^{t,\o}_{\tau  } \Big)
  \n = \n \b1_{\{ \tau \land \ga    < s\}} R^{t,\o} ( \t, \ga) \n \in \n \cF^t_s .
  \eeas
  Then one    can deduce from   \eqref{eq:ff011}, \eqref{eq:da041},
   \eqref{eq:ab015} and  \eqref{eq:uxu170} that
        \beas
   ~ \;    &  \dn \dn &  \dn \dn \hspace{-1cm}    \hE_{\wh{\hP}_\l}  \big[   R^{t,\o} ( \t, \wh{\ga}_\l )  \big]
       \n = \n  \hE_{\wh{\hP}_\l}  \big[ (\b1_{\{ \tau \land \ga    < s\}}
         \n   +  \n   \b1_{\{\tau  \land \ga  \ge  s\} \cap \cA_0} ) R^{t,\o} (\t, \ga  )  \big]
         \n + \n  \sum^\l_{j=1} \hE_{\wh{\hP}_\l}  \big[ \b1_{\{\t  \land \ga  \ge  s\} \cap \cA_j}
         R^{t,\o}  (\tau   , \wp^n_j )   \big]   \\
    & \dn \dn \le&  \dn \dn   \hE_\hP  \bigg[ \b1_{\{\tau \land \ga    < s\}}  R^{t,\o} (\t, \ga  )
    \n +  \n    \b1_{\{\t \land \ga  \ge s\}}  \Big( ( V^n_s )^{t,\o}     \n + \n  \int_t^s g^{t,\o}_r   dr  \Big)
      \n + \n   \b1_{\{\t \land \ga \ge   s\} \cap \cA_0} \Big( R^{t,\o} (\t, \ga  )  \n - \n  ( V^n_s )^{t,\o}
       \n - \n  \int_t^s g^{t,\o}_r   dr
         \Big)  \bigg]  \n + \n   \e  \nonumber \\
           & \dn \dn \le &  \dn \dn
            \hE_\hP  \bigg[ \b1_{\{\t \land \ga   < s\}} R^{t,\o} (\t, \ga  )
            \n +  \n    \b1_{\{\t \land \ga  \ge s\}}  \Big( ( V^n_s )^{t,\o}     \n + \n  \int_t^s g^{t,\o}_r   dr  \Big)
            \n + \n  \b1_{ \cA_0  }
     \bigg( 2 \int_t^T \big| g^{t,\o}_r   \big| dr + 2 \Psi^{t,\o}_*  \bigg)     \bigg]  \n +  \n  \e  .
        \eeas

     Taking  supremum over $ \t  \n \in \n  \cT^t(n)$  yields that
           \beas
    V^n_t (\o) 
       \ls   \underset{\t \in \cT^t (n) }{\sup} \,  \hE_\hP
          \bigg[ \b1_{\{\t \land \ga   < s\}} R^{t,\o} (\t, \ga  )
        \n  +  \n    \b1_{\{\t \land \ga  \ge s\}}
         \Big( ( V^n_s )^{t,\o}  \n + \n  \int_t^s \n  g^{t,\o}_r   dr  \Big)  \bigg]
        \+   2  \hE_\hP   \bigg[ \b1_{ \cA_0   }
    \bigg(   \int_t^T \n  \big| g^{t,\o}_r   \big| dr \+   \Psi^{t,\o}_*  \bigg)     \bigg]   \n +  \n   \e .
           \eeas
    Then  taking  infimum over $ \ga  \n \in \n  \cT^t$ on the right-hand-side, we obtain
             \beas
          V^n_t (\o)   \n \le \dn   \underset{\ga \in \cT^t }{\inf}   \underset{\t \in \cT^t(n) }{\sup}
       \hE_\hP     \bigg[ \b1_{\{\t \land \ga   < s\}} R^{t,\o} (\t, \ga  )
        \n + \n  \b1_{\{\t \land \ga  \ge s\}}    \Big( ( V^n_s )^{t,\o}  \n + \dn  \int_t^s g^{t,\o}_r   dr  \Big)  \bigg]
        \dn  +   \n 2 \hE_\hP     \bigg[ \b1_{  \big(  \underset{ j=1}{\overset{\l}{\cup}} \cA_j  \big)^c   }
        \bigg(  \n   \int_t^T  \n  \big| g^{t,\o}_r   \big| dr  \n + \n    \Psi^{t,\o}_*  \bigg)    \bigg]
               \n +  \n   \e .
           \eeas

 Since          $       \underset{j \in \hN }{\cup}  \cA_j   \= \underset{j \in \hN }{\cup} O^s_\d (\wh{\o}^t_j)
 \supset \underset{j \in \hN }{\cup} O^T_\d (\wh{\o}^t_j)  \=  \O^t $
     and since
     \bea \label{eq:uxu183}
      \hE_\hP \bigg[ \int_t^T \big| g^{t,\o}_r   \big| dr +   \Psi^{t,\o}_* \bigg]   \n < \n  \infty
      \eea
       by \eqref{eq:xxx111},
         letting $\l \to \infty$,  one can deduce from
       the dominated convergence theorem  that
  \beas
   V^n_t (\o)       \le  \underset{\ga \in \cT^t }{\inf} \,   \underset{\t \in \cT^t(n) }{\sup} \, \hE_\hP
  \bigg[  \b1_{\{\t \land \ga   < s\}} R^{t,\o} (\t, \ga  )
        +    \b1_{\{\t \land \ga  \ge s\}}  \Big( ( V^n_s )^{t,\o}  \n + \n  \int_t^s g^{t,\o}_r   dr  \Big)  \bigg]  +\e  .
   \eeas
 Eventually,  taking infimum over $ \hP \in \cP (t,\o) $ on the right-hand-side and then letting $\e \to 0$
 yield  \eqref{eq:bb013}. \qed

 \if{0}
 \begin{lemm} \label{lem_limit}
 Assume \(A\)  and Assumption \ref{assum_fP_Phi}.
 For any $(t,\o) \in [0,T] \times \O$, let $\tau, \ga \in \cT^t$ and
 let $\{\tau_n\}_{n \in \hN} $ be  a sequence of $\cT^t$ such that
 $ \lmtd{n \to \infty} \tau_n (\wt{\o}) = \tau (\wt{\o}) $, $\fa \wt{\o} \in \O^t$.
 Then it holds  for any $\hP \in \wh{\fP}_t$ that
 $ \hE_\hP \big[ R^{t,\o} (\tau,\ga ) \big] \le \linf{n \to \infty} \hE_\hP \big[ R^{t,\o} (\tau_n ,\ga ) \big] $.
 \end{lemm}

   \no {\bf Proof:} For any $n \in \hN$,  as $L \le U$ by (A),
 \bea
 R^{t,\o} (\tau_n ,\ga ) & \tn = & \tn
 \int_t^{\tau_n \land \ga} \n  g^{t,\o}_s ds  \n + \n  \b1_{\{ \tau_n \le \ga \}} L^{t,\o}_{\tau_n}
  \n + \n  \b1_{\{ \ga < \tau_n   \}} U^{t,\o}_{\ga}
  \n \ge  \n  \int_t^{\tau_n \land \ga}  \n g^{t,\o}_s ds  \n + \n  \b1_{\{ \tau_n \le \ga \}} L^{t,\o}_{\tau_n}
  \n + \n  \b1_{\{ \ga < \tau    \}} U^{t,\o}_{\ga}
  \n + \n  \b1_{\{ \tau \le \ga < \tau_n   \}} L^{t,\o}_{\ga} \nonumber \\
  & \tn = & \tn  \int_t^{\tau_n \land \ga} g^{t,\o}_s ds + \b1_{\{ \tau_n \le \ga \}} L^{t,\o}_{\tau_n}
 + \b1_{\{ \ga < \tau    \}} U^{t,\o}_{\ga} + \b1_{\{ \tau \le \ga < \tau_n   \}}
 L^{t,\o}_{(\tau \vee \ga) \land  \tau_n  } : = \xi_n .  \label{eq:da311}
 \eea
 Set   $ A : = \underset{n \in \hN}{\cup} \{ \tau_n \le \ga \} \subset \{ \tau  \le \ga \} $.
   The  continuity of $L$   implies that
 \bea \label{eq:da313}
 \lmt{n \to \infty}  \xi_n  =
 \int_t^{\tau  \land \ga} g^{t,\o}_s ds +  \b1_A     L^{t,\o}_{\tau }
 + \b1_{\{ \ga < \tau    \}} U^{t,\o}_{\ga} + (\b1_{\{ \tau \le \ga      \}} -  \b1_A  )
 L^{t,\o}_{ \tau  }
 = R^{t,\o} (\tau, \ga) .
 \eea
 Since $|\xi_n| \n \le \n  \int_t^T  \n |g^{t,\o}_s| ds    +    \Psi^{t,\o}_*  $, $\fa n \in \hN$ and since
 $\hE_\hP \big[ \int_t^T \n  |g^{t,\o}_s| ds    +    \Psi^{t,\o}_* \big]  \n < \n  \infty $ by \eqref{eq:xxx111},
 we can deduce from the dominated convergence theorem, \eqref{eq:da311} and \eqref{eq:da313} that
 $ \hE_\hP \big[ R^{t,\o} (\tau,\ga ) \big]   = \lmt{n \to \infty} \hE_\hP  [ \xi_n  ]
 \le \linf{n \to \infty} \hE_\hP  \big[ R^{t,\o} (\tau_n ,\ga )  \big] $.    \qed

 \fi

   \no {\bf Proof of Proposition \ref{prop_DPP2}:}
  Let  $0 \n \le \n   t  \n \le \n  s  \n \le \n  T $ and $  \o  \n \in \n  \O $. It suffices to show
  for a given  $\hP  \n \in \n  \cP(t,\o)$ that
      \bea
  \underset{\ga  \in \cT^t }{\inf} \,  \underset{\t  \in \cT^t  }{\sup} \,  \hE_\hP  \big[  R^{t,\o} (\t, \ga )  \big]
  \gs   \underset{\ga  \in \cT^t  }{\inf}
  \underset{\t  \in \cT^t  }{\sup}    \hE_\hP   \bigg[ \b1_{\{\tau \land \ga < s\}} R^{t,\o} (\tau, \ga )
  \n +  \n     \b1_{\{\tau \land \ga \ge s\}}
  \Big(  \ol{V}^{t,\o}_s \dn + \dn \int_t^s \n  g^{t,\o}_r   dr   \Big) \bigg]   .  \label{eq:uxu159}
   \eea
   Fix 
   $\e \n > \n 0$.
   There exists a $\wh{\ga} \= \wh{\ga} (\e) 
   \n \in \n  \cT^t$ such that
 \bea \label{eq:da321}
 \underset{\t  \in \cT^t  }{\sup} \,   \hE_\hP  \big[  R^{t,\o} (\t,\wh{\ga})    \big]
 \le  \underset{\ga  \in \cT^t }{\inf} \,
 \underset{\t  \in \cT^t  }{\sup} \,   \hE_\hP  \big[  R^{t,\o} (\t,\ga)    \big] + \e / 5  .
 \eea

 \no {\bf 1)}  Set   $\wh{\ga}' \n : = \n  \wh{\ga}  \n \vee \n  s  \n \in \n  \cT^t_s$.
 {\it In the first step, we   use a ``dense" countable subset
   of $\cT^s $ and Proposition \ref{prop_shift1} to show that  }
   \bea
     \ol{V}^{t,\o}_s \n + \n \int_t^s g^{t,\o}_r   dr
       \le  \underset{\t \in \cT^t_s }{\esssup} \,
 \hE_\hP \big[ R^{t,\o} (\t,\wh{\ga}') \big|\cF^t_s\big]        \n + \n \frac35 \e ,
     \q \pas    \label{eq:aa117}
     \eea

  As in the proof of  \cite[Proposition 4.1]{ROSVU} \big(see part (2a) and (2c) therein\big),
  we can construct a  dense  countable subset
 $\G$ of $\cT^s $ in  sense that for any $\d > 0$, $\z \in \cT^s$ and    $\wt{\hP} \in \fP_s$,
  \bea \label{eq:ad031}
  \hb{     $ \exists \, \{ \vs_n \}_{n \in \hN}  \n \subset \n  \G$ such that
  $ \lmtd{n \to \infty} \vs_n (\wh{\o}) \n = \n  \z (\wh{\o}) $, $ \fa \wh{\o} \n \in \n \O^s  $ and that
  $ \wt{\hP}  \{ \vs_n  \n \ne \n  \z_n \}  \n < \n  \d $, $ \fa  n \n \in \n \hN$, }
  \eea
  where  $   \z_n \n := \n  \sum^{\lfloor 2^n T \rfloor}_{i=\lfloor 2^n s \rfloor}
  \b1_{\{i 2^{-n} \le \z < (i+1) 2^{-n} \}}   \big( \frac{i+1}{ 2^n }   \n  \land \n  T \big)  \n \in \n  \cT^s $.

  \ss   Since   $\z ( \Pi^t_s )  \n \in \n  \cT^t_s   $    for any $\z  \n \in \n  \cT^s$
  by \eqref{eq:xax001},  it holds except  on   a $\hP-$null set  $\cN    $ that
 \bea  \label{eq:aa111}
 \hE_\hP \big[ R^{t,\o} \big( \z ( \Pi^t_s ) , \wh{\ga}'      \big)   \big|\cF^t_s\big]
  \n \le \n  \underset{\t \in \cT^t_s }{\esssup} \,
 \hE_\hP \big[ R^{t,\o} (\t,\wh{\ga}') \big|\cF^t_s\big]   ,  \q \fa \z \in \G .    \q
 \eea

 By Proposition \ref{prop_shift0} (2), $ \ga_{\wt{\o}} \df (\wh{\ga}')^{s,\wt{\o}} \ins \cT^s$.
 In light of      \eqref{eq:f475}, there exists  a  $\hP-$null set $\wt{\cN}   $
 such that for any $\wt{\o} \in \wt{\cN}^c $,
     \bea
     \hE_\hP \big[ R^{t,\o} \big( \z ( \Pi^t_s ) , \wh{\ga}'   \big) \big|\cF^t_s\big]  (\wt{\o})
    \n =  \n  \hE_{\hP^{s,\wt{\o}}}  \Big[ \big( R^{t,\o}  ( \z ( \Pi^t_s ) , \wh{\ga}'  ) \big)^{s,\wt{\o}}  \Big]
    \n  =  \n  \hE_{\hP^{s,\wt{\o}}}  \big[   R^{s,\o \otimes_t \wt{\o}}
     (\z,   \ga_{\wt{\o}}  )  \big]
     \n + \n \int_t^s g^{t,\o}_r (\wt{\o}) dr  ,  ~  \fa \z \in \G . \q
  \label{eq:aa113}
  \eea
    Here we used an analogy to \eqref{eq:ab027} that
   $    \big( R^{t,\o}  ( \z ( \Pi^t_s ) , \wh{\ga}'  ) \big)^{s,\wt{\o}}
   \n = \n       R^{s, \o \otimes_t  \wt{\o}} \big( \z , \ga_{\wt{\o}}  \big)
   \n + \n \int_t^s g^{t,\o}_r (\wt{\o}) dr $.


 By (P2),   there exist  an extension $(\O^t,\cF',\hP')$ of $(\O^t,\cF^t_T,\hP)$ and $\O' \in \cF'$ with $\hP'(\O') = 1$
 such that for any $\wt{\o} \in \O'$, $\hP^{s,   \wt{\o}} \in \cP (s,  \o  \otimes_t \wt{\o} ) $.
 Let $ \ol{\cN} $ be the $ \cF^t_T-$measurable set    containing $\cN \cup \wt{\cN}$ and with  $\hP(\ol{\cN})=0$.

 Now, fix   $\wt{\o} \in    \O' \cap \ol{\cN}^c   \in \cF'  $.
 There exists   a $ \z_{\wt{\o}} \in \cT^s $ such that
   \bea      \label{eq:cc133b}
   \underset{\z  \in \cT^s }{\sup} \, \hE_{\hP^{s,\wt{\o}}}
   \big[  R^{s, \o \otimes_t \wt{\o}}  (\z, \ga_{\wt{\o}}  )   \big]
  \n \le \n  \hE_{\hP^{s,\wt{\o}}}  \big[  R^{s, \o \otimes_t \wt{\o}}
  \big(\z_{\wt{\o}} , \ga_{\wt{\o}} \big)  \big]  \n + \n  \e / 5   \,   .
   \eea
  As $\hP^{s,\wt{\o}} \n \in \n   \cP(s,\o \otimes_t \wt{\o}) $, \eqref{eq:xxx111} shows that
 \bea \label{eq:uxu179}
  \hE_{\hP^{s,\wt{\o}}} \bigg[
     \int_s^T \n  \big| g^{s, \o \otimes_t \wt{\o}}_r \big| dr \n + \n \Psi^{s, \o \otimes_t \wt{\o}}_* \bigg]
      \n < \n  \infty   .
      \eea
 So for some $   \d_{\wt{\o}}  \n > \n  0$,
 \bea   \label{eq:bb435}
 \hb{ $ \dis  \hE_{\hP^{s,\wt{\o}}}
     \bigg[ \b1_A  \Big(  \int_s^T \big| g^{s, \o \otimes_t \wt{\o}}_r \big| dr
     \n + \n \Psi^{s, \o \otimes_t \wt{\o}}_*  \Big) \bigg] < \e / 5 $ \; for any $A \in \cF^s_T$
  with $ \hP^{s,\wt{\o}} (A) < \d_{\wt{\o}} $. }
  \eea

  Applying \eqref{eq:ad031} with $(\d,\z,\wt{\hP}) \= \big( \d_{\wt{\o}} ,\z_{\wt{\o}}, \hP^{s,\wt{\o}} \big)$, there exist
  $ \big\{ \vs^k_{\wt{\o}} \big\}_{k \in \hN}  \n \subset \n  \G$ such that
  $ \lmtd{k \to \infty} \vs^k_{\wt{\o}} (\wh{\o}) \n = \n  \z_{\wt{\o}}  (\wh{\o})  $, $ \wh{\o} \n \in \n \O^s  $ and that
  $  \hP^{s,\wt{\o}}   \{ \vs^k_{\wt{\o}} \ne \z^k_{\wt{\o}} \}  \n < \n  \d_{\wt{\o}} $, $\fa k \n \in \n \hN$,
  where  $   \z^k_{\wt{\o}} \n := \n  \sum^{\lfloor 2^k T \rfloor}_{i=\lfloor 2^k s \rfloor}
  \b1_{\{i 2^{-k} \le \z_{\wt{\o}}  < (i+1) 2^{-k} \}}   \big( \frac{i+1}{ 2^k }   \n  \land \n  T \big)  \n \in \n  \cT^s $.

     Given $k \n \in \n \hN$,   \eqref{eq:bb435} and \eqref{eq:ab015} imply   that
       \beas
    \hE_{\hP^{s,\wt{\o}}}  \Big[  \big| R^{s, \o \otimes_t \wt{\o}} \big( \z^k_{\wt{\o}} , \ga_{\wt{\o}} \big)
   \n - \n  R^{s, \o \otimes_t \wt{\o}}  \big( \vs^k_{\wt{\o}}   , \ga_{\wt{\o}} \big)   \big|  \Big]
    & \tn =  & \tn    \hE_{\hP^{s,\wt{\o}}}  \Big[ \b1_{\{ \z^k_{\wt{\o}} \ne \vs^k_{\wt{\o}} \}}
      \big| R^{s, \o \otimes_t \wt{\o}} \big( \z^k_{\wt{\o}} , \ga_{\wt{\o}} \big)
     \n - \n  R^{s, \o \otimes_t \wt{\o}}  \big( \vs^k_{\wt{\o}}   , \ga_{\wt{\o}} \big) \big|  \Big] \\
       & \tn  \le  & \tn   2 \hE_{\hP^{s,\wt{\o}}}
     \bigg[ \b1_{\{ \z^k_{\wt{\o}} \ne \vs^k_{\wt{\o}} \}} \Big(    \int_s^T \big| g^{s, \o \otimes_t \wt{\o}}_r \big| dr
     \n + \n \Psi^{s, \o \otimes_t \wt{\o}}_* \Big) \bigg] < \frac25 \e  ,
    \eeas
    which together with     \eqref{eq:aa111} and  \eqref{eq:aa113} shows that
        \beas
    \q      \hE_{\hP^{s,\wt{\o}}}  \big[  R^{s, \o \otimes_t \wt{\o}} \big( \z^k_{\wt{\o}} , \ga_{\wt{\o}} \big)  \big]
         \<   \hE_{\hP^{s,\wt{\o}}}  \big[  R^{s, \o \otimes_t \wt{\o}}  \big( \vs^k_{\wt{\o}}   , \ga_{\wt{\o}} \big)  \big]
         \n + \n  \frac25 \e  \ls   \underset{\t \in \cT^t_s }{\esssup} \,
 \hE_\hP \big[ R^{t,\o} (\t,\wh{\ga}') \big|\cF^t_s\big] (\wt{\o})
 \n - \n \int_t^s g^{t,\o}_r (\wt{\o}) dr   \n + \n    \frac25 \e   .
 \eeas
 As one can deduce from $\z_{\wt{\o}}  \n = \n  \lmtd{k \to \infty} \z^k_{\wt{\o}}$ and the continuity of  $L$ that
 \bhe \bea
R^{s,\o  \otimes_t  \wt{\o}} (\z_{\wt{\o}},\ga_{\wt{\o}} )
\ls  \lmt{k \to \infty} R^{s,\o  \otimes_t  \wt{\o}} (\z^k_{\wt{\o}},\ga_{\wt{\o}} ) , \label{eq:uxu181}
 \eea \ehe
   \eqref{eq:ab015},  \eqref{eq:uxu179},
  the dominated convergence theorem and    \eqref{eq:cc133b} imply  that
  \beas
  \q  \ol{V}^{t,\o}_s (  \wt{\o}  )  & \tn \dn = & \tn  \dn   \ol{V}_s ( \o  \n \otimes_t \n  \wt{\o})
   \n \le \n \underset{\z  \in \cT^s }{\sup} \, \hE_{\hP^{s,\wt{\o}}}
   \big[  R^{s, \o \otimes_t \wt{\o}}  (\z, \ga_{\wt{\o}}  )   \big]
  \n \le \n  \hE_{\hP^{s,\wt{\o}}}  \big[  R^{s, \o \otimes_t \wt{\o}}
  \big(\z_{\wt{\o}} , \ga_{\wt{\o}} \big)  \big]  \n + \n  \e / 5  \\
   & \tn \dn  = & \tn  \dn  \lmt{k \to \infty}
     \hE_{\hP^{s,\wt{\o}}}  \big[  R^{s, \o \otimes_t \wt{\o}} \big( \z^k_{\wt{\o}} , \ga_{\wt{\o}} \big)  \big]
     \n + \n  \e / 5
    \n \le \n    \underset{\t \in \cT^t_s }{\esssup} \,
 \hE_\hP \big[ R^{t,\o} (\t,\wh{\ga}') \big|\cF^t_s\big] (\wt{\o})
 \n - \n \int_t^s g^{t,\o}_r (\wt{\o}) dr  \n + \n \frac35 \e  , \q \fa \wt{\o}  \n \in  \n   \O' \cap \ol{\cN}^c  ,
 \eeas
 This shows   $\O'  \n \cap \n  \ol{\cN}^c  \n \subset \n  \ol{A}  \n := \n
  \Big\{   \ol{V}^{t,\o}_s \n + \n \int_t^s g^{t,\o}_r   dr
        \n \le \n   \underset{\t \in \cT^t_s }{\esssup} \,
 \hE_\hP \big[ R^{t,\o} (\t,\wh{\ga}') \big|\cF^t_s\big]  \n + \n \frac35 \e   \Big\}    $.
 As Remark \ref{rem_V_adapted} and  Proposition \ref{prop_shift0} (1) imply that
     $ \ol{V}^{t,\o}_s \n + \n \int_t^s g^{t,\o}_r   dr
      \n = \n  \big(  \ol{V}_s \n + \n \int_t^s g_r   dr \big)^{t,\o}  \n \in \n  \cF^t_s$,
   we see that $ \ol{A}  \n \in \n  \cF^t_s $ and thus
   $   \hP (\ol{A}) \n = \n  \hP'(\ol{A})  \n \ge \n  \hP' \big(\O' \cap  \ol{\cN}^c \big)  \n = \n 1 $.
   Therefore,      \eqref{eq:aa117} holds.

   Moreover, one can find a sequence  $\{\tau_n\}_{n \in \hN} $     in $ \cT^t_s $ such that
 \bhe \bea \label{eq:uxu175}
\underset{\t \in \cT^t_s }{\esssup} \, \hE_\hP \big[R^{t,\o} (\tau , \wh{\ga}')\big|\cF^t_s\big]
= \lmtu{n \to \infty } \hE_\hP \big[R^{t,\o} (\tau_n , \wh{\ga}')\big|\cF^t_s\big]  , \q \pas
 \eea \ehe

   \no  {\bf 2)}
       Next, let $\tau \ins \cT^t$ and $n \ins \hN $. Since
        \bhe \bea \label{eq:uxu177}
          \ol{\t}_n \df  \b1_{\{\t \land \wh{\ga} <s\}} \t
   \+  \b1_{\{\t \land \wh{\ga}  \ge s\}} \tau_n
   \eea \ehe
   defines an $\bF^t-$stopping time,
  \eqref{eq:aa117} and \eqref{eq:uxu170} show that
      \bea
  && \hspace{-1.5cm} \hE_\hP   \bigg[ \b1_{\{\tau \land \wh{\ga} < s\}} R^{t,\o} (\tau, \wh{\ga})
  \n +  \n     \b1_{\{\tau \land \wh{\ga} \ge s\}}
  \Big(  \ol{V}^{t,\o}_s \n + \n \int_t^s g^{t,\o}_r   dr   \Big) \bigg] \nonumber \\
  && \le      \hE_\hP    \Big[ \b1_{\{\tau \land \wh{\ga} < s\}} R^{t,\o} ( \ol{\tau}_n , \wh{\ga})
  \n +  \n    \b1_{A_n \cap \{\tau \land \wh{\ga} \ge s\}}
   \Big( \hE_\hP \big[R^{t,\o} (\tau_n ,  \wh{\ga}')\big|\cF^t_s\big] \+ \frac45 \e \Big) \Big]
  \+ \a_n  ,
  \qq  \label{eq:da323}
  \eea
  where $A_n \n := \n \Big\{  \underset{\t \in \cT^t_s }{\esssup} \, \hE_\hP \big[R^{t,\o} (\tau , \wh{\ga}')\big|\cF^t_s\big]
 \n < \n  \hE_\hP \big[R^{t,\o} (\tau_n , \wh{\ga}')\big|\cF^t_s\big] \+ \e / 5 \Big\} \ins \cF^t_s $ and
  $\a_n \df  \hE_\hP    \Big[   \b1_{A^c_n}
  \big(   \int_t^T  | g^{t,\o}_r  |   dr  \n + \n \Psi^{t,\o}_*    \big) \Big]  $.

  Also, we can deduce from \eqref{eq:ab015} that
  \beas
 \q   &&  \hspace{-1.2cm}  \hE_\hP    \Big[     \b1_{A_n \cap \{\tau \land \wh{\ga} \ge s\}}
   \hE_\hP \big[R^{t,\o} (\tau_n ,  \wh{\ga}')\big|\cF^t_s\big]  \Big]
   =   \hE_\hP    \Big[
   \hE_\hP \big[   \b1_{A_n \cap \{\tau \land \wh{\ga} \ge s\}}
    R^{t,\o} (\tau_n ,  \wh{\ga}')\big|\cF^t_s\big]  \Big]
       =       \hE_\hP \big[   \b1_{A_n \cap \{\tau \land \wh{\ga} \ge s\}}
    R^{t,\o} (\tau_n ,  \wh{\ga}) \big]  \\
  &&   =  \n    \hE_\hP    \big[      \b1_{\{\tau \land \wh{\ga} \ge s\}}
    R^{t,\o} (\tau_n ,  \wh{\ga})
     \n -  \n        \b1_{A^c_n \cap \{\tau \land \wh{\ga} \ge s\}}
    R^{t,\o} (\tau_n ,  \wh{\ga})   \big]
     \n \le  \n    \hE_\hP    \big[      \b1_{\{\tau \land \wh{\ga} \ge s\}}
    R^{t,\o} (\ol{\tau}_n ,  \wh{\ga})        \big] \n + \n \a_n ,
   \eeas
   which together with \eqref{eq:da323} and \eqref{eq:da321}    leads to that
         \beas
    \hE_\hP   \bigg[ \b1_{\{\tau \land \wh{\ga} < s\}} R^{t,\o} (\tau, \wh{\ga})
  \n +  \n     \b1_{\{\tau \land \wh{\ga} \ge s\}}
  \Big(  \ol{V}^{t,\o}_s \n + \n \int_t^s g^{t,\o}_r   dr   \Big) \bigg]
  & \tn \le  & \tn      \hE_\hP    \big[   R^{t,\o} ( \ol{\tau}_n , \wh{\ga})   \big]
  \+ 2 \a_n  \+ \frac45 \e
    \le   \underset{\t  \in \cT^t  }{\sup} \,   \hE_\hP  \big[  R^{t,\o} (\t,\wh{\ga})    \big]
  \+ 2 \a_n  \+ \frac45 \e  \\
     & \tn  \le  & \tn  \underset{\ga  \in \cT^t }{\inf} \,   \underset{\t  \in \cT^t  }{\sup} \,
    \hE_\hP  \big[  R^{t,\o} (\t, \ga)    \big]   + 2 \a_n  +   \e   .
  \eeas

  Since  $ \lmtu{n \to \infty } \hP (A_n)  \n = \n  1 $ by \eqref{eq:uxu175},
    we see from \eqref{eq:uxu183} and the  dominated convergence theorem   that
  $\lmtd{n \to \infty} \a_n \= 0$ and thus
           \bea  \label{eq:da351}
     \hE_\hP   \bigg[ \b1_{\{\tau \land \wh{\ga} < s\}} R^{t,\o} (\tau, \wh{\ga})
  \n +  \n     \b1_{\{\tau \land \wh{\ga} \ge s\}}
  \Big(  \ol{V}^{t,\o}_s \n + \n \int_t^s g^{t,\o}_r   dr   \Big) \bigg]
  \le  \underset{\ga  \in \cT^t }{\inf} \,   \underset{\t  \in \cT^t  }{\sup} \,
    \hE_\hP  \big[  R^{t,\o} (\t, \ga )    \big]    +   \e  , \q \fa \tau \ins \cT^t .
  \eea
    Taking supremum over $\t \n \in \n  \cT^t$ on the left-hand-side 
    \if{0}
      \bea
  \q   && \hspace{-1.5cm} \underset{\ga  \in \cT^t  }{\inf}
  \underset{\t  \in \cT^t  }{\sup}    \hE_\hP   \bigg[ \b1_{\{\tau \land \ga < s\}} R^{t,\o} (\tau, \ga )
  \n +  \n     \b1_{\{\tau \land \ga \ge s\}}
  \Big(  \ol{V}^{t,\o}_s \dn + \dn \int_t^s \n  g^{t,\o}_r   dr   \Big) \bigg]
   \n \le  \n  \underset{\t  \in \cT^t  }{\sup} \,
   \hE_\hP   \bigg[ \b1_{\{\tau \land \wh{\ga} < s\}} R^{t,\o} (\tau, \wh{\ga})
  \n +  \n     \b1_{\{\tau \land \wh{\ga} \ge s\}}
  \Big(  \ol{V}^{t,\o}_s \dn + \dn \int_t^s \n  g^{t,\o}_r   dr   \Big) \bigg] \nonumber \\
 && \le  \underset{\ga  \in \cT^t }{\inf} \,   \underset{\t  \in \cT^t  }{\sup} \,
    \hE_\hP  \big[  R^{t,\o} (\t, \ga )    \big]    +   \e .  \label{eq:uxu160}
   \eea
 \fi
  and then  letting $\e \n \to \n  0$ lead to \eqref{eq:uxu159}.   \qed

   \no {\bf Proof of Proposition \ref{prop_grid}:}
 Let $n \in \hN$, $t \in [0,T]$, $\a > 0$ and $ \o \in O^t_\a (\bz)$.

 We fix $\hP  \n \in \n  \cP(t,\o)$ and $\ga,  \tau  \n \in \n  \cT^t$.
 Set $ \{t^n_i\}^{2^n}_{i=0}  $ as in \eqref{eq:uxu185} and define
  $\dis \tau_n  \n : = \n  \b1_{\{\tau = t\}} t  \n   + \n  \sum^{2^n}_{i=1}
 \b1_{ \{ t^n_{i-1} < \tau \le t^n_i  \} } t^n_i   \n   \in \n  \cT^t (n) $.
  One can deduce that
 \bea
 R^{t,\o} (\tau,\ga) \n -  \n  R^{t,\o}  (\tau_n ,\ga)
 & \tn \= & \tn   \- \int_{\tau  \land \ga}^{\tau_n \land \ga} \n g^{t,\o}_r dr
 \+ \b1_{\{ \tau \le \ga  \}} \big(  L^{t,\o}_\tau \-  \b1_{\{ \tau_n \le \ga  \}} L^{t,\o}_{\tau_n}
 \- \b1_{\{ \ga < \tau_n   \}}    U^{t,\o}_\ga   \big)
 \+ \b1_{\{ \ga < \tau   \}}  ( U^{t,\o}_\ga \- U^{t,\o}_\ga )  \nonumber  \\
 & \tn \= & \tn  - \n \int_{\tau  \land \ga}^{\tau_n \land \ga} \n g^{t,\o}_r dr
 \+ \sum^{2^n}_{i = 1}  \n  \Big( \b1_{\{t^n_{i-1} < \tau \le t^n_i \le \ga  \}}
 \big( L^{t,\o}_\tau \- L^{t,\o}_{t^n_i}       \big)
 \n + \n \b1_{\{t^n_{i-1} < \tau \le \ga < t^n_i  \}} \big( L^{t,\o}_\tau \- U^{t,\o}_\ga  ) \Big) .
 \qq \qq  \label{eq:da361}
 \eea

 Given $i = 1, \cds, 2^n$,   \eqref{eq:aa211} shows that for any $ \wt{\o} \in \{t^n_{i-1} < \tau \le t^n_i \le \ga  \} $
 \bea
 \big| L^{t,\o}_\tau (\wt{\o})  \n - \n  L^{t,\o}_{t^n_i} (\wt{\o})   \big|
  & \tn  \dn =  & \tn \dn   \big|  L \big( \t  (  \wt{\o}), \o  \n \otimes_t \n  \wt{\o}  \big) \n  - \n
   L \big(  t^n_i , \o  \n \otimes_t \n  \wt{\o}  \big)  \big|
   \n \le \n   \rho_0 \Big(  \big( t^n_i   \n - \n  \t  (  \wt{\o})  \big)
   \n + \n  \underset{r \in  [0 , T  ] }{\sup}
  \big| (\o  \n \otimes_t \n  \wt{\o})   \big( r  \n \land \n  \t  (  \wt{\o}) \big)
   \n - \n  (\o  \n \otimes_t \n  \wt{\o}) ( r  \n \land \n  t^n_i )   \big|  \Big) \nonumber  \\
 & \tn  \dn \le & \tn \dn   \rho_0 \Big(  2^{-n}
  \n + \n  \underset{r \in  [\t  ( \wt{\o}) , t^n_i ] }{\sup}
  \big| \wt{\o}   (r)  \n - \n  \wt{\o} ( \t  (  \wt{\o}) )   \big|  \Big)
   \n \le  \n  \rho_0 \Big(  2^{-n}  \n + \n  \underset{ \tau (\wt{\o}) \le  r \le
  ( \tau (\wt{\o}) + 2^{-n} ) \land T }{\sup}
  \big| B^t_r (\wt{\o})  \n - \n  B^t_\tau (\wt{\o})   \big| \Big)  . \qq   \qq    \label{eq:da363}
 \eea
    Similarly, it holds for any $ \wt{\o} \in \{t^n_{i-1} < \tau \le \ga < t^n_i  \} $ that
 \bea
 \hspace{-5mm}
  \big|  U^{t,\o}_\tau   \n - \n  U^{t,\o}_\ga    \big| ( \wt{\o} )
   \ls     \rho_0 \Big(  \big( \ga (\wt{\o})   \n - \n  \t  (  \wt{\o})  \big)
   \n + \n  \underset{r \in  [\t  ( \wt{\o}) , \ga (\wt{\o}) ] }{\sup}
  \big| \wt{\o}   (r)  \n - \n  \wt{\o} ( \t  (  \wt{\o}) )   \big|  \Big)
       \ls   \rho_0 \Big(  2^{-n}  \n + \n  \underset{ \tau (\wt{\o}) \le  r \le
  ( \tau (\wt{\o}) + 2^{-n} ) \land T }{\sup}
  \big| B^t_r (\wt{\o})  \n - \n  B^t_\tau (\wt{\o})   \big| \Big)    . \q  \label{eq:da367}
 \eea
 Moreover, another analogy to \eqref{eq:da363} shows that for any $(s,\wt{\o}) \in [t,T] \times \O^t$
\bea
   \big| g^{t,\o}_s (\wt{\o})  \n - \n  g_t (\o) \big|
  \n \le  \n   \big|  g \big( s , \o \otimes_t \wt{\o}  \big)  \n - \n    g  (  t , \o  )  \big|
   \n \le  \n  \rho_0 \Big( s   \n - \n  t  \n + \n  \underset{r \in [t,s] }{\sup}
  \big| \wt{\o} (r)     \big| \Big)
    \n \le \n   \rho_0  \Big( T   \n - \n  t  \n + \n  \underset{r \in [t , T]   }{\sup}
  \big| B^t_r (\wt{\o})  \n - \n  B^t_t (\wt{\o})    \big| \Big)  ,  \q \label{eq:da369}
\eea
 where we used the fact that   $B^t_t = 0 $ in the last inequality.
 Plugging  \eqref{eq:da363}$-$\eqref{eq:da369}  back into \eqref{eq:da361} leads to that
 \beas
  R^{t,\o} (\tau,\ga) \n  - \n   R^{t,\o}  (\tau_n ,\ga)
  \n \le  \n  2^{-n} \bigg[ | g_t (\o)|
   \n + \n  \rho_0 \Big( T   \n - \n  t  \n + \n \underset{r \in [t , T]   }{\sup}
  \big| B^t_r (\wt{\o})  \n - \n  B^t_t (\wt{\o})    \big|  \Big)  \bigg]
  \n + \n   \rho_0 \Big(  2^{-n}  \n + \n  \underset{r \in [\tau ,
  ( \tau   + 2^{-n} ) \land T] }{\sup}
   | B^t_r   \n - \n  B^t_\tau      | \Big)   .
\eeas

 Taking expectation $\hE_\hP [~]$, we see from \eqref{eq:aa213b} that
 \beas
 \hE_\hP \big[  R^{t,\o} (\tau,\ga) \big]
   \ls  \hE_\hP \big[  R^{t,\o} (\tau_n,\ga) \big]
    \n + \n  I^n_\a
      \ls   \underset{\tau' \in \cT^t (n)}{\sup} \hE_\hP \big[  R^{t,\o} (\tau',\ga) \big]
      \n + \n  I^n_\a  ,
 \eeas
 where $I^n_\a \df  \rho_\a  (  2^{-n}  )  \n + \n   2^{-n} \big( | g_t (\o)| \n + \n \rho_\a  (T \- t)  \big) $.
 Taking supremum over $\tau \in \cT^t$ on the left-hand-side yields that
 \beas
  \underset{\tau \in \cT^t  }{\sup} \hE_\hP \big[  R^{t,\o} (\tau,\ga) \big]
    \le  \underset{\tau \in \cT^t (n)}{\sup} \hE_\hP \big[  R^{t,\o} (\tau,\ga) \big]
      \n + \n  I^n_\a  .
 \eeas
 Eventually, taking  infimum over $ \ga  \n \in \n  \cT^t$ and   $ \hP \in \cP (t,\o) $ leads to
 \eqref{eq:da365}.   \qed

     \no {\bf Proof of Proposition \ref{prop_conti_V}:}
    Fix $n  \n \in \n  \hN \cup \{\infty\}$, $ \o   \n \in \n   \O$ and set $\a  \n := \n  1  \n + \n  \|\o\|_{0,T}$.
  Let  $0  \n \le \n  t  \n < \n  s  \n \le \n  T$ such that
  $\d_{t,s} \n := \n (s \n - \n t) \vee
  \underset{t\le r < r' \le s  }{\sup} \big| \o(r' ) \n - \n  \o(  r) \big|
    \n \le \n  T $.

   \no {\bf 1a)} {\it We first utilize Proposition \ref{prop_DPP} and \eqref{eq:aa213b} to show that}
    \bea  \label{eq:db111}
   V^n_{t}(\o) \-  V^n_{s}(\o) \n \le \n (s \n - \n t) \underset{r \in [0,T]}{\sup} | g_r (\o)  |
   \n + \n (2 \n + \n s \n - \n t) \rho_\a (\d_{t,s})  .
   \eea

     Let  $\hP \n \in \n \cP(t,\o)$.
 Applying \eqref{eq:bb013} and taking $\ga = s  $   show   that
 \bea
  V^n_{t}(\o) -  V^n_{s}(\o)
  & \tn   \le    & \tn
  \underset{\t \in \cT^t (n)}{\sup} \,   \hE_\hP \bigg[ \b1_{\{\tau   < s\}}   R^{t,\o} (\tau, s )
   + \b1_{\{\tau   \ge s\}} \Big( (V^n_s)^{t,\o} + \int_t^s g^{t,\o}_r dr \Big)     \bigg]   -  V^n_{s}(\o)
     \nonumber   \\
     & \tn   =    & \tn
  \underset{\t \in \cT^t (n)}{\sup} \,   \hE_\hP \bigg[ \b1_{\{\tau   < s\}}   L^{t,\o}_\tau
   + \b1_{\{\tau   \ge s\}}  (V^n_s)^{t,\o} -  V^n_{s}(\o) + \int_t^{\tau \land s} g^{t,\o}_r dr  \bigg]  .
    \label{eq:bb439}
 \eea

 Then,  let $ \tau \in \cT^t (n) $.  For any $\wt{\o} \in \{\tau < s \} $,    \eqref{eq:aa211} implies  that
 \bea
   \big|   L^{t,\o}_\tau (\wt{\o}) \n - \n    L^{t,\o}_s (\wt{\o})  \big|
  & \tn \dn = & \tn \dn  \big|   L  \big( \tau  (\wt{\o}),   \o \otimes_t \wt{\o} \big)
 -  L  \big( s,   \o \otimes_t \wt{\o} \big)  \big|
 \n \le \n  \rho_0 \Big( (s  \n - \n t)  \n + \n   \underset{r \in [t,T]}{\sup}
  \big| \wt{\o} \big( r  \n \land \n   \tau  (\wt{\o}) \big)
    \n - \n  \wt{\o}  ( r  \n \land \n  s  ) \big|  \Big) \nonumber \\
 & \tn \dn  \le  & \tn \dn
      \rho_0 \Big( (s  \n - \n t)  \n + \n
      \underset{r \in  [\tau  (\wt{\o}), s]   }{\sup}
   \big|  \wt{\o} (r) \n - \n \wt{\o} (\tau  (\wt{\o}))  \big|  \Big)
     \n   \le   \n   \rho_0 \Big( (s  \n - \n t)  \n + \n
     \underset{r \in  [\tau  (\wt{\o}), (\tau  (\wt{\o}) +s-t) \land T]   }{\sup}
  \big|  B^t_r (\wt{\o}) \n - \n  B^t_\tau ( \wt{\o})  \big|  \Big)
   .  \qq \q   \label{eq:af017}
 \eea
 Similarly, using \eqref{eq:aa211} again  and
 applying \eqref{eq:bb421} with $\eta = g_t \in \cF_t $ yields that for any $\wt{\o} \in \O^t$
  \bea
 \bigg|    \int_t^{  \tau  (\wt{\o}) \land s} \n  g^{t,\o}_r (\wt{\o}) dr    \bigg|
 & \tn \le   & \tn   \int_t^s  \n \big| g^{t,\o}_r (\wt{\o}) \big| dr
  \n \le  \n      \int_t^s \big( \big| g^{t,\o}_t (\wt{\o}) \big|
  \n +  \n     \big| g^{t,\o}_r (\wt{\o})  \n - \n  g^{t,\o}_t (\wt{\o}) \big| \big) dr \nonumber \\
   & \tn  \le   & \tn       \int_t^s \Big(  | g_t (\o)  |
  \n + \n   \rho_0 \Big( (s  \n - \n t)  \n + \n   \underset{r \in [t , s]   }{\sup}
  \big|  B^{t}_r (\wt{\o}) \n - \n  B^{t}_t ( \wt{\o})  \big|  \Big) \Big) dr     .  \label{eq:af025}
 \eea
 Also,   \eqref{eq:aa213} shows that  for any $ \wt{\o}  \n \in \n  \O^{t} $
  \bea
  \big|  V^n_{s}(\o)  \n   -  \n   (V^n_s)^{t,\o}  (\wt{\o}) \big|
     & \tn  \dn   =  & \tn  \dn    \big|  V^n_{s}(\o)   \n    -  \n   V^n (s, \o \otimes_{t} \wt{\o} ) \big|
      \n  \le \n  \rho_1  \big(  \| \o  \n - \n  \o \otimes_{t} \wt{\o} \|_{0,s} \big)
      \n = \n  \rho_1 \Big( \, \underset{r \in [t, s]}{\sup} \big| \o(r) \n  - \n \o(t) \n  - \n  \wt{\o}(r)
           \big| \Big)    \nonumber  \\
    & \tn  \dn  \le & \tn  \dn  \rho_1 \Big( \,
         \underset{r \in [t, s]}{\sup} \big|  \o(r) - \o(t)     \big|
     \n  +  \underset{r \in [t, s]}{\sup} \big| \wt{\o}(r)    \big| \Big)
     \n  \le \rho_1 \Big( \, \d_{t,s} +   \underset{r \in [t , (t+\d_{t,s}) \land T ]    }{\sup}
  \big|  B^{t}_r (\wt{\o}) \n - \n  B^{t}_t ( \wt{\o})  \big|    \Big)  .
     \qq \qq  \q  \label{eq:bb441}
 \eea

 Since $ \|\o\|_{0,t} \le \|\o\|_{0,T} < \a$,
 we can deduce from \eqref{eq:af017}, \eqref{eq:af025}, \eqref{eq:uxu170}, \eqref{eq:aa213b} and \eqref{eq:bb441}  that
   \beas
 && \hspace{-0.7cm}  \hE_\hP \bigg[ \b1_{\{\tau   < s\}}   L^{t,\o}_\tau
   + \b1_{\{\tau   \ge s\}}  (V^n_s)^{t,\o} -  V^n_{s}(\o) + \int_t^{\tau \land s} g^{t,\o}_r dr  \bigg]
   \n - \n (s \n - \n t)  | g_t (\o)  | \\
 & &  \le  \n
      \hE_\hP \bigg[ \b1_{\{\tau   < s\}}   L^{t,\o}_s
   \n + \n  \b1_{\{\tau   \ge s\}}  (V^n_s)^{t,\o}  \n - \n   V^n_{s}(\o)
    \n + \n   \rho_1 \Big( (s  \n - \n t)  \n + \n
     \underset{r \in  [\tau   , (\tau    +s-t) \land T]   }{\sup}
  \big|  B^t_r   \n - \n  B^t_\tau    \big|  \Big)
    \n + \n   (  s \n - \n t) \rho_1 \Big( (s  \n - \n t)  \n + \n  \underset{r \in [t , s]   }{\sup}
  \big|  B^t_r   \n - \n  B^t_t    \big|  \Big)    \bigg]  \\
    & &  \le  \n
      \hE_\hP   \big[   (V^n_s)^{t,\o} \n - \n   V^n_{s}(\o) \big]
   \n + \n   (1 \n + \n s \n - \n t) \rho_\a (s-t)
     \n \le \n  (2 \n + \n s \n - \n t) \rho_\a (\d_{t,s}) .
 \eeas
 Taking supremum over $ \tau \in \cT^t (n) $ on the left-hand-side, we obtain \eqref{eq:db111} from \eqref{eq:bb439}.

  \no {\bf 1b)} {\it Next, we show that for $\ol{V}$ the inequality \eqref{eq:db111} can be strengthened as }
  \bea   \label{eq:db114}
  \big| \ol{V}_{s}(\o) \n - \n  \ol{V}_{t}(\o) \big|    \le
  (s \n - \n t) \underset{r \in [0,T]}{\sup} | g_r (\o)  | \n +  \n  (2 \n + \n s \n - \n t) \rho_\a ( \d_{t,s} )     .
  \eea

   Fix $\e \n > \n  0$.    We can find a
  $ \hP  \=  \hP (\e)   \n \in \n  \cP(t,\o)  $ such that
 $  \ol{V}_{t}(\o) \n + \n  \e / 2
   \n \ge \n   \underset{\ga \in \cT^t }{\inf} \underset{\t  \in \cT^t }{\sup} \,   \hE_\hP
 \big[    R^{t,\o} (\t, \ga )       \big]   $.
 By  \eqref{eq:da351},  there exists some $\wh{\ga} \= \wh{\ga} (\e)   \n \in \n  \cT^t$ such that
           \beas
  \hE_\hP   \bigg[ \b1_{\{\tau \land \wh{\ga} < s\}} R^{t,\o} (\tau, \wh{\ga})
  \n +  \n     \b1_{\{\tau \land \wh{\ga} \ge s\}}
  \Big(  \ol{V}^{t,\o}_s \n + \n \int_t^s g^{t,\o}_r   dr   \Big) \bigg]
  \le  \underset{\ga  \in \cT^t }{\inf} \,   \underset{\t  \in \cT^t  }{\sup} \,
    \hE_\hP  \big[  R^{t,\o} (\t, \ga )    \big]    +   \e / 2 , \q \fa \tau \ins \cT^t  .
  \eeas
 In particular,  taking $\tau  \n = \n  s$ on the left-hand-side gives that
 \bea
  \ol{V}_{t}(\o)    \+    \e
 \gs    \hE_\hP  \bigg[ \b1_{\{  \wh{\ga} < s\}}  R^{t,\o} (s, \wh{\ga})   \n +  \n
   \b1_{\{  \wh{\ga}  \ge s\}} \Big( \ol{V}^{t,\o}_s \+ \int_t^s \n g^{t,\o}_r dr   \Big) \bigg] 
    \=   \hE_\hP \bigg[ \int_t^{\wh{\ga} \land s} g^{t,\o}_r dr  \n +  \n
 \b1_{\{  \wh{\ga} < s\}}  U^{t,\o}_{  \wh{\ga} }   \n +  \n
 \b1_{\{  \wh{\ga}  \ge s\}}   \ol{V}^{t,\o}_s      \bigg]   . \q    \label{eq:af021}
 \eea

 An analogy to \eqref{eq:af017} and \eqref{eq:af025} shows that
  \beas
  \big|  U^{t,\o}_{  \wh{\ga} } (\wt{\o}) \n  -  \n    U^{t,\o}_s  (\wt{\o})  \big|
   & \tn  \le  & \tn   \rho_0 \Big( (s  \n - \n t)  \n + \n
     \underset{r \in  [\wh{\ga}  (\wt{\o}), (\wh{\ga}  (\wt{\o}) +s-t) \land T]   }{\sup}
  \big|  B^t_r (\wt{\o}) \n - \n  B^t_{\wh{\ga}} ( \wt{\o})  \big|  \Big) ,
  \q \fa \wt{\o} \in \{  \wh{\ga}  < s\} \q \hb{and} \\
  \bigg|    \int_t^{  \wh{\ga}  (\wt{\o}) \land s} \n  g^{t,\o}_r (\wt{\o}) dr    \bigg|
  & \tn  \le   & \tn       (s \n - \n t) \Big(  | g_t (\o)  |
  \n + \n   \rho_0 \Big( (s  \n - \n t)  \n + \n  \underset{r \in [t , s ]  }{\sup}
  \big|  B^t_r (\wt{\o}) \n - \n  B^t_t ( \wt{\o})  \big|  \Big) \Big)      , \q   \fa \wt{\o} \in \O^t .
 \eeas
 As $ \|\o\|_{0,t} \n \le \n  \|\o\|_{0,T}  \n < \n  \a$,
 plugging them back to \eqref{eq:af021} and applying \eqref{eq:bb441} with $n  \n = \n  \infty$,
 we can deduce from \eqref{eq:aa213b} and \eqref{eq:uxu170}   that
  \beas
&& \hspace{-0.8cm} \ol{V}_{t}(\o) \n- \n  \ol{V}_s (\o)   \n + \n  \e \n + \n  (s \n - \n t)  | g_t (\o)  |
 \n  \ge   \n      \hE_\hP \Big[ \b1_{\{  \wh{\ga} < s\}}  U^{t,\o}_s
    \n +  \n    \b1_{\{  \wh{\ga}  \ge s\}}   \ol{V}^{t,\o}_s  \n- \n  \ol{V}_s (\o) \Big]
     \n - \n  (1 \n + \n s \n - \n t) \rho_\a (s-t)
        \\
  &&    \ge    \n     \hE_\hP \big[    \ol{V}^{t,\o}_s  \n- \n  \ol{V}_s (\o) \big]
   \n - \n  (1 \n + \n s \n - \n t) \rho_\a (s-t)
  \n    \ge    \n   -   (2 \n + \n s \n - \n t) \rho_\a (\d_{t,s})    .
 \eeas
 Letting $\e \to 0$ and taking \eqref{eq:db111} with $n \= \infty$  yield   \eqref{eq:db114}.

   Since $ \lmtd{t \nearrow s} \d_{t,s} = \lmtd{s \searrow t} \d_{t,s} = 0$, we can deduce from
   \eqref{eq:db111} and \eqref{eq:db114} that  each path  of   $V^n$ is both  left-upper-semicontinuous and right-lower-semicontinuous, in particular, each path of    $\ol{V}$ is continuous.

  \no {\bf 2)}   Given $(t,\o) \ins [0,T] \ti  \O$, Remark \ref{rem_V_adapted},
   Proposition \ref{prop_shift0} (4) and Part 1 show  that
 $ \ol{V}^{t,\o} $ is an $\bF^t-$adapted process with all continuous paths.
    For any $\hP \ins \cP(t,\o)$,  \eqref{eq:uxu170} and \eqref{eq:xxx111} imply     that
 $ \hE_\hP \Big[\ol{V}^{t,\o}_* \Big]
  \n \le \n  \hE_\hP \big[\Psi^{t,\o}_* \big]  \n < \n  \infty  $.
  So $\ol{V}^{t,\o} \ins \hS  ( \bF^t,   \hP)$.   \qed

  \subsection{Proofs of the results in Section \ref{sec:RDG}}\label{pfsec4}

        \no {\bf Proof of \eqref{eq:uxu173}:}
    Fix $n  \n \in \n  \hN \cup \{\infty\}$ and $\tau   \n \in \n  \cT$. We let
    $(t,\o)  \n \in \n  [0,T]  \n \times \n  \O$ and $ \hP  \n \in \n  \cP(t,\o) $.
    Since $ V^n_\tau \in \cF_T $ and $ \int_0^\tau g_r dr  \in \cF_T  $ by Remark \ref{rem_V_adapted},
    Proposition \ref{prop_shift0} (1)
    shows that both $ \big( V^n_\tau \big)^{t,\o}   $ and
    $ \big( \int_0^\tau g_r dr \big)^{t,\o} $ belong to $ \cF^t_T $.

    \no {\bf 1)}        If $\wh{t}  \n := \n  \tau (\o   )     \n \le \n  t$,
   Proposition \ref{prop_shift0} (3) 
   shows that $\tau ( \o  \n \otimes_t \n  \O^t )   \n \equiv \n  \wh{t}    $.
    Applying \eqref{eq:bb421}   to $\eta  \n = \n  V^n_{\wh{t}}  \n \in \n  \cF_{\wh{t}}  \n  \subset \n  \cF_t  $
    and to $\eta \n = \n  \int_0^{\wh{t}} \n  g_r dr   \n \in \n  \cF_{\wh{t}}   \n \subset \n  \cF_t  $ yields that
    for any $ \wt{\o}  \n \in \n  \O^t $
   \bea \label{eq:dc311}
     \big( V^n_{\tau  } \big)^{t,\o} ( \wt{\o} )
   \n = \n  V^n \big(     \tau  (\o  \n \otimes_t \n   \wt{\o})   , \o  \n \otimes_t \n  \wt{\o} \big)
    \n = \n  V^n \big(  \wh{t},      \o  \n \otimes_t \n  \wt{\o} \big)
    \n = \n  V^n \big(  \wh{t},      \o \big)  ,
    \eea
    and $ \big( \int_0^{\tau  }  \n  g_r dr  \big)^{t,\o} ( \wt{\o} )
    \n = \n  \int_0^{\tau  (\o    \otimes_t     \wt{\o})} \n  g_r (\o  \n \otimes_t \n  \wt{\o}) dr
     \n = \n  \int_0^{\wh{t}} g_r (\o  \n \otimes_t \n  \wt{\o}) dr
     \= \int_0^{\wh{t}} g_r (\o) dr $.
     Both only depend on $\o$.

    \no {\bf 2)}   Next, suppose that   $  \tau   \n > \n  t$.
    Proposition \ref{prop_shift0} (3)  also  shows that
     $ \tau (\o    \otimes_t   \wt{\o})     \>     t   $, $ \fa \wt{\o} \ins  \O^t $
        and that $\z \df \tau^{t,\o}   $ is a $\cT^t-$stopping time.
    It follows that
    $     \big( V^n_{\tau  } \big)^{t,\o} ( \wt{\o} )
    \n = \n  V^n \big(     \tau  (\o  \n \otimes_t \n   \wt{\o})   , \o  \n \otimes_t \n  \wt{\o} \big)
    \n = \n  V^n \big(  \tau^{t,\o} (\wt{\o}),      \o  \n \otimes_t \n  \wt{\o} \big)
    \n = \n  (V^n)^{t,\o} \big(  \z (\wt{\o}) , \wt{\o} \big)  $, $ \fa \wt{\o}  \n \in \n  \O^t $.
   By the first equality of \eqref{eq:xax141},   we also have
  $  \big( \n \int_0^\tau \n g_r dr  \big)^{t,\o} ( \wt{\o} )
     \= \n  \int_0^{\tau  (\o    \otimes_t     \wt{\o})} \n  g_r (\o   \n \otimes_t \n    \wt{\o}) dr
     \= \n  \int_0^t \n  g_r (\o   ) dr  \n + \n
     \int_t^{\z ( \wt{\o})} \n  g^{t,\o}_r (    \wt{\o}) dr $.
  Then \eqref{eq:uxu170}  and \eqref{eq:xxx111} imply that
   \beas
  \hE_\hP \bigg[ \big|( V^n_\tau    )^{t,\o} \big| \n + \n
   \bigg| \Big( \int_0^\tau  \n  g_r dr   \Big)^{t,\o} \bigg| \bigg]
   \n \le \n  \hE_\hP \bigg[ \big|( V^n    )^{t,\o}_\z \big|
   \n + \dn  \int_t^{\z  } \n \big| g^{t,\o}_r \big|  dr  \bigg]
   \n + \dn  \int_0^t \n  |g_r (\o   )| dr
    \n \le \n    \hE_\hP \bigg[  \Psi^{t,\o}_*  \n + \dn     \int_t^T \n \big| g^{t,\o}_r \big|  dr  \bigg]
   \n + \dn  \int_0^t \n  |g_r (\o   )| dr   \n < \n  \infty .
   \eeas

 \no {\bf Proof of   Theorem \ref{thm_RDG}:   }
  Define $\ol{\U}_t  \df   \ol{V}_t \+ \int_0^t g_r dr $, $ t \in [0,T]$ as in Lemma \ref{lem_submg0}.

  Given $(t,\o) \ins [0,T] \ti \O$ and $n \ins \hN$, since Remark \ref{rem_V_adapted},
 Proposition \ref{prop_shift0} (4) and Proposition
  \ref{prop_conti_V} show that   $ (V^n)^{t,\o}     \n - \n  L^{t,\o}  $ is an $\bF^t-$adapted process with  left-upper-semicontinuous paths
  and that    $\ol{V}^{t,\o}   \n - \n  L^{t,\o}  $ is
  an $\bF^t-$adapted process with all continuous paths,
   we can deduce from \eqref{eq:ek011}   that
   \beas
      \tau^{n,\d}_{(t,\o)} \n : = \n
   \inf\big\{ s  \n \in \n  [t,T] \n :  (V^n)^{t,\o}_s   \n < \n  L^{t,\o}_s \n + \n \d  \big\} , \q  \fa   \d  \n > \n  0
  \eeas
  are all $\bF^t-$optional times and that
     \beas 
     \tau^*_{(t,\o)} \n  : = \n     \inf\big\{ s  \n \in \n  [t,T] \n :
      \ol{V}^{t,\o}_s   \n = \n  L^{t,\o}_s      \big\}
  \n = \n    \inf \big\{ s  \n \in \n  [t,T] \n :  \ol{V}^{t,\o}_s   \n \le  \n  L^{t,\o}_s      \big\}
  \eeas
  is an $\bF^t-$stopping time.

 \no {\bf 1)}
Let  $(t,\o)  \n \in \n  [0,T]  \n \times \n  \O $ and  $  \ga  \n \in \n  \cT^t $.
Since    $\ga (\Pi^0_t) \ins \cT_t $  by \eqref{eq:xax001},
 Taking $ t' \= t $ and $ \z   \=   \ga (\Pi^0_t)   $
  in \eqref{eq:dc251} of Lemma \ref{lem_submg0} shows that
 \bea  \label{eq:dc253}
  \ol{V}_t (\o) + \int_0^t g_r (\o) dr =
    \ol{\U}_t (\o) \le  \underset{\hP \in \cP(t , \o)}{\inf} \,    \hE_\hP \bigg[ \Big( \,
 \ol{\U}_{   \big( \tau^*_{(t,\o)} (\Pi^0_t)     \land    \ga (\Pi^0_t)  \big) \vee t  }  \Big)^{t, \o} \bigg] .
 \eea
 For any $\wt{\o} \in \O^t$,   \eqref{eq:uxu170} and the first equality in \eqref{eq:xax141} imply  that
 \beas
 && \hspace{-1cm}  \Big( \,  \ol{\U}_{   \big( \tau^*_{(t,\o)} (\Pi^0_t)
  \land  \ga (\Pi^0_t)  \big) \vee t  }  \Big)^{t, \o} (\wt{\o})
 \= \ol{\U} \Big( \Big(\tau^*_{(t,\o)} \big(\Pi^0_t (\o \oti_t \wt{\o})\big)
   \n  \land  \n    \ga \big(\Pi^0_t (\o \oti_t \wt{\o})\big)  \Big) \vee t,
  \o \oti_t \wt{\o} \Big)
  \= \ol{\U} \big(  \tau^*_{(t,\o)}  (  \wt{\o} )
    \n \land \n     \ga (  \wt{\o} )     ,   \o \oti_t \wt{\o} \big) \\
  &&   \= \ol{V}^{t,\o} \big(  \tau^*_{(t,\o)}  (  \wt{\o} )
    \n \land \n     \ga (  \wt{\o} )     ,     \wt{\o} \big) \+ \int_0^{\tau^*_{(t,\o)}  (  \wt{\o} )
   \land    \ga (  \wt{\o} ) } g_r (\o \oti_t \wt{\o}) dr \\
  &&   \le \b1_{\{\tau^*_{(t,\o)}  (  \wt{\o} )   \le     \ga (  \wt{\o} )\}}
   L^{t,\o} \big(  \tau^*_{(t,\o)}  (  \wt{\o} )     ,     \wt{\o} \big)
   \+ \b1_{\{ \ga (  \wt{\o} ) < \tau^*_{(t,\o)}  (  \wt{\o} )      \}}
   U^{t,\o} \big( \ga (  \wt{\o} )     ,     \wt{\o} \big)
   \+ \int_0^t g_r (\o) dr \+ \int_t^{\tau^*_{(t,\o)}  (  \wt{\o} )
   \land    \ga (  \wt{\o} ) } g^{t,\o}_r (  \wt{\o}) dr \\
  &&   \= \big( R^{t,\o}  ( \tau^*_{(t,\o)}  ,   \ga    ) \big) (  \wt{\o} ) \+ \int_0^t g_r (\o) dr  .
 \eeas
 Plugging this into \eqref{eq:dc253}  yields that
 $ \ol{V}_t (\o)  \le  \underset{\hP \in \cP(t , \o)}{\inf} \,
    \hE_\hP \Big[ R^{t,\o} \big( \tau^*_{(t,\o)}  ,   \ga   \big) \Big] $.
    Taking infimum over $\ga \in \cT^t$ leads to that
 \beas
 \q \ol{V}_t (\o)  \ls \underset{\ga \in \cT^t}{\inf} \, \underset{\hP \in \cP(t , \o)}{\inf} \,
    \hE_\hP \big[ R^{t,\o} \big( \tau^*_{(t,\o)}  ,   \ga   \big) \big]
    \ls  \underset{\tau \in \cT^t}{\sup} \,  \underset{\ga \in \cT^t}{\inf} \, \underset{\hP \in \cP(t , \o)}{\inf} \,
    \hE_\hP \big[ R^{t,\o}  ( \tau  ,   \ga    ) \big]  \= \ul{V}_t (\o) \ls  \ol{V}_t (\o) , \q
    \hb{proving \eqref{eq:xax015}} .
 \eeas

      \no {\bf 2)} Let $  \zeta  \n \in \n  \cT $ and $(t,\o)  \n \in \n  [0,T]  \n \times \n  \O $.
       If $\wh{t}  \n := \n  \tau_*(\o   )   \n \land \n  \z (\o   )  \n \le \n  t$,
       similar to \eqref{eq:dc311}, we can deduce from
   Proposition \ref{prop_shift0} (3), the $\bF-$adaptedness of $\U $   by Remark \ref{rem_V_adapted}
    as well as  \eqref{eq:bb421}  that
   $  \big( \U_{\tau_* \land \z  } \big)^{t,\o} ( \wt{\o} )
    \n = \n  \U \big(  \wh{t},      \o \big) $,
   $\fa \wt{\o}  \n \in \n  \O^t$.
   Then
   \bea  \label{eq:dc249}
   \ul{\sE}_t \big[ \,  \U_{\tau_* \land \z  }  \big] (\o)
   = \underset{\hP \in \cP(t,\o)}{\inf} \hE_\hP \Big[ \big( \U_{\tau_* \land \z  } \big)^{t,\o}  \Big]
   = \underset{\hP \in \cP(t,\o)}{\inf} \hE_\hP \big[ \, \U (  \wh{t},      \o  ) \big]
   = \U \big(  \wh{t},      \o \big)= \U  \big(   \tau_*(\o) \land \z (\o) \land t  ,   \o \big) .
   \eea
     On the other hand, if $  \tau_* (\o   )    \land    \z (\o   )  \n > \n  t$,
     applying  Proposition \ref{prop_shift0} (3) once again shows that
      $  \o   \oti_t   \O^t   \sb   \{\tau_*    \land   \z     \>   t   \}    $.
  So it holds for any $ \wt{\o}  \n \in \n  \O^t$  that
     $  \big( \U_{\tau_*    \land    \z   } \big)^{t,\o} (\wt{\o})
   \n = \n  \U_{\tau_* \land \z   } \big( \o    \otimes_t    \wt{\o}  \big)
   \n = \n  \U_{  ( \tau_* \land \z )  \vee t   } \big( \o    \otimes_t    \wt{\o}  \big)
   \n = \n   \big( \U_{  ( \tau_* \land \z )  \vee t   } \big)^{t,\o} (\wt{\o})   $.
  As $\tau_* \= \tau^*_{(0,\bz)} \= \tau^*_{(0, \o)} $,
  taking $ t'  \=  0 $ in \eqref{eq:dc251}   yields that
  \beas
    \U_{\tau_*   \land \zeta   \land  t}       ( \o) = \U_t (\o) \le
    \underset{\hP \in \cP(t , \o)}{\inf} \,
    \hE_\hP \Big[ \big( \U_{  ( \tau_*     \land    \z  )  \vee t  } \big)^{t, \o} \Big]
    = \underset{\hP \in \cP(t , \o)}{\inf} \,
    \hE_\hP \Big[ \big( \U_{   \tau_*     \land    \z     } \big)^{t, \o} \Big]
    =  \ul{\sE}_t  \big[ \, \U_{  \tau_*  \land \zeta  }   \big] ( \o)  ,
  \eeas
   which together with \eqref{eq:dc249} proves \eqref{eq:cc761}. \qed

\subsection{Proof of Proposition \ref{prop_P1P2P3_Ass}}

For any $\a,\d \n \in \n (0,\infty)$,
we define $ \Phi (\a,\d) : = \vr_0   (  \d \n +  \n   \d^{1/4}    )  \n + \n
   \k    (1 \n + \n 2^{\varpi-1} \d^\varpi) \vf_1 (\a) \d^{1/4}
  \n + \n  \k 2^{\varpi-1}  \vf_{\varpi+1} (\a) \d^{\varpi/2+1/4} $.

 \no {\bf 1)} {\it   we first  show that the  probability class
$\{\cP(t,\o)\}_{(t,\o)    \in    [0,T]    \times    \O }$ satisfies \(P1\) and \(P2\).}

  Let $(t,\o) \n \in \n  [0,T]  \n \times \n  \O$ and   $\mu  \n \in \n  \sU_t$.
We  set $(\hP,\fp, \cX) \n : = \n  \big(\hP^{t,\o,\mu},\fp^{t,\o,\mu}, X^{t,\o,\mu}  \big)$.
  Given $\wt{\o}  \n \in \n  \O^t$, \eqref{eq:ab025} shows that
   \bea \label{eq:xax111}
 ~ \; \big|\Psi^{t,\bz}_r (\cX (\wt{\o})) \n - \n  \Psi_r (\bz) \big|
   \n = \n  \big|\Psi_r ( \bz  \n \otimes_t \n  \cX (\wt{\o}))  \n - \n  \Psi_r (\bz) \big|
    \n \le \n  \vr_0 \big( \|\bz  \n \otimes_t \n  \cX (\wt{\o})\|_{0,r} \big)
   \n \le \n  \k \big( 1 \n + \n  \|\cX (\wt{\o})\|^\varpi_{t,r} \big),
   ~ \fa r  \n \in \n  [t,T] .
   \eea
   It follows that
    $
  \Psi^{t,\bz}_* (\cX (\wt{\o}))
  = \underset{r \in [t,T]}{\sup} \big| \Psi^{t,\bz}_r (\cX (\wt{\o})) \big|
   \le \k \big( 1+ \|\cX (\wt{\o})\|^\varpi_{t,T} \big) + M^\Psi_0  $,
    where    $ M^\Psi_0 := \underset{r \in [t,T]}{\sup}
   \big|  \Psi_r (\bz) \big| < \infty$ by the continuity of path $\Psi_\cd (\bz)$.
   Since $\Psi^{t,\bz}$ is an $\bF^t-$adapted process by Proposition \ref{prop_shift0} (4),
   applying \eqref{eq:xxx153} 
   yields  that
    \beas
      \hE_{\hP} \big[  \Psi^{t,\bz}_* \big]   \n =  \n    \hE_{\fp} \big[  \Psi^{t,\bz}_* \big]
    \n =  \n  \hE_t \big[  \Psi^{t,\bz}_* (\cX) \big]
   \n  \le  \n    \k  \big( 1 \n + \n  \hE_t \big[\|\cX  \|^{  \varpi}_{t,T} \big] \big)  \n    + \n
     M^\Psi_0
     \n  \le     \n     \k  \big( 1 \n + \n   \vf_{    \varpi } \big( \|\o\|_{0,t} \big) \, T^{\varpi/2} \big)
      \n   + \n     M^\Psi_0   \n < \n  \infty  .
      \eeas
   Namely, $ \Psi^{t,\bz} \in  \hS  (\bF^t, \hP) $. Similar to \eqref{eq:xax111},
   one can deduce from \eqref{eq:aa211} that
   $ \big|g^{t,\bz}_r (\cX (\wt{\o})) \n - \n  g_r (\bz) \big|
   \n \le \n  \k \big( 1 \n + \n  \|\cX (\wt{\o})\|^\varpi_{t,r} \big) $ for any
   $ r  \n \in \n  [t,T] $.
 Then    Fubini's Theorem and \eqref{eq:xxx153} imply that
   \bea
   \hE_{\hP} \int_t^T |g^{t,\bz}_r| dr
   & \tn \dn =  & \tn  \dn
    \hE_{\fp} \int_t^T |g^{t,\bz}_r| dr \n = \n \hE_t  \int_t^T \n |g^{t,\bz}_r (\cX)| dr
    \le    \k  \int_t^T  \n  \big( 1 \n + \n \hE_t [ \|\cX  \|^\varpi_{t,T} ] \big) dr
     \n + \n    \int_t^T \n  | g_r (\bz) |  dr \nonumber \\
     & \tn  \dn  \le  & \tn  \dn  \k T \big( 1 \n + \n   \vf_{    \varpi } \big( \|\o\|_{0,t} \big) \, T^{\varpi/2} \big)
    \n + \n   \int_t^T \n  | g_r (\bz) |  dr  \n < \n  \infty . \q \hb{Hence $ \hP   \in \wh{\fP}_t $.}
    \label{eq:uxu217}
   \eea

 For any $t \ins [0,T]$ and $\o_1,\o_2  \n \in \n  \O$ with $\o_1|_{[0,t]}  \n = \n  \o_2|_{[0,t]} $,
 since the SDE \eqref{FSDE1}   depends only on $\o|_{[0,t]}$ for a given path $\o \ins \O$, we see that
 $ X^{t,\o_1,\mu} = X^{t,\o_2,\mu}$  and thus
 $\hP^{t,\o_1,\mu}= \hP^{t,\o_2,\mu}$ for any $\mu \in \sU_t$.
 It follows  that $ \cP(t,\o_1) = \cP(t,\o_2) $. So Assumption  \(P1\) is satisfied.
 Also, Proposition 6.3 of \cite{ROSVU} has already shown that
 the probability class $\{\cP(t,\o)\}_{(t,\o) \in [0,T] \times \O}$
 satisfies  (P2).

 \no {\bf 2)}   {\it The verification that   the  probability class
$\{\cP(t,\o)\}_{(t,\o)    \in    [0,T]    \times    \O }$ satisfies \(P3\) is relatively lengthy.
We split it  into several steps. }

 \no {\bf 2a)} {\it Let us first quote some knowledge on the inverse   mapping of $X^{t,\o,\mu}$ from \cite{ROSVU},
 which has already verified \(P3\) \(i\), \(ii\) for $\{\cP(t,\o)\}_{(t,\o)    \in    [0,T]    \times    \O }$. }

   Given $(t,\o) \n \in \n  [0,T]  \n \times \n  \O$ and $\mu  \n \in \n  \sU_t$,
 according to   \cite{ROSVU}  (see the context around (7.62) and (7.63)  therein),
 there exists an   $\bF^t-$progressively measurable  process $W^{t,\o,\mu}  $ such that
 for all  $ \wt{\o} \in \O^t $ except on a   $\hP^t_0-$null set $\cN_{ t,\o,\mu }$
  \beas 
   B^t_s ( \wt{\o} ) =  W^{t,\o,\mu}_s \big(  X^{t,\o,\mu} ( \wt{\o} ) \big) , \q  \fa  s \in [t,T] ,
  \eeas
 and that the  $\fp^{t,\o,\mu}$ probability of set
     $A_{t,\o,\mu}  \n  := \n  
   \{\wt{\o}'  \n \in \n  \O^t \n :  \cN^c_{t,\o,\mu}  \n \cap \n  ( X^{t,\o,\mu})^{-1} (\wt{\o}')  \n \ne \n  \es  \}
  $ is $1$, i.e.,     $ A^{\,\raisebox{0.2ex}{\scriptsize $c$}}_{t,\o,\mu}  \n \in \n
    \sN^{\fp^{t,\o,\mu}} \n := \n  \big\{   A
      \n \in \n  \cG^{X^{t,\o,\mu} }_T  \n :
        \fp^{t,\o,\mu}  ( A )  \n = \n  0  \big\} $.  
  For any $r  \n \in \n  [t,T]$,    \eqref{eq:xxx439} and Lemma A.3   (2) of \cite{ROSVU} show that
  $\fF^{t,\o,\mu}_r \n := \n  \si \big(\cF^t_r \cup \sN^{\fp^{t,\o,\mu}}   \big)  \n \subset \n  \cG^{X^{t,\o,\mu}}_r$.

  We   see from  the context around (7.67)$-$(7.69) of \cite{ROSVU}   that  $\wt{W}^{t,\o,\mu}_r (\wt{\o})
    \n := \n  \b1_{\{\wt{\o} \in A_{t,\o,\mu} \}} W^{t,\o,\mu}_r (\wt{\o}) $,
   $ (r,\wt{\o})  \n \in \n  [t,T]  \n \times \n  \O^t $ is an  $ \{\fF^{t,\o,\mu}_r\}_{r   \in   [t,T]}-$adapted
   process such that all its paths belong to $ \O^t $,   that
  \bea   \label{eq:xax069}
   \wt{\o} = B^t (\wt{\o})   = W^{t,\o,\mu}  \big( X^{t,\o,\mu}  (\wt{\o}) \big)
 = \wt{W}^{t,\o,\mu}  \big( X^{t,\o,\mu}  (\wt{\o}) \big), \q \fa \wt{\o} \in \cN^c_{t,\o,\mu} ,
  \eea
 and that
    \bea \label{eq:xax047}
   \big(\wt{W}^{t,\o,\mu}\big)^{-1}(A') \in      \fF^{t,\o,\mu}_r , \q    \fa    A' \in \ol{\cF}^t_r, ~ \fa r \in [t,T] .
   \eea

  Fix $0 \n \le \n  t  \n < \n  s  \n \le \n  T$,  $\o  \n \in \n  \O$ and $\mu  \n \in \n  \sU_t$,
  $ \d  \n \in \n  \hQ_+   $ and $\l  \n \in \n  \hN$.
   We consider    a $\cF^t_s-$partition $\{\cA_j\}^\l_{j=0}$ of $\O^t$ such that for $j \n = \n 1,\cds \n , \l$,
  $\cA_j  \n \subset \n  O^s_{\d_j} (\wt{\o}_j)$ for some
  $\d_j  \n \in \n  \big( (0,\d]  \n \cap \n  \hQ \big)  \n \cup \n  \{\d\}$ and  $\wt{\o}_j  \n \in \n  \O^t $,
  and let   $ \{\mu^j\}^\l_{j=1}  \n \subset  \n    \sU_s    $.
  We will simply set
  \bea \label{eq:xax153}
   (   \hP,\fp,  \cX,    \cW,   \fF_\cd   ) :=
      \big( \hP^{t,\o,\mu},\fp^{t,\o,\mu},   X^{t,\o,\mu} ,   \wt{W}^{t,\o,\mu},   \fF^{t,\o,\mu}_\cd \big)   .
    \eea

    Given $j=1,\cds \n , \l $, \eqref{eq:xxx439} shows that $\cA^\cX_j \n := \n  \cX^{-1}(\cA_j)  \n \in \n  \ol{\cF}^t_s$.
   So there exists an $A_j   \n \in \n  \cF^t_s$   such that $ \cA^\cX_j  \, \D \,  A_j   \n  \in \n  \ol{\sN}^t    $
       (see e.g. Problem 2.7.3 of \cite{Kara_Shr_BMSC}).
       Following similar arguments to those used in the proof of Proposition 6.3 of \cite{ROSVU},
       one can  show  that

        \no (u1) The set
        $ \wt{A}_j  \n := \n  A_j \big\backslash   \underset{j'<j}{\cup}  A_{j'}  \n \in \n  \cF^t_s   $ satisfies
     $  \cA^\cX_j \D \wt{A}_j \in \ol{\sN}^t $ \big(see (7.70) of \cite{ROSVU}\big).

     \no (u2) The pasted control
 $  \wh{\mu}_r (\wt{\o}) \n := \n
 \b1_{\{r \in [t,s)\}} \mu_r (\wt{\o})  \n + \n  \b1_{\{r \in [ s,T]\}}
 \Big( \b1_{\{\wt{\o} \in \wt{A}_0\}} \mu_r (\wt{\o})
  \n + \n  \sum^\l_{j=1}  \b1_{\{\wt{\o} \in \wt{A}_j\}} \mu^{j}_r (\Pi^t_s \big(\wt{\o})\big) \Big) $,
 $ \fa (r, \wt{\o})  \n \in \n  [t,T]  \n \times \n  \O^t $
 belongs to $ \sU_t $, where $\wt{A}_0 := \Big( \underset{j=1}{\overset{\l}{\cup}} \wt{A}_j \Big)^c \in \cF^t_s $
 \big(see (7.71) of \cite{ROSVU}\big).
 Set
 \beas
  \big(   \wh{\hP}, \wh{\fp},  \wh{\cX}, \wh{\cW},   \wh{\fF}_\cd , \wh{\cN}  \big) \n := \n
   \big( \hP^{t,\o,\wh{\mu}},\fp^{t,\o,\wh{\mu}},   X^{t,\o,\wh{\mu}} ,   \wt{W}^{t,\o,\wh{\mu}},
  \fF^{t,\o,\wh{\mu}}_\cd , \cN_{t,\o,\wh{\mu}} \big)   .
  \eeas

     \no (u3) There exists a $\hP^t_0-$null set
     $\wt{\cN}_j  $  such that for any $\wt{\o} \in \wt{A}_j \cap \wt{\cN}^c_j $,
     \beas
  \q   \cN_{\wt{\o}} \df  \big\{  \wh{\o} \ins \O^s \n :  \wh{\cX}_r    (\wt{\o} \oti_s \wh{\o})
       \n \ne \n  \big( \cX  (\wt{\o}  )     \oti_s     X^{s, \o \otimes_t  \cX  ( \wt{\o} ), \mu^j } (\wh{\o})\big)(r)          \hb{ for some } r  \n \in \n  [t,T] \big\} \hb{ belongs to } \ol{\sN}^s  ~
        \hb{\big(see (7.78) of \cite{ROSVU}\big). }
     \eeas

  \no    (u4) For any $A \in \cF^t_s$,  $  \cX^{-1} (A) \D \wh{\cX}^{-1} (A) \in \ol{\sN}^t  $
  \big(see (7.74) of \cite{ROSVU}\big).

   Also,  analogous   to    part (2b) of \cite[Proposition 6.3]{ROSVU},
    we can use the uniqueness of controlled SDE \eqref{FSDE1}
 to  show that the equality  $ \wh{\mu} = \mu$
 over $\big([t,s] \times \O^t \big) \cup \big( [s,T] \times  \wt{A}_0 \big) $
 implies the equality  $  \wh{\cX}    = \cX    $
 over $\big([t,s] \times \O^t \big) \cup \big( [s,T] \times  \wt{A}_0 \big) $, and thus that
 $ \wh{\hP}$  satisfies \(P3\) \(i\), \(ii\).

   \no {\bf 2b)}  {\it    To show that $\wh{\hP}$ satisfies   \eqref{eq:xxx617},
   we make some technical setting and preparation first.  }

     Proposition \ref{prop_shift0} (4) shows that  $ \cY^1_r := g^{t,\o}_r $,
 $ \cY^2_r := L^{t,\o}_r    $ and $ \cY^3_r := U^{t,\o}_r    $, $   r \in [t,T]$
 are three $ \bF^t-$adapted processes with all continuous paths.  For $\ell=1,2,3$,    \eqref{eq:xxx439} implies that
  $\cY^\ell\big(\wh{\cX}\big)$   is an $\ol{\bF}^t-$adapted  process with all continuous paths.
 Applying Lemma A.2 (3) of \cite{ROSVU} with
 $(\hP,X) = (\hP^t_0,B^t)$ shows that $\cY^\ell \big(\wh{\cX}\big)$ has an $(\bF^t,\hP^t_0)-$version $\sY^\ell $.
 More precisely, $\sY^\ell$'s   are   $\bF^t-$progressively   measurable  processes such that
 \bea \label{eq:xax123}
  \cN_R \n := \n \underset{\ell=1}{\overset{3}{\cup}}
  \big\{\wt{\o}  \n \in \n  \O^t \n : \sY^\ell_r (\wt{\o})  \n \ne \n  \cY^\ell_r \big(\wh{\cX} (\wt{\o})\big)
  \hb{ for some } r \n \in \n  [t,T]  \big\}  \n \in \n  \ol{\sN}^t  .
  \eea
  By Lemma \ref{lem_null_sets}, it holds for all  $\wt{\o}  \n \in \n  \O^t$ except on  an
 $\wt{\cN}_R  \n \in \n  \ol{\sN}^t $ that $ \big(\cN_R \n \cup \n \wh{\cN}\big)^{s,\wt{\o}}
  \n \in \n  \ol{\sN}^s $.

    We see from Proposition \ref{prop_shift0} (4)   that the random variables
  \bea \label{eq:xax155}
   \xi_m  \n := \n  \underset{t' \in [t,T]}{\sup}  \int_{t'}^{(t'+2^{-m}) \land T} \n  \big| g^{t,\o}_r \big|  dr ,
   \q \fa m \n \in \n  \hN
   \eea
 are $\cF^t_T-$measurable.
 Since $\lmtd{m \to \infty} \xi_m  \n = \n  0$,   \eqref{eq:xxx111} and the dominated convergence theorem
 show that $ \lmtd{m \to \infty} \hE_{\wh{\hP}} [ \xi_m ]  \n = \n  0 $. So there exists
 $\fm  \n \in \n  \hN$ such that $\hE_{\wh{\hP}} [ \xi_\fm ]  \n \le \n  \d / 2 $
 and  $ \Phi \big(\|\o\|_{0,t}, 2^{-\fm} \big)  \n \le \n  \d / 2 $. Set $ \fra   \n := \n 2^{-\fm}  $.

  Now, fix $n \n \in \n  \hN  \n \cup \n  \{\infty\} $, $  \wp  \n \in \n  \cT^s$
 and let     $j \n = \n 1,\cds  \n ,\l $.  We   set
  \beas
  (\hP^j, \fp^j, \cX^j,  \cW^j, \fF^j_\cd, \cN_{\cX^j} )
  \n := \n  \big( \hP^{  s, \o \otimes_t \wt{\o}_j , \mu^j },
  \fp^{  s, \o \otimes_t \wt{\o}_j , \mu^j },    X^{  s, \o \otimes_t \wt{\o}_j , \mu^j } ,
  \wt{W}^{  s, \o \otimes_t \wt{\o}_j , \mu^j } ,   \fF^{  s, \o \otimes_t \wt{\o}_j , \mu^j }_\cd ,
  \cN_{  s, \o \otimes_t \wt{\o}_j , \mu^j }   \big)
  \eeas
  and define
 \bhe \bea \label{eq:uxu214}
   \wp_j  \n : = \n  \wp  ( \cX^j ) \ins \ol{\cT}^s  , \q
   \nu_j  \n : = \n  \wp_j ( \Pi^t_s ) \ins  \ol{\cT}^t_s  , \q
   \wh{\ga}_j  \n : = \n  \nu_j \big( \wh{\cW} \big)   ,
 \eea \ehe
 where $\wh{\ga}_j$ is a $\wh{\fF}-$stopping time  that takes values in $[s,T]$.

    Given $i  \n   =    \n    0 , \cds   \n   , 2^\fm  $,
    we set    $s_i    \n  :=  \n   s  \n \vee \n  (i 2^{-\fm} T)     $ and
 $D_i  \n := \n  \{ s_{i-1}    \n < \n  \wh{\ga}_j
  \n \le \n  s_i  \}  \n \in \n \wh{\fF}_{ s_i }     $ with $s_{-1}  \n := \n  -1 $.
  By e.g.   Problem 2.7.3 of \cite{Kara_Shr_BMSC},
  there exists an $\wt{D}_i   \n \in \n  \cF^t_{ s_i }$   such that
  $ D_i \, \D \,  \wt{D}_i   \n \in \n  \sN^{\wh{\fp}}  $.
     Define    $\cD_i  \n := \n  \wt{D}_i \backslash \underset{i' < i}{\cup}
     \wt{D}_{i'}  \n \in \n  \cF^t_{ s_i }  $
        and $ \ol{\cD}  \n := \n  \underset{i=0}{\overset{2^\fm}{\cup}}  \cD_i
      \n = \n  \underset{i=0}{\overset{2^\fm}{\cup}}  \wt{D}_i  \n \in \n  \cF^t_T $.
     Then $  \ga'_j   \n := \n  \sum^{2^\fm}_{i=0} \b1_{D_i}   s_i$
     is a $\wh{\fF}-$stopping time while
     $  \ga_j   \n := \n  \sum^{2^\fm}_{i=0} \b1_{\cD_i}   s_i
        +    \b1_{\ol{\cD}^c} T$ defines an $\cT^t_s-$stopping time.
     Clearly, $\ga'_j   $ coincides with $\ga_j$ over
     $\underset{i=1}{\overset{2^\fm}{\cup}} \big(  D_i \cap \cD_i \big)$, whose complement
     $ \underset{i=1}{\overset{2^\fm}{\cup}} \big(  D_i \backslash \cD_i \big) $
     belongs to $\sN^{\wh{\fp}}$      because
     $$
        D_i \backslash \cD_i \=  D_i \cap
      \Big[ \big( \wt{D}_i \big)^c \cup \Big( \underset{i' < i}{\cup} \wt{D}_{i'} \Big)  \Big]
       \= \big( D_i \backslash \wt{D}_i \big) \cup \Big( \underset{i' < i}{\cup} \big(  \wt{D}_{i'} \cap D_i \big)  \Big)
       \sb \big( D_i \D \wt{D}_i \big) \cup \Big( \underset{i' < i}{\cup} \big( \wt{D}_{i'} \cap D^c_{i'}  \big)  \Big)
      \sb   \underset{i' \le i}{\cup} \big( D_{i'} \D \wt{D}_{i'} \big) \ins  \sN^{\wh{\fp}} .
     $$
 for $i \=  1, \cds  \n , 2^\fm$.  To wit, we have
  \bea \label{eq:xax143}
   \ga'_j = \ga_j    ,  \q \wh{\fp}-a.s.
  \eea

   \no {\bf 2c)}  {\it Now, fix $A \in \cF^t_s$,  $\t \in \cT^t_s (n)$ and
 set   $\wh{\tau} \n := \n  \tau \big( \wh{\cX} \big) $. We   show an auxiliary inequality: }
 \bea \label{eq:xxx849}
  \sum^\l_{j=1}  \hE_{\wh{\hP}} \big[ \b1_{A \cap \cA_j} R^{t,\o}  \big(\tau, \ga_j   \big) \big]
 \n  \le \n  \sum^\l_{j=1}  \hE_t \big[ \b1_{\wh{\cX}^{-1} (  A \cap \cA_j )} \Xi_j \big]  \n + \n  \d  ,
 \eea
 where $ \Xi_j \n := \n \int_t^{\wh{\tau}  \land \nu_j  } \sY^1_r    dr
 + \b1_{\{ \wh{\tau}  \le \nu_j  \}} \sY^2_{   \wh{\tau}  }
 + \b1_{\{ \nu_j  < \wh{\tau}    \}} \sY^3_{ \nu_j   } $.

 For any $r \in [s,T]$, an analogy to \eqref{eq:xax119}
    shows that $ \{\wh{\tau}  \n \le \n  r \} \n = \n \wh{\cX}^{-1} \big(\{\tau \le r \}\big)
    \in \ol{\cF}^t_r $,  So $\wh{\tau} \n \in \n \ol{\cT}^t_s   $.
 By Lemma 2.5 (3) in the ArXiv version of \cite{ROSVU},
 it holds    for all $\wt{\o} \n \in \n \O^t$ except on a
 $ \cN_\tau \ins  \ol{\sN}^t$ that   $  \wh{\t}^{s, \wt{\o}} \n \in \n \ol{\cT}^s  $.
 For $j \n = \n 1,\cds  \n ,\l $,
 since $\sY^\ell$'s   are   $\bF^t-$progressively   measurable  processes
 and since $\nu_j  $ is a  $ \ol{\cT}^t_s -$stopping time,
 we see that    $\Xi_j$ is an   $\ol{\cF}^t_T-$measurable random variable.

 Let     $j \n = \n 1,\cds  \n ,\l $. By  \eqref{eq:xax143},
 \bea \label{eq:xax125}
 \hE_{\wh{\hP}} \big[ \b1_{A \cap \cA_j} R^{t,\o}  \big(\tau, \ga_j   \big) \big]
 \n = \n  \hE_{\wh{\fp}} \big[ \b1_{A \cap \cA_j} R^{t,\o}  \big(\tau, \ga_j   \big) \big]
  \n = \n  \hE_{\wh{\fp}} \big[ \b1_{A \cap \cA_j} R^{t,\o}  \big(\tau, \ga'_j   \big) \big]
  \n = \n  \hE_t \Big[ \b1_{\wh{\cX}^{-1} (  A \cap \cA_j )}  R^{t,\o}  \big(\tau, \ga'_j   \big) \big(\wh{\cX}\big) \Big] .
 \eea
   Given $\wt{\o}  \n \in \n  \O^t$,
 since $ 0  \n \le \n  \ga'_j (\wt{\o})  \n - \n  \wh{\ga}_j (\wt{\o})   \n  < \n  \fra  $,
 \eqref{eq:aa211} implies that
\beas
&& \hspace{-1cm} R^{t,\o}  \big(\tau, \ga'_j   \big)   (\wt{\o})
\n - \n  R^{t,\o}  \big(\tau, \wh{\ga}_j   \big)   (\wt{\o})
  \n = \n   \int_{\tau (\wt{\o}) \land  \wh{\ga}_j (\wt{\o})}^{\tau (\wt{\o}) \land  \ga'_j (\wt{\o})}
 g^{t,\o}_r (\wt{\o}) dr
  \n + \n  \b1_{\{\wh{\ga}_j (\wt{\o}) < \tau (\wt{\o}) \le  \ga'_j (\wt{\o})\}}
 \big( L^{t,\o} (\tau (\wt{\o}), \wt{\o})  \n - \n  U^{t,\o} (\wh{\ga}_j (\wt{\o}), \wt{\o}) \big) \nonumber \\
&&       +    \b1_{\{\ga'_j (\wt{\o}) < \tau (\wt{\o})    \}}
 \big( U^{t,\o} (\ga'_j (\wt{\o}), \wt{\o})  \n - \n  U^{t,\o} (\wh{\ga}_j (\wt{\o}), \wt{\o}) \big) \nonumber  \\
&& \hspace{-3mm} \le \n  \xi_\fm (\wt{\o})  \n + \n  \b1_{\{\wh{\ga}_j (\wt{\o}) < \tau (\wt{\o}) \le  \ga'_j (\wt{\o})\}}
 \vr_0 \Big( (\tau (\wt{\o})  \n - \n  \wh{\ga}_j (\wt{\o}) ) \n + \n
 \underset{r \in [0,T]}{\sup} \big|(\o  \n \otimes_t \n  \wt{\o}) \big(r  \n \land \n  \tau (\wt{\o}) \big)
  \n - \n  (\o  \n \otimes_t \n  \wt{\o}) \big(r  \n \land \n  \wh{\ga}_j (\wt{\o}) \big) \big| \Big)  \nonumber \\
&&     +   \b1_{\{\ga'_j (\wt{\o}) < \tau (\wt{\o})    \}}  \vr_0 \Big( (\ga'_j (\wt{\o})  \n - \n  \wh{\ga}_j (\wt{\o}) )
 \n + \n
 \underset{r \in [0,T]}{\sup} \big|(\o  \n \otimes_t \n  \wt{\o}) \big(r  \n \land \n  \ga'_j (\wt{\o}) \big)
  \n - \n  (\o  \n \otimes_t \n  \wt{\o}) \big(r  \n \land \n  \wh{\ga}_j (\wt{\o}) \big) \big| \Big)  \nonumber  \\
&&  \hspace{-3mm}  \le \n  \xi_\fm (\wt{\o})  \n + \n  \b1_{\{\wh{\ga}_j (\wt{\o}) < \tau (\wt{\o}) \le  \ga'_j (\wt{\o})\}}
 \vr_0 \Big( \fra  \n + \dn
 \underset{r \in [\wh{\ga}_j (\wt{\o}) , \tau (\wt{\o})]}{\sup}
 \big|  \wt{\o} (r)  \n - \n   \wt{\o} ( \wh{\ga}_j (\wt{\o}))  \big| \Big)
  \dn +  \n   \b1_{\{\ga'_j (\wt{\o}) < \tau (\wt{\o})    \}}  \vr_0 \Big( \fra  \n + \dn
 \underset{r \in [\wh{\ga}_j (\wt{\o}) , \ga'_j (\wt{\o})]}{\sup}
 \big|  \wt{\o} (r)  \n - \n   \wt{\o} ( \wh{\ga}_j (\wt{\o}))
  \big| \Big) \qq \qq  \nonumber   \\
&&  \hspace{-3mm}   \le \n  \xi_\fm (\wt{\o})  \n + \n  \vr_0 \Big( \fra +
 \underset{r \in \big[\nu_j  \big( \wh{\cW}(\wt{\o}) \big) ,
 \big(\nu_j  \big( \wh{\cW}(\wt{\o}) \big) +\fra\big) \land T \big]}{\sup}
 \Big|  \wt{\o} (r) -  \wt{\o} \big( \nu_j  \big( \wh{\cW}(\wt{\o}) \big) \big)
  \Big| \Big) . 
\eeas
Taking $\wt{\o} \n = \n  \wh{\cX} (\wt{\o}')$,
one can deduce  from \eqref{eq:xax069}   that for $\hP^t_0-$a.s.  $\wt{\o}' \n \in \n \O^t$,
\bea \label{eq:xax127}
R^{t,\o}  \big(\tau, \ga'_j   \big)   \big(\wh{\cX} (\wt{\o}') \big)
\n - \n  R^{t,\o}  \big(\tau, \wh{\ga}_j   \big)   \big(\wh{\cX} (\wt{\o}') \big)
\le \n  \xi_\fm \big(\wh{\cX} (\wt{\o}') \big)
  \n + \n  \vr_0 \Big( \fra + \underset{r \in [\nu_j (\wt{\o}')  , (\nu_j (\wt{\o}')  +\fra) \land T]}{\sup}
 \big|  \wh{\cX}_r (\wt{\o}')    \n  -  \n  \wh{\cX}_{   \nu_j } (\wt{\o}')  \big|  \Big) .
\eea
 Also, \eqref{eq:xax069} and  \eqref{eq:xax123}    show that for any
 $\wt{\o}' \n \in \n \big(\cN_R \n \cup \n \wh{\cN}\big)^c  $
\bea
&& \hspace{-1.5cm} R^{t,\o}  \big(\tau, \wh{\ga}_j   \big)   \big(\wh{\cX} (\wt{\o}') \big)
 \n = \n   \int_t^{\wh{\tau} (\wt{\o}') \land \nu_j (\wt{\o}') }
 \n \cY^1_r \big(\wh{\cX} (\wt{\o}') \big) dr
 \n + \n  \b1_{\{ \wh{\tau} (\wt{\o}') \le \nu_j (\wt{\o}') \}}
\cY^2 \big( \wh{\tau} (\wt{\o}') , \wh{\cX} (\wt{\o}') \big)
 \n + \n  \b1_{\{ \nu_j (\wt{\o}') < \wh{\tau} (\wt{\o}')   \}}
\cY^3 \big( \nu_j (\wt{\o}') , \wh{\cX} (\wt{\o}') \big) \nonumber  \\
& &  =   \int_t^{\wh{\tau} (\wt{\o}') \land \nu_j (\wt{\o}') } \n  \sY^1_r  ( \wt{\o}'   ) dr
 \n + \n  \b1_{\{ \wh{\tau} (\wt{\o}') \le \nu_j (\wt{\o}') \}}
\sY^2 \big( \wh{\tau} (\wt{\o}') ,  \wt{\o}'  \big)
 \n + \n  \b1_{\{ \nu_j (\wt{\o}') < \wh{\tau} (\wt{\o}')   \}}
\sY^3 \big( \nu_j (\wt{\o}') , \wt{\o}' \big)
 \n = \n  \Xi_j (\wt{\o}') .  \label{eq:xax129}
\eea

 Since $ \wh{\cX}^{-1} (  A \cap \cA_j ) \in \ol{\cF}^t_s $, $j \= 0, \cds, \l$ by \eqref{eq:xxx439}
 and since $\nu_j  $'s are    $ \ol{\cT}^t_s -$stopping times,
$\ol{\nu} \df \b1_{\wh{\cX}^{-1} (   \cA_0 )} T  \+ \sum^\l_{j=1} \b1_{\wh{\cX}^{-1} (   \cA_j )} \nu_j$
is also a $\ol{\cT}^t_s-$stopping time.
Set $\ol{\eta} \df \underset{r \in [\ol{\nu}   , (\ol{\nu}   +\fra) \land T]}{\sup}
 \big|  \wh{\cX}_r      \n  -  \n  \wh{\cX}_{ \ol{\nu} }    \big|$.
 Using the inequality $(a\+b)^\varpi \ls  2^{\varpi-1} (a^\varpi \+ b^\varpi)$, $\fa a,b \> 0$,
 one can deduce from \eqref{eq:xax129},    \eqref{eq:xax127}  and \eqref{eq:xxx153}   that
\bea
&& \hspace{-1 cm} \sum^\l_{j=1} \hE_t \Big[ \b1_{\wh{\cX}^{-1} (  A \cap \cA_j )}
\Big( R^{t,\o}  \big(\tau, \ga'_j   \big) \big(\wh{\cX}\big)  \n - \n \Xi_j \Big) \Big]
\n \le \n   \sum^\l_{j=1} \hE_t \bigg[ \b1_{\wh{\cX}^{-1} (  A \cap \cA_j )} \bigg( \xi_\fm  \big(\wh{\cX}  \big)
 \+ \vr_0 \Big( \fra + \underset{r \in [\nu_j    , (\nu_j    +\fra) \land T]}{\sup}
 \big|  \wh{\cX}_r      \n  -  \n  \wh{\cX}_{   \nu_j }    \big|  \Big) \bigg) \bigg] \nonumber \\
 && \hspace{-0.5cm} =  \sum^\l_{j=1} \hE_t \Big[ \b1_{\wh{\cX}^{-1}(  A \cap \cA_j )} \Big( \xi_\fm   \big(\wh{\cX}   \big)
 \+ \vr_0 \big(\fra \+ \ol{\eta}\big) \Big) \Big]
 \ls \hE_t \Big[ \xi_\fm   \big(\wh{\cX}   \big) \+ \vr_0 \big(\fra \+ \ol{\eta}\big) \Big]  \nonumber \\
 &&  \hspace{-0.5cm} \le  \hE_{\wh{\fp}}  [  \xi_\fm   ] \+  \hE_t \big[
    \b1_{ \{ \ol{\eta} \le \fra^{\frac14}   \}} \vr_0   (  \fra \n +  \n   \fra^{\frac14}    )  \n + \n
 \k \b1_{ \{ \ol{\eta} > \fra^{\frac14}   \}} \big( 1 \+ (\fra \n + \n \ol{\eta})^\varpi \big) \big]
      \ls  \hE_{\wh{\hP}}  [  \xi_\fm   ] \+ \vr_0   (  \fra \n +  \n   \fra^{\frac14}    )  \n + \n  \k \fra^{-1/4} \hE_t
 \Big[   (1\+2^{\varpi-1} \fra^\varpi) \ol{\eta}  \+ 2^{\varpi-1} \ol{\eta}^{\varpi+1}   \Big]   \nonumber  \\
  &&  \hspace{-0.5cm}  \le  \n \d/2 \+ \vr_0   (  \fra \n +  \n   \fra^{\frac14}    )  \n + \n
   \k    (1\+2^{\varpi-1} \fra^\varpi) \vf_1 (\|\o\|_{0,t}) \fra^{\frac14}
 \+ \k 2^{\varpi-1}  \vf_{\varpi+1} (\|\o\|_{0,t}) \fra^{\varpi/2+1/4}
  \=  \d/2 \+  \Phi \big(\|\o\|_{0,t}, 2^{-\fm} \big)  \ls  \d .      \label{eq:xax157}
\eea
 Then we see from   \eqref{eq:xax125}   that
 \beas
 \qq \sum^\l_{j=1}  \hE_{\wh{\hP}} \big[ \b1_{A \cap \cA_j} R^{t,\o}  \big(\tau, \ga_j   \big) \big]
  \= \sum^\l_{j=1} \hE_t \Big[ \b1_{\wh{\cX}^{-1} (  A \cap \cA_j )}
  R^{t,\o}  \big(\tau, \ga'_j   \big) \big(\wh{\cX}\big)   \Big]
 \n  \le \n  \sum^\l_{j=1}  \hE_t \big[ \b1_{\wh{\cX}^{-1} (  A \cap \cA_j )} \Xi_j \big]  \n + \n  \d ,
 \q \hb{proving \eqref{eq:xxx849}. }
 \eeas

 \no {\bf 2d)} {\it We are ready to use  \eqref{eq:da025}   and
  the estimate  \eqref{eq:xxx151} to verify  \eqref{eq:xxx617} for $\wh{\hP}$.   }

 Let $j \= 1, \cds \n , \l$ again.
 As $ \wh{\hP} \ins \wh{\fP}_t  $ by \eqref{eq:uxu217},
   \eqref{eq:xax129},  \eqref{eq:ab015} and \eqref{eq:xxx111}   imply  that
 \beas
 \hE_t \big[ |\Xi_j| \big] \le
 \hE_t \bigg[ \int_t^T \big|g^{t,\o}_r \big(\wh{\cX}\big)\big| dr  +   \Psi^{t,\o}_* \big(\wh{\cX}\big)  \bigg]
 = \hE_{\wh{\fp}} \bigg[ \int_t^T \big|g^{t,\o}_r  \big| dr  +   \Psi^{t,\o}_*    \bigg]
 = \hE_{\wh{\hP}} \bigg[ \int_t^T \big|g^{t,\o}_r  \big| dr  +   \Psi^{t,\o}_*    \bigg] < \infty .
 \eeas
  Since $ \wh{\cX}^{-1} (  A \cap \cA_j ) \in \ol{\cF}^t_s $,
 applying Lemma A.2 (1) of \cite{ROSVU}
 with $(\hP,X,\xi) \n = \n  \big(\hP^t_0,B^t, \Xi_j \big)$,
   using (u4) with $A  \n = \n  A \cap \cA_j$
   and applying Proposition 2.3 in the ArXiv version of \cite{ROSVU}
   with $(\hP, \xi) \n = \n  \big(\hP^t_0,  \Xi_j \big)$,
 we can  deduce from    Proposition \ref{prop_shift1} (1)  and  (u1) that
 \bea
  \hE_t \big[ \b1_{\wh{\cX}^{-1} (  A \cap \cA_j )} \Xi_j \big]
   & \tn \dn = & \tn \dn     \hE_t \Big[ \b1_{ \wh{\cX}^{-1} (  A \cap \cA_j ) }
 \hE_t \big[ \Xi_j \big| \ol{\cF}^t_s \big] \Big]
 \n = \n    \hE_t \Big[ \b1_{ \wh{\cX}^{-1} (  A \cap \cA_j ) }  \hE_t \big[ \Xi_j \big|  \cF^t_s \big] \Big]
 \=  \hE_t \Big[ \b1_{  \cX^{-1} (  A \cap \cA_j ) }   \hE_t \big[ \Xi_j \big|  \cF^t_s \big] \Big]  \nonumber  \\
  & \tn \dn   =   & \tn \dn
       \hE_t \Big[ \b1_{\{\wt{\o} \in \cX^{-1} (  A   ) \cap \cA^\cX_j  \}}
  \hE_s \big[   \Xi_j^{s,\wt{\o}}  \big] \Big]
  =    \hE_t \Big[ \b1_{\{\wt{\o} \in \cX^{-1} (  A   ) \cap \cA^\cX_j \cap \wt{A}_j  \}}
  \hE_s \big[  \Xi_j^{s,\wt{\o}}  \big] \Big]  .  \qq \q \label{eq:xax151}
 \eea

  Let $ \wt{\o} \n \in \n  \cA^\cX_j  \n \cap \n  \wt{A}_j  \n \cap    \wt{\cN}^c_j
    \cap    \wt{\cN}^c_R    \cap    \cN^c_\tau $.
  As  $  \wh{\t}^{s, \wt{\o}} \n \in \n \ol{\cT}^s  $,
  similar to $ \wh{\ga}_j  \n  = \n  \nu_j ( \wh{\cW})  $,
 $\z_{\wt{\o}} \n : = \n  \wh{\tau}^{s,\wt{\o}} (\cW^j)$ is   a $\fF^j-$stopping time.
  Let $\wh{\o} \n  \in  \n  \O^s$ such that $\wh{\o}$ is not in the $\hP^s_0-$null set
  $ \big(\cN_R   \cup   \wh{\cN}\big)^{s,\wt{\o}}     \cup    \cN_{\cX^j}   \cup   \cN_{\wt{\o}}   $,
  and define $ \D X^j_{\wt{\o}} (\wh{\o}) \df \big\|  X^{  s, \o \otimes_t \cX(\wt{\o}) , \mu^j }(\wh{\o})
   \- \cX^j (\wh{\o}) \big\|_{s,T} $\,.
 Taking $\wt{\o}' \= \wt{\o} \otimes_s \wh{\o} \n \in \n \big(\cN_R \n \cup \n \wh{\cN}\big)^c $ in \eqref{eq:xax129},
 we see from \eqref{eq:da021}, \eqref{eq:xax069}, (u3),  \eqref{eq:da025}
 as well as  an analogy to the second equality of \eqref{eq:xax141} that
 \beas
 \Xi^{s,\wt{\o}}_j (\wh{\o}) 
 & \tn \dn  =  & \tn \dn   R^{t,\o}  \big(\tau, \wh{\ga}_j   \big)   \big(\wh{\cX} (\wt{\o} \otimes_s \wh{\o}) \big)
 = R \Big( t, \tau \big(\wh{\cX} (\wt{\o} \otimes_s \wh{\o}) \big) ,
 \wh{\ga}_j \big(\wh{\cX} (\wt{\o} \otimes_s \wh{\o}) \big) ,
 \o \otimes_t  \big(\wh{\cX} (\wt{\o} \otimes_s \wh{\o}) \big) \Big) \\
  & \tn \dn   = & \tn \dn   R \Big( t, \wh{\tau}   (\wt{\o} \otimes_s \wh{\o})   ,
 \nu_j   (\wt{\o} \otimes_s \wh{\o})   ,
 \o \otimes_t  \big(\wh{\cX} (\wt{\o} \otimes_s \wh{\o}) \big)  \Big)
    = R \Big( t, \wh{\tau}^{s,\wt{\o}}   (  \wh{\o})   ,
 \wp_j   ( \wh{\o})   ,
 \o \otimes_t  \big(\wh{\cX} (\wt{\o} \otimes_s \wh{\o}) \big)  \Big) \\
    & \tn \dn   = & \tn \dn   R \Big( t, \z_{\wt{\o}}    \big( \cX^j ( \wh{\o})\big)   ,
 \wp   \big( \cX^j ( \wh{\o})\big)   ,
 \o \otimes_t    \big( \cX  (\wt{\o}  )  \n  \otimes_s  \n
  X^{s, \o \otimes_t  \cX  ( \wt{\o} ), \mu^j } (\wh{\o})\big)   \Big) \\
  & \tn \dn   \le & \tn \dn   R \Big( t, \z_{\wt{\o}}    \big( \cX^j ( \wh{\o})\big)   ,
 \wp   \big( \cX^j ( \wh{\o})\big)   ,
 \big(\o \otimes_t      \cX  (\wt{\o}  ) )  \n  \otimes_s  \n  ( \cX^j (\wh{\o})\big)  \Big)
 \n + \n  (1 \n + \n T ) \vr_0 \big( \D X^j_{\wt{\o}} (\wh{\o}) \big) \\
 & \tn \dn  =  & \tn \dn   R \Big( s, \z_{\wt{\o}}    \big( \cX^j ( \wh{\o})\big)   ,
 \wp   \big( \cX^j ( \wh{\o})\big)   ,
 \big(\o  \n \otimes_t  \n      \cX  (\wt{\o}  ) )  \n  \otimes_s  \n  ( \cX^j (\wh{\o})\big)  \Big)
  \n + \n  \int_t^s \n  g_r
  \big( (\o  \n \otimes_t  \n      \cX  (\wt{\o}  ) )  \n  \otimes_s  \n  ( \cX^j (\wh{\o})) \big) dr
  \n + \n  (1 \n + \n T ) \vr_0 \big( \D X^j_{\wt{\o}} (\wh{\o}) \big)  \\
   & \tn \dn  =  & \tn \dn   \big( R^{s,\o  \otimes_t  \cX  (\wt{\o}  )}
    (\z_{\wt{\o}},\wp )\big) \big( \cX^j ( \wh{\o})\big)
  \n + \n  \int_t^s \n  g_r \big( \o \otimes_t      \cX  (\wt{\o}  )  \big) dr
  \n + \n  (1 \n + \n T ) \vr_0 \big( \D X^j_{\wt{\o}} (\wh{\o}) \big) .
 \eeas
 Since $ \vr_0 \big( \D X^j_{\wt{\o}} (\wh{\o}) \big)
 \le  \b1_{\big\{\D X^j_{\wt{\o}} (\wh{\o}) \le \d^{1/2}\big\}} \vr_0 \big(\d^{1/2}\big)
     \n + \n  \b1_{\big\{\D X^j_{\wt{\o}} (\wh{\o}) > \d^{1/2}\big\}} \k \d^{-1/2}  \Big(\D X^j_{\wt{\o}} (\wh{\o})
      \n + \n  \big(\D X^j_{\wt{\o}} (\wh{\o})\big)^{\varpi+1} \Big) $,
      \eqref{eq:xxx151} shows that
 \bea
 \hE_s \big[ \Xi^{s,\wt{\o}}_j \big]
 & \tn \le & \tn  \hE_s \Big[ \big( R^{s,\o  \otimes_t  \cX  (\wt{\o}  )}
    (\z_{\wt{\o}},\wp )\big)  ( \cX^j   ) \Big] + \int_t^s \n  g^{t,\o}_r \big(\cX(\wt{\o})\big) dr
    + (1 \n + \n T ) \vr_0 \big(\d^{1/2}\big) \nonumber \\
  & \tn  & \tn    + (1 \n + \n T )    \k \d^{-1/2} \big(C_1 T \|\o  \n \otimes_t \n  \cX(\wt{\o})
     \n -  \n  \o  \n \otimes_t \n  \wt{\o}_j \|_{0,s}
     \n +  \n  C_{\varpi+1} T^{\varpi+1} \|\o  \n \otimes_t \n    \cX(\wt{\o})
       \n -  \n  \o  \n \otimes_t \n  \wt{\o}_j \|^{\varpi+1}_{0,s}  \big) . \label{eq:uxu219}
 \eea

  Set $\wh{\vr}_0(\d)  \n := \n  \d  \n +  \n  (1 \n + \n T ) \vr_0 \big(\d^{1/2}\big)
    \n +  \n (1 \n + \n T ) \k  \big(C_1 T  \d^{1/2}
     \n +  \n  C_{\varpi+1} T^{\varpi+1}  \d^{\varpi+1/2}  \big)$.
    As $ \wt{\o}  \n \in \n  \cA^\cX_j  \n = \n  \cX^{-1} (\cA_j) $, i.e. $ \cX(\wt{\o})  \n  \in \n  \cA_j
    \n \subset \n  O^s_{\d_j} (\wt{\o}_j) $, one has
    $\|\o  \n \otimes_t  \n   \cX(\wt{\o})  \n - \n  \o  \n \otimes_t \n  \wt{\o}_j \|_{0,s}
    \n = \n  \|\cX(\wt{\o})  \n - \n    \wt{\o}_j \|_{t,s}  \n < \n  \d_j  \n \le \n  \d $.
    It follows from \eqref{eq:uxu219} that
   \bhe  \bea
   \hE_s \big[ \Xi^{s,\wt{\o}}_j  \big]
  & \tn \le  & \tn   \hE_{\fp^j} \Big[ R^{s,\o  \otimes_t  \cX  (\wt{\o}  )}
    (\z_{\wt{\o}},\wp )  \Big]
   + \int_t^s \n  g^{t,\o}_r \big(\cX(\wt{\o})\big) dr   + \wh{\vr}_0(\d) - \d  \nonumber \\
   & \tn  \le  & \tn    \underset{\vs \in \cT^s (n)}{\sup} \hE_{\hP^j}
   \Big[ R^{s,\o  \otimes_t  \cX  (\wt{\o}  )}  (\vs ,\wp )  \Big]
    + \int_t^s \n  g^{t,\o}_r \big(\cX(\wt{\o})\big) dr  + \wh{\vr}_0(\d) - \d  .  \label{eq:xax149}
 \eea \ehe
 Plugging this back into   \eqref{eq:xax151},  we see from \eqref{eq:xxx849} and (u1) that
 \beas
 \hspace{-3mm}
 \sum^\l_{j=1} \hE_{\wh{\hP}} \big[ \b1_{A \cap \cA_j} R^{t,\o}  \big(\tau, \ga_j   \big) \big]
 & \tn \dn \le & \tn  \dn  \sum^\l_{j=1}  \hE_t \Big[ \b1_{\{\wt{\o} \in \cX^{-1} (  A   ) \cap \cX^{-1} ( \cA_j )    \}}
  \Big( \, \underset{\vs \in \cT^s (n)}{\sup} \hE_{\hP^j}
   \Big[ R^{s,\o  \otimes_t  \cX  (\wt{\o}  )}  (\vs ,\wp )  \Big]
   \+ \int_t^s \n  g^{t,\o}_r  \big(\cX(\wt{\o})\big)  dr  \+ \wh{\vr}_0(\d) \- \d  \Big) \Big] \+ \d  \\
    & \tn  \dn  =  & \tn  \dn
   \sum^\l_{j=1}   \hE_\fp \Big[ \b1_{\{\wt{\o} \in    A     \cap \cA_j   \}}
  \Big( \, \underset{\vs \in \cT^s (n)}{\sup} \hE_{\hP^j}
   \Big[ R^{s,\o  \otimes_t   \wt{\o}   }  (\vs ,\wp )  \Big]
   + \int_t^s \n  g^{t,\o}_r   (\wt{\o})   dr \+ \wh{\vr}_0(\d) \- \d  \Big) \Big]  \+ \d   \\
       & \tn  \dn  =  & \tn  \dn
   \sum^\l_{j=1}   \hE_\hP \Big[ \b1_{\{\wt{\o} \in    A     \cap \cA_j   \}}
  \Big( \, \underset{\vs \in \cT^s (n)}{\sup} \hE_{\hP^j}
   \Big[ R^{s,\o  \otimes_t   \wt{\o}   }  (\vs ,\wp )  \Big]
   + \int_t^s \n  g^{t,\o}_r   (\wt{\o})   dr   \Big) \Big]
   \+ \hP  ( A \n \cap \n  \cA^c_0  ) (\wh{\vr}_0(\d) \- \d) \+ \d   .
 \eeas
 In the last equality,    we used the fact that the mapping $ \wt{\o} \to \underset{\vs \in \cT^s (n) }{\sup}
   \hE_{\hP^j}   \Big[ R^{s,\o \otimes_t     \wt{\o} } (\vs , \wp ) \Big]$  is continuous
    under   norm    $\|~\|_{t,T}$ and thus $\cF^t_T-$measurable by Remark \ref{rem_P3} (2).
 Therefore,    \eqref{eq:xxx617} holds for $\wp^n_j = \ga_j$, $j \= 1, \cds \n \l$.

   \no {\bf 3)}     {\it In this part, we still use    \eqref{eq:da025}   and
  the estimate \eqref{eq:xxx151}  to show that  $\{\cP(t,\o)\}_{(t,\o) \in [0,T] \times \O }$
  satisfies  Assumption \ref{assum_V_conti}. }

  Fix $n \n \in \n  \hN  \n \cup \n  \{\infty\}$, $t  \n \in \n  [0,T]$, $\o, \o'  \n \in \n  \O$,
  $\mu  \n \in \n  \sU_t$ and   set $\d  \n := \n  \|\o' \n - \n \o\|_{0,t}$.
  We still take the notation \eqref{eq:xax153} and set
  $(   \hP',\fp',  \cX',    \cW',   \fF'_\cd   ) :=
      \big( \hP^{t,\o',\mu},\fp^{t,\o',\mu},   X^{t,\o',\mu} ,   \wt{W}^{t,\o',\mu},   \fF^{t,\o',\mu}_\cd \big) $.

 Fix $\e > 0$. We still define $\xi_m$'s as in \eqref{eq:xax155} and can find a $\fk \n \in \n \hN$ such that
 $\hE_{ \hP'} [ \xi_{\fk} ]  \n \le \n  \e / 2 $
 and  $ \Phi \big(\|\o'\|_{0,t}, 2^{-\fk} \big)  \n \le \n  \e / 2 $.
 Also,  fix $ \ga  \n  \in \n  \cT^t   $ and $\tau \n \in \n  \cT^t(n)$.
 Similar to $\wh{\tau} \n = \n \tau \big(\wh{\cX}\big)$  in part 2c), 
 $\tau  (  \cX' ) $ belongs to $\ol{\cT}^t $;  and   analogous to
 $\wh{\ga}_j  \n = \n  \nu_j \big( \wh{\cW} \big)$, \eqref{eq:xax047}
 implies that $ \wt{\tau}  \n : = \n  \tau \big(     \cX'  (  \cW  )  \big) $
 is a $\fF-$stopping time. Symmetrically, $\ga (\cX)$ belongs to $\ol{\cT}^t $ and
 $ \wt{\ga} \n := \n   \ga \big(     \cX  (  \cW '  )  \big)  $  defines a $\fF'-$stopping time.

 Set    $t_i    \n  :=  \n   t   \n  \vee  \n   (i 2^{-\fk} T)     $,
      $i  \n   =    \n    0 , \cds   \n   , 2^{\fk}  $.
 Then      $  \wt{\ga}'_\fk  \n := \n  \sum^{2^{\fk}}_{i=0} \b1_{\{ t_{i-1}  <  \wt{\ga}
    \le    t_i  \}} \, t_i$ defines a $ \fF'-$stopping time,
    where $t_{-1}  \n := \n  -1 $.
    By similar arguments to  those that lead to \eqref{eq:xax143},
    one can construct a $\cT^t-$stopping time $\wt{\ga}_\fk$
     valued in $\{t_i\}^{2^{\fk}}_{i=0}$ such that $\wt{\ga}'_\fk  \n = \n  \wt{\ga}_\fk$, $\fp'-$a.s.
     Analogous to \eqref{eq:xax127}, we can deduce that  for $\hP^t_0-$a.s.  $\wt{\o} \n \in \n \O^t$,
\beas
R^{t,\o}  \big(\tau, \wt{\ga}'_\fk   \big)   \big(\cX' (\wt{\o}) \big)
\n - \n  R^{t,\o}  \big(\tau, \wt{\ga}    \big)   \big(\cX' (\wt{\o}) \big)
\le \n  \xi_\fk \big(\cX' (\wt{\o}) \big)   \n + \n  \vr_0 \big( 2^{-\fk} + \eta' (\wt{\o}) \big) ,
\eeas
 where $ \eta'  := \underset{r \in [ \ga (\cX)   , (\ga (\cX)   +2^{-\fk}) \land T]}{\sup}
 \big|  \cX'_r     \n  -  \n  \cX'_{ \ga (\cX) }   \big| $.
 And similar to \eqref{eq:xax157}, \eqref{eq:xxx153} implies that
 \bea
 \hE_{\fp'} \big[ R^{t,\o}  \big(\tau, \wt{\ga}_\fk   \big) \n - \n  R^{t,\o}  \big(\tau, \wt{\ga}    \big)   \big]
 & \tn = & \tn \hE_{\fp'} \big[ R^{t,\o}  \big(\tau, \wt{\ga}'_\fk   \big)
 \n - \n  R^{t,\o}  \big(\tau, \wt{\ga}    \big)   \big]
 \n = \n \hE_t \big[ R^{t,\o}  \big(\tau, \wt{\ga}'_\fk   \big) (\cX')
 \n - \n  R^{t,\o}  \big(\tau, \wt{\ga}    \big) (\cX')  \big]
 \nonumber   \\
  & \tn \le & \tn  \hE_t \big[ \xi_\fk (\cX') + \vr_0 \big( 2^{-\fk} \+ \eta'   \big) \big]
   \le  \hE_{\hP'} [\xi_\fk] + \Phi \big(\|\o'\|_{0,t}, 2^{-\fk} \big) \le \e .  \label{eq:xax159}
 \eea

 Since \eqref{eq:xax069} shows that
 $ \tau \big( \cX'(\wt{\o})\big) \n = \n  \tau \big( \cX'\big( \cW ( \cX (\wt{\o}))\big)\big)
  \n = \n  \wt{\tau} \big( \cX (\wt{\o})\big) $
 and $ \wt{\ga} ( \cX'(\wt{\o}))  \n = \n   \ga \big(     \cX  \big(  \cW ' ( \cX'(\wt{\o})) \big)  \big)
  \n = \n  \ga \big(     \cX ( \wt{\o} ) \big) $ hold for $\hP^t_0-$a.s. $\wt{\o} \in \O^t$,
 we see from    \eqref{eq:da021} and  \eqref{eq:da025}    that
  for $\hP^t_0-$a.s. $\wt{\o} \in \O^t$
 \beas
 && \hspace{-1cm} \big( R^{ t , \o' }  ( \tau  , \wt{\ga}  ) \big)  \big( \cX'(\wt{\o})\big)
  \n - \n  \big( R^{t,\o} (\wt{\tau} , \ga) \big) \big(\cX(\wt{\o})\big)
   \n = \n   R \big( t,  \tau ( \cX'(\wt{\o})), \wt{\ga} ( \cX'(\wt{\o})), \o'  \n \otimes_t \n  \cX'(\wt{\o})  \big)
  \n - \n R \big( t,  \wt{\tau} ( \cX(\wt{\o})),  \ga  ( \cX(\wt{\o})), \o  \n \otimes_t \n  \cX(\wt{\o})  \big)     \\
  && = \n  R \big( t,  \wt{\tau} ( \cX (\wt{\o})),  \ga  ( \cX (\wt{\o})), \o'  \n \otimes_t \n  \cX'(\wt{\o})  \big)
  \n - \n R \big( t,  \wt{\tau} ( \cX(\wt{\o})),  \ga  ( \cX(\wt{\o})), \o  \n \otimes_t \n  \cX(\wt{\o})  \big)   \\
 &&  \le (1 \n + \n T )   \vr_0 \big( \|\o'  \n \otimes_t \n  \cX'(\wt{\o})
   \n - \n  \o  \n \otimes_t \n  \cX(\wt{\o}) \|_{0, T }  \big)
    \ls  (1 \n + \n T) \vr_0 \big( \| \o'  \n - \n  \o \|_{0,t}
    \n + \n  \| \cX' (\wt{\o})  \n -  \n  \cX (\wt{\o}) \|_{t,T} \big)
    \n = \n (1 \n + \n T) \vr_0 \big( \d  \n + \n  \D X (\wt{\o}) \big)  \\
 &&    \n \le  \n  \b1_{\{\D X (\wt{\o}) \le  \d^{1/2} \}} (1 \n + \n T) \vr_0 \big( \d  \n + \n  \d^{1/2} \big)
      \n + \n   \b1_{\{ \D X (\wt{\o}) > \d^{1/2} \}} \k (1 \n + \n T)
  \d^{-1/2} \big( ( 1  \n + \n   2^{\varpi-1} \d^\varpi ) \D X (\wt{\o})
  \n +  \n  2^{\varpi-1} (\D X (\wt{\o}) )^{\varpi+1} \big) ,
 \eeas
 with $ \D X (\wt{\o}) := \| \cX' (\wt{\o}) \n -  \n  \cX (\wt{\o}) \|_{t,T} $.
 Then  \eqref{eq:xax159} and \eqref{eq:xxx151} show  that
 \bea
  \hE_{\hP'} \big[ R^{t,\o}  \big(\tau, \wt{\ga}_\fk   \big) \big]
  & \tn =  & \tn  \hE_{\fp'} \big[ R^{t,\o}  \big(\tau, \wt{\ga}_\fk   \big) \big]
  \le   \hE_{\fp'}  \big[  R^{t,\o'} (\tau, \wt{\ga})    \big] \n + \n \e   \n  =    \n
      \hE_t \big[ \big( R^{t,\o'} (\tau, \wt{\ga}) \big) ( \cX' ) \big] \n + \n \e     \nonumber  \\
    & \tn  \le  & \tn  \hE_t \big[ \big( R^{t,\o} (\wt{\tau} , \ga) \big) (\cX)\big]  \n + \n   \vr_1  (\d) \n + \n \e
    \n = \n   \hE_{\fp} \big[  R^{t,\o} (\wt{\tau} , \ga)  \big]
     \n + \n   \vr_1  (\d)  \n + \n \e  ,  \label{eq:xxx751}
 \eea
  where   $\vr_1(\d)  \n : = \n     (1 \n + \n T)   \vr_0  (  \d  \n + \n    \d^{1/2}    )
   \n + \n  \k (1 \n + \n T) \big( ( 1 \n + \n
  2^{ \varpi-1 } \d^{\varpi} )  C_{1 } T \d^{1/2}
   \n + \n   2^{ \varpi-1 }   C_{\varpi+1} T^{\varpi+1} \d^{ \varpi +1  /2} \big) \n \ge \n \vr_0 (\d)  $.

  Similar to \eqref{eq:xax149}, one can deduce that
 $  \hE_{\fp} \big[  R^{t,\o} (\wt{\tau} , \ga)  \big]
     \n  \le  \n    \underset{\vs \in \cT^t (n)}{\sup} \hE_{\hP }
   \big[ R^{t,\o }  (\vs ,\ga )  \big]   $.
 So it  follows from \eqref{eq:xxx751} that
 \beas
  \hE_{\hP'}  \big[  R^{t,\o'} (\tau, \wt{\ga}_\fk)    \big]
  \le  \underset{\vs \in \cT^t (n)}{\sup} \hE_{\hP }
   \big[ R^{t,\o }  (\vs ,\ga )  \big]  \n + \n   \vr_1  (\d)  \n + \n  \e  .
 \eeas
 Taking supremum over $\t \in \cT^t (n) $ on the left-hand-side yields  that
 \beas
 \underset{\z \in \cT^t}{\inf} \,
  \underset{\tau \in \cT^t (n) }{\sup}  \hE_{\hP'}  \big[  R^{t,\o'} (\tau, \z)    \big]
  \le  \underset{\tau \in \cT^t (n) }{\sup}  \hE_{\hP'}  \big[  R^{t,\o'} (\tau, \wt{\ga}_\fk)    \big]
    \le  \underset{\vs \in \cT^t (n) }{\sup}  \hE_{\hP} \big[  R^{t,\o} (\vs , \ga)  \big]
     \n + \n   \vr_1  (\d)   \n + \n  \e .
 \eeas
 Then  taking infimum over $\ga \in \cT^t$ on the right-hand-side, we obtain   that
 \beas
 \underset{\z \in \cT^t}{\inf} \,
  \underset{\tau \in \cT^t (n) }{\sup}  \hE_{\hP^{t,\o',\mu}}  \big[  R^{t,\o'} (\tau, \z)    \big]
    \le  \underset{\ga  \in \cT^t}{\inf} \,
     \underset{\vs \in \cT^t (n) }{\sup}  \hE_{\hP^{t,\o,\mu}} \big[  R^{t,\o} (\vs , \ga)  \big]
     \n + \n   \vr_1  (\d)   \n + \n  \e .
 \eeas
 Letting $\e \to 0$  and   taking infimum over $\mu \in \sU_t$ on both sides lead to  that
 \beas
 \hspace{-3mm}
 V^n_t(\o') \= \n  \underset{\mu \in \sU_t}{\inf} \,
  \underset{\z  \in \cT^t}{\inf} \, \underset{\t \in \cT^t (n) }{\sup}
   \hE_{\hP^{t,\o',\mu}} \n  \big[  R^{t,\o'} ( \tau, \z )    \big]
  \ls  \underset{\mu \in \sU_t}{\inf} \,  \underset{\ga  \in \cT^t}{\inf} \,
  \underset{\vs \in \cT^t (n) }{\sup} \hE_{\hP^{t,\o,\mu}} \n  \big[  R^{t,\o} (\vs, \ga )  \big]
       \n + \n  \vr_1 \big( \|\o' \dn - \n \o\|_{0,t} \big)
 \=    V^n_t(\o)    \n + \n  \vr_1 \big( \|\o' \dn - \n \o\|_{0,t} \big)   .
\eeas
 Exchanging the roles of $\o'$ and $\o$ shows that $\{\cP(t,\o)\}_{(t,\o) \in [0,T] \times \O }$
  satisfies  \eqref{eq:aa213}.

   \no {\bf 4)}  To verify    Assumption  \ref{assum_V_conti_2} for $\{\cP(t,\o)\}_{(t,\o) \in [0,T] \times \O }$,
  we fix $\a  \n > \n  0$ and $\d  \n \in \n  (0,T]$.

Let $t  \n \in \n  [0,T)$, $\o \n \in \n O^t_\a (\bz) $,  $\mu  \n \in \n  \sU_t$ and $\z   \n   \in  \n    \cT^t$.
We take the notation \eqref{eq:xax153}  again.  Similar to $\wh{\tau} \n = \n \tau \big(\wh{\cX}\big)$  in part 2c),
 $ \wt{\z}  \n := \n  \z  (  \cX  ) $ is a $\ol{\cT}^t -$stopping time.
 Set   $\wt{\eta}  \n : = \n    \underset{r \in \big[\wt{\z}, (\wt{\z}+\d) \land T \big]}{\sup}
   \big|  \cX_r  \n - \n   \cX_{\wt{\z}}    \big|  $. Analogous to \eqref{eq:xax157},
   one can deduce from \eqref{eq:xxx153} that
 \beas
 && \hspace{-1cm} \hE_\hP  \bigg[  \vr_1  \Big(   \d \n + \n   \underset{r \in [\z, (\z+\d) \land T] }{\sup}
  \big|  B^t_r  \n - \n  B^t_\z   \big| \Big) \bigg]
   \n = \n  \hE_\fp  \bigg[  \vr_1  \Big(   \d \n + \n   \underset{r \in [\z, (\z+\d) \land T] }{\sup}
  \big|  B^t_r  \n - \n  B^t_\z   \big| \Big) \bigg]
   \n = \n
   \hE_t  \bigg[  \vr_1  \Big(   \d \n + \n   \underset{r \in \big[\wt{\z}, (\wt{\z}+\d) \land T \big] }{\sup}
  \big|  \cX_r  \n - \n  \cX_{\wt{\z}}   \big| \Big) \bigg]
  \n = \n  \hE_t  \big[  \vr_1   (   \d \n + \n \wt{\eta}  ) \big] \\
  &&   \le  \n  \vr_1   (  \d \n +  \n   \d^{1/4}    )  \n + \n
   \k    (1 \+ 2^{\varpi-1} \d^\varpi) \vf_1 (\|\o\|_{0,t}) \d^{1/4}
 \+ \k 2^{\varpi-1}  \vf_{\varpi+1} (\|\o\|_{0,t}) \d^{\varpi/2+1/4} \le \vr_\a ( \d) ,
 \eeas
 where $  \vr_\a ( \d) : = \vr_1   (  \d \n +  \n   \d^{1/4}    )  \n + \n
   \k    (1 \n + \n 2^{\varpi-1} \d^\varpi) \vf_1 (\a) \d^{1/4}
  \n + \n  \k 2^{\varpi-1}  \vf_{\varpi+1} (\a) \d^{\varpi/2+1/4} $.
 Taking supremum over $ \z   \n   \in  \n    \cT^t $
 and then taking supremum over $ \mu  \n \in \n  \sU_t $ and $\o \n \in \n O^t_\a (\bz) $ yield \eqref{eq:aa213b}. \qed

  \subsection{Proof of Theorem \ref{thm_triplet}}


\if{0}
Since the path-independent probability class
$\{\cP_t\}_{t \in [0,T]}$ satisfies (P1)$-$(P3), Assumption \ref{assum_V_conti} and Assumption \ref{assum_V_conti_2},
we see from Theorem \ref{thm_RDG}, Remark \ref{rem_V_adapted}, Proposition \ref{prop_conti_V},
\eqref{eq:uxu170} and (A$'$) that item (1)$-$(3) hold.

 It remains to   find an optimal triplet for $V_0$:
\fi

 If $  V_0 = L_0  $,    then $\tau_* = 0$ and it thus holds
 for any $(\hP, \ga) \in \cP \times  \cT$ that
 $ \hE_\hP [ R( \tau_* , \ga  ) ] = \hE_\hP [ R(0, \ga  ) ] = \hE_\hP [L_0] = L_0 = V_0$.

 Next, let us assume that $  V_0 > L_0  $.
 Theorem \ref{thm_RDG} (1), Proposition \ref{prop_conti_V} (1),  (A$'$) and the proof of Remark \ref{rem_V_conti}
 imply  that  the process $\sX_t : = V_t - L_t$, $t \in [0,T]$
 has all continuous paths and satisfies
 \beas
 |\sX_t (\o) - \sX_t (\o') | \le  | V_t (\o) - V_t (\o') | +  | L_t (\o) - L_t (\o') |
 \le 2 \rho_0 \big(\|\o-\o'\|_{0,t}\big)   , \q \fa t \in [0,T], ~ \fa  \o,\o' \in \O .
 \eeas
 Then   applying Theorem 3.1   of \cite{RDOSRT} with payoff processes
 $\cL \n := \n -U$, $\cU \n := \n -L$ and random maturity
  $\tau_0 \n = \n  \inf\{t  \n \in \n  [0,T] \n : \sX_t  \n \le \n  0\}  \n \land \n  T
   \n = \n  \inf\{t  \n \in \n  [0,T] \n : V_t  \n = \n  L_t \}  \n = \n  \tau_* $ shows that
   (In particular, (H4) implies (P4) of \cite{RDOSRT}
  by Remark 3.1 (3) therein) for some $(\hP_*, \ga_*)  \n \in \n  \cP  \n \times \n  \cT$,
  $  \underset{(\hP,\ga) \in \cP \times \cT}{\sup} \hE_\hP \big[ \b1_{\{\ga < \tau_*\}} \cL_{\ga} +
  \b1_{\{ \tau_* \le \ga   \}} \cU_{\tau_*} \big] = \hE_{\hP_*} \big[ \b1_{\{\ga_* < \tau_*\}} \cL_{\ga_*} +
  \b1_{\{ \tau_* \le \ga_*   \}} \cU_{\tau_*} \big] $.
 Multiplying $-1$ on both sides, we see from \eqref{eq:xax015} that
 $ V_0 = \underset{(\hP,\ga) \in \cP \times \cT}{\inf} \,
 \hE_\hP [ R (\tau_*, \ga) ]  =  \hE_{\hP_*} \n  \big[ R (\tau_*,  \ga_*) \big] $.
 \qed

\appendix
\renewcommand{\thesection}{A}
\refstepcounter{section}
\makeatletter
\renewcommand{\theequation}{\thesection.\@arabic\c@equation}
\makeatother

  \section{Appendix} 

 \subsection{A Technical Lemma}

 \begin{lemm} \label{lem_submg0}
 Define $\ol{\U}_t  \df   \ol{V}_t \+ \int_0^t g_r dr $, $ t \in [0,T]$
 Given $\z  \n \in \n  \cT$,
 it holds for any $(t,\o)  \n \in \n  [0,T]  \n \times \n  \O$ and $ t'   \n \in \n  [0,t]  $    that
  \bea  \label{eq:dc251}
   \ol{\U}_t (\o) \le  \underset{\hP \in \cP(t , \o)}{\inf} \,    \hE_\hP \bigg[ \Big( \,
 \ol{\U}_{   \big( \tau^*_{(t',\o)} (\Pi^0_{t'})     \land    \z  \big) \vee t  }  \Big)^{t, \o} \bigg]  .
 \eea
 \end{lemm}

 \no {\bf Proof of Lemma \ref{lem_submg0}:}
 Fix    $0 \le t' \le t \le T  $, $\o \in \O$,   $\z  \n \in \n  \cT$
 and set $\a  \n : = \n  1  \n + \n  \|\o   \|_{0 , T}  $.


\no {\bf 1)}   When $t \= T$, one has
  $  \underset{\hP \in \cP(T , \o)}{\inf} \,    \hE_\hP \bigg[ \Big( \,
 \ol{\U}_{   \big( \tau^*_{(t',\o)} (\Pi^0_{t'})     \land    \z  \big) \vee T  }  \Big)^{T, \o} \bigg]
  \=  \underset{\hP \in \cP(T , \o)}{\inf} \,    \hE_\hP \Big[   (\ol{\U}_T)^{T, \o} \Big]
  \=  \underset{\hP \in \cP(T , \o)}{\inf} \,    \hE_\hP \Big[   \ol{\U}_T (\o) \Big]
  \=   \ol{\U}_T (\o)     $.

   \no {\bf 2)} Next, suppose that $t \< T$ and $\ol{V}_t(\o) \= L_t(\o) $. Then
 \beas
 \tau^*_{(t',\o)} \big(\Pi^0_{t'} (\o)\big) \n   = \n     \inf\Big\{ s  \n \in \n  [t',T] \n :
      \ol{V}^{t',\o}_s  \big(\Pi^0_{t'} (\o)\big) \n = \n  L^{t',\o}_s \big(\Pi^0_{t'} (\o)\big)      \Big\}
       \n   = \n     \inf\big\{ s  \n \in \n  [t',T] \n :
      \ol{V}_s    (\o)  \n = \n  L_s  (\o)       \big\} \le t ,
 \eeas
 which means that  $\o \ins \big(\Pi^0_{t'}\big)^{-1} (A') $ with
 $A' \df \big\{ \o' \ins \O^{t'} \n :  \tau^*_{(t',\o)} (\o') \ls t \big\} \ins \cF^{t'}_t $.
 Since Lemma A.1 of \cite{ROSVU} shows that
 \bea  \label{eq:xax003}
  \hb{$\Pi^0_{t'}$ is an $\cF_r / \cF^{t'}_r-$measurable mapping, }  \q \fa r  \n \in \n  [t',T] ,
   \eea
 we see  that $\big(\Pi^0_{t'}\big)^{-1} (A') \ins  \cF_t$. It follows from  Lemma \ref{lem_element}   that
 \bea \label{eq:uxu203}
 \o \otimes_t \O^t \subset \big(\Pi^0_{t'}\big)^{-1} (A') \q \hb{or} \q
 \tau^*_{(t',\o)} \big( \Pi^0_{t'} (\o \otimes_t  \wt{\o}) \big) \le t , \q \fa \wt{\o} \in \O^t .
 \eea

 Remark \ref{rem_V_adapted}  and Proposition \ref{prop_conti_V} (1) show that
   $ \ol{\U} $ is an $\bF-$adapted process with all continuous paths.
   Applying \eqref{eq:bb421} to $\ol{\U}_t \ins  \cF_t$ and using \eqref{eq:uxu203} yield that
\beas
  \Big( \,
 \ol{\U}_{   \big( \tau^*_{(t',\o)} (\Pi^0_{t'})     \land    \z  \big) \vee t  }  \Big)^{t, \o} (\wt{\o})
 \= \ol{\U} \Big( \big( \tau^*_{(t',\o)} (\Pi^0_{t'} (\o \otimes_t \wt{\o})) \n \land  \n  \z (\o \otimes_t \wt{\o}) \big)
  \n \vee \n  t ,  \o \otimes_t \wt{\o} \Big) \= \ol{\U}_t (\o \otimes_t \wt{\o})
 \=  \ol{\U}_t (\o )   , \q \fa \wt{\o} \ins \O^t ,
 \eeas
 Thus we still obtain \eqref{eq:dc251} as an equality.

 \no {\bf 3)}   The discussion of the case $t \< T$ with $\ol{V}_t(\o)  \n > \n  L_t(\o) $
 is relatively lengthy. We split it  into several steps.
 Since $\lmtu{n \to \infty} V^n_t (\o)  \n = \n  \ol{V}_t (\o)$
 by \eqref{eq:da371} and Proposition \ref{prop_grid},
 there exists an integer $N \n = \n N(t,\o) \n > \n   \log_2 \big( \frac{T}{T   -   t} \big)$ such that
 $ V^n_t (\o)  \n > \n   L_t(\o)  $ for any $n  \n \ge \n  N$.

  Fix $\d  \n > \n  0$ and    $k , n  \n \in \n  \hN$  with   $k  \n \ge \n  n    \n > \n  N $.
  For any $r \ins [t',T]$,
  as $A_r  \n := \n  \big\{\wt{\o}  \n \in \n  \O^{t'} \n : \tau^{n,\d}_{(t',\o)} (\wt{\o})  \n < \n  r \big\}
  \n \in \n  \cF^{t'}_r $,   \eqref{eq:xax003} implies that
 $ \Big\{\o' \ins \O \n : \tau^{n,\d}_{(t',\o)} (\Pi^0_{t'} (\o')) \< r \Big\}
 = \big\{\o' \ins \O \n :  \Pi^0_{t'} (\o') \ins A_r \big\} \= (\Pi^0_{t'})^{-1} (A_r) \ins \cF_r $.
 So   $\tau^{n,\d}_{(t',\o)} (\Pi^0_{t'})$ is an $\bF-$optional time valued in $[t',T]$, and
 it follows that $\nu^{n,\d}  \n : = \n    \big( \tau^{n,\d}_{(t',\o)} (\Pi^0_{t'})  \n  \land  \n  \z  \big)
  \n \vee \n  t $ is  an $\bF-$optional time valued in $[t,T]$.

   Let $i_k$ be the largest integer such that $ i_k 2^{-k}T  \n \le \n  t $.
 As $k  \n >   \log_2 \big( \frac{T}{T   -   t} \big)  $, one can deduce  that $i_k  \n < \n  2^k  \n - \n  1 $.
  Set $t_{i_k} \df t$ and $ t_i   \df   i 2^{-k}T $
  for   $i =i_k \+ 1, \cds, 2^k   $.

 \no {\bf 3a)} {\it In the first step, we derive from  Proposition \ref{prop_DPP}   an auxiliary inequality:}
 \bea \label{eq:da411}
 V^n_{t_i} (\o) \le L_{t_i} (\o) \vee
   \ul{\sE}_{t_i}   \bigg[ V^n_{t_{i+1}}
    \n + \n \int_{t_i}^{t_{i+1}} g_r dr  \bigg] (\o) , \q  i= i_k , \cds, 2^k -1    .
 \eea

 Let $i \= i_k , \cds, 2^k \- 1 $.   Applying \eqref{eq:bb013}  with $(t,s) \= (t_i,t_{i+1})$
 and taking $\ga \= t_{i+1}$ yield  that
 \bea
  V^n_{t_i} (\o)
  \ls  \underset{\hP \in \cP  ( t_i , \, \o  )}{\inf} \,
  \underset{\t \in \cT^{t_i} \n  (n)}{\sup} \,
   \hE_\hP \Bigg[ \b1_{\{\tau   < t_{i+1} \}}   R^{t_i , \,  \o}   (\tau, t_{i+1})
   + \b1_{\{\tau   \ge {t_{i+1}}\}}
   \bigg( \big( V^n_{t_{i+1}}  \big)^{t_i \n , \, \o}
   \n + \n \int_{t_i}^{t_{i+1}} g^{t_i \n , \, \o}_r dr  \bigg)     \Bigg]  \, .      \label{eq:ae015}
 \eea
  For any $\t  \n \in \n  \cT^{t_i}    (n) $, it    takes values in
  $\{t_i\} \cup \{  j 2^{-n} T \}^{2^n}_{j=j_0}$, where $j_0$ is the smallest integer such that
  $ t_i \< j_0 2^{-n} T $.
  As $n \ls k$, one has   $ t_{i+1} \ls j_0 2^{-n} T  $,
  so $\{\tau    \n < \n  t_{i+1} \} \n  = \n  \{\tau   \n = \n  t_i \}  \n \in \n
 \cF^{\raisebox{0.2ex}{\scriptsize $t_i$}}_{\raisebox{-0.2ex}{\scriptsize $t_i$}}
 \= \{ \es, \O^{t_i} \}$.
 To wit,  we have either $ \{\tau   \< t_{i+1} \}  \= \{\tau  \= t_i \}
 \= \O^{t_i}  $
 or $ \{\tau   \gs  t_{i+1} \}  \= \O^{t_i}  $. Since $ R^{t_i , \,  \o}   (t_i, t_{i+1})
 \= L^{t_i , \,  \o}_{t_i}  \=  L (t_i , \,  \o )   $
 by \eqref{eq:uxu205},
 we see from \eqref{eq:ae015} that
 \beas
 \q    V^n_{t_i} (\o)     \ls    \underset{\hP \in \cP  ( t_i , \, \o  )}{\inf}   \bigg(  L_{t_i} (\o)  \n \vee  \n
    \hE_\hP \bigg[  \big( V^n_{t_{i+1}}  \big)^{t_i \n , \, \o}
    \n + \n \int_{t_i}^{t_{i+1}} g^{t_i \n , \, \o}_r dr  \bigg] \bigg)
    \n = \n  L_{t_i} (\o)  \n \vee    \ul{\sE}_{t_i}
    \bigg[ V^n_{t_{i+1}}
    \n + \n \int_{t_i}^{t_{i+1}} g_r dr  \bigg] (\o) , ~ \;  \hb{proving \eqref{eq:da411}. }
 \eeas

\no {\bf 3b)} {\it In the next step, we will show that over time grids $\{t_i\}^{2^k}_{i_k}$,
the $\bF-$adapted process $\U^n_t   \n := \n
   V^n_t  \n + \n  \int_0^t    g_r dr $, $ t  \n \in \n  [0,T]$    
   is an $\ul{\sE}-$submartingale    up to time
   $    \nu^{n,\d}_k    \n := \n
 \sum^{2^k}_{i = i_k + 1} \b1_{\{ t_{i-1}  \le  \nu^{n,\d} < t_i \}} t_i
 \n + \n  \b1_{\{ \nu^{n,\d}  =  T\}} T   $, i.e.    }
    \bea    \label{eq:db014}
 \U^n_{\nu^{n,\d}_k \land t_i}   (\o) \le     \ul{\sE}_{t_i}
   \Big[ \U^n_{\nu^{n,\d}_k \land t_{i+1}}  \Big] (\o) , \q
    i = i_k , \cds , 2^k - 1   .
 \eea

  For any $r  \n \in \n  \big[ t_{i^k+1},T \big)$,  let
 $j_r$ be the largest integer such that $ t_{j_r}    \n \le \n  r $. Since
 $\nu^{n,\d}   $ is  an $\bF-$optional time, one can deduce that
 $\{ \nu^{n,\d}_k  \n \le \n  r\}  \n = \n  \underset{i=i_k+1}{\overset{j_r}{\cup}}
 \{ \nu^{n,\d}_k  \n = \n  t_i  \}  \n = \n
  \underset{i=i_k+1}{\overset{j_r}{\cup}} \{ t_{i-1}   \n \le  \n  \nu^{n,\d}  \n < \n  t_i \}
  \n = \n   \{ \nu^{n,\d}  \n < \n  t_{j_r} \}     \n \in \n  \cF_{ t_{j_r}  }    \n \subset \n  \cF_r  $.
 So $ \nu^{n,\d}_k $ is a $ \cT_t (k) -$stopping time.

 \no {\bf (i)}
   Let   $i \= i_k$ first. We simply denote $t_{i_k+1}$ by $s$.
   Since $  V^n_t (\o)  \n > \n   L_t(\o)   $,
   applying \eqref{eq:da411} with $i  \n = \n  i_k$ yields that
   \bea \label{eq:db011b}
     V^n_{t} (\o)  \n \le  \n   \ul{\sE}_{t}   \bigg[ V^n_{s}
    \n + \n \int_{t}^{s} g_r dr  \bigg] (\o)  .
    \eea
    As $ \nu^{n,\d}_k \gs t_{i_k+1} \= s > t_{i_k } \= t $,
    the first equality in \eqref{eq:xax141} shows   that
    \beas
     \big( \U^n_{\nu^{n,\d}_k \land s} \big)^{t, \o} (\wt{\o})
  & \tn =  & \tn   \U^n \big( \nu^{n,\d}_k ( \o  \n \otimes_{t} \n  \wt{\o} )   \n \land \n  s ,
   \o  \n \otimes_{t} \n  \wt{\o} \big)
    \n = \n  \U^n \big(   s ,  \o  \n \otimes_{t} \n  \wt{\o} \big)
     \n = \n    V^n \big(   s ,  \o  \n \otimes_{t} \n  \wt{\o} \big)
    \n + \n \int_{t}^{s}  \n g_r   \big(     \o  \n \otimes_{t} \n  \wt{\o} \big)   dr
    \n + \n \int_0^{t}  \n g_r   \big(     \o   \big)   dr   \nonumber  \\
      & \tn  =  & \tn       \bigg( V^n_{s}
    \n + \n \int_{t}^{s}  \n g_r dr \bigg)^{t, \o} (\wt{\o})
    \n + \n \int_0^{t}  \n g_r   \big(     \o   \big)   dr     ,
    ~ \q \fa \wt{\o}  \n \in \n  \O^{t} .
    \eeas
    Taking expectation $\hE_\hP [~]$ and then  taking infimum over $\hP \in \cP(t,\o)$,
    we see from \eqref{eq:db011b}    that
    \beas
     \ul{\sE}_{t}   \Big[ \U^n_{\nu^{n,\d}_k \land s}  \Big] (\o)
    \=   \ul{\sE}_{t}   \bigg[ V^n_{s}
    \n + \n \int_{t}^{s} g_r dr  \bigg] (\o) \+ \int_0^{t} g_r (\o) dr
    \gs   \U^n_{  t}   (\o) \= \U^n_{\nu^{n,\d}_k \land t}   (\o)   ,
    \q   \hb{proving \eqref{eq:db014} for $   i = i_k$.} 
    \eeas

 \no {\bf (ii)}   Next, let   $   i  \n = \n  i_k  \n + \n  1 , \cds  \n , 2^k  \n - \n  1 $.
   Given $\o \ins \{\nu^{n,\d}_k      \n  \le \n  t_i\}$,
  applying Proposition \ref{prop_shift0} (3) with $(t,s,\tau) = \big(0,t_i, \nu^{n,\d}_k \big)$ shows that
   $      \nu^{n,\d}_k  \big(\o  \n \otimes_{t_i} \n  \O^{t_i}\big)  \n \equiv \n
   \nu^{n,\d}_k   (\o) \df  \wh{t}    $.
   As $ \U^n_{\wh{t}}  \n \in \n  \cF_{\wh{t}}  \n  \subset \n  \cF_{t_i}   $,  using \eqref{eq:bb421}
   with $(t,s,\eta) \= \big(0,t_i,  \U^n_{\wh{t}} \big) $  yields that for any $\wt{\o} \in    \O^{t_i}$
   $    \Big( \U^n_{\nu^{n,\d}_k \land t_{i+1} } \Big)^{t_i,\o} ( \wt{\o} )
   \= \U^n \big(    \nu^{n,\d}_k (\o \oti_{t_i}  \wt{\o}) \land t_{i+1}  , \o \oti_t \wt{\o} \big)
   \= \U^n \big(  \wh{t} \land t_{i+1} ,  \o \oti_{t_i} \wt{\o} \big)
   \=   \U^n \big(  \wh{t}   ,  \o \oti_{t_i} \wt{\o} \big) \= \U^n \big(  \wh{t},      \o \big)
   $.
   It follows that
   \bea   \label{eq:db017}
   \ul{\sE}_{t_i} \big[ \,  \U^n_{\nu^{n,\d}_k \land t_{i+1}  }  \big] (\o)
   \n = \n  \underset{\hP \in \cP(t_i,\o)}{\inf} \hE_\hP
   \Big[ \big( \U^n_{\nu^{n,\d}_k \land t_{i+1}  } \big)^{t_i,\o}  \Big]
    \n = \n  \underset{\hP \in \cP(t_i,\o)}{\inf} \hE_\hP \big[ \, \U^n (  \wh{t},      \o  ) \big]
    \n = \n  \U^n \big(  \wh{t},      \o \big)  \n = \n  \U^n \big(   \nu^{n,\d}_k (\o) \n \land \n  t_i  ,  \o \big) . \q
   \eea

 Then we let $\o \ins \{\nu^{n,\d}_k   \> t_i\}$.
      Proposition \ref{prop_shift0} (3)   shows that
      \bea    \label{eq:cc601}
      \o  \n \otimes_{t_i} \n  \O^{t_i}  \subset \big\{ \nu^{n,\d}_k       >   t_i   \big\}
        =    \big\{ \nu^{n,\d}_k       \ge   t_{i+1}    \big\}  ,
      \eea
     and  one can deduce that
     $ \nu^{n,\d}  (\o   ) \n \ge \n  t_i  \n \ge \n  t_{i_k+1}  \n > \n  t_{i_k}  \n = \n   t$.
     By the definition of $\nu^{n,\d}$, one has
     $ t_i \n \le  \n \nu^{n,\d}  (\o   )  \n = \n  \tau^{n,\d}_{(t',\o)} (\Pi^0_{t'} (\o)) \land    \z (\o)
      \n \le \n  \tau^{n,\d}_{(t',\o)} \big(\Pi^0_{t'} (\o)\big) $ and it follows that
      \bhe \bea \label{eq:uxu211}
      V^n ( t_i , \o)  \gs  L   ( t_i , \o) \+ \d    .
      \eea \ehe
 This together with \eqref{eq:da411} shows   that
 $    V^n_{t_i} (\o) \n  \le   \n
   \ul{\sE}_{t_i}   \Big[ V^n_{t_{i+1}}
    \n + \n \int_{t_i}^{t_{i+1}}  \n  g_r dr  \Big] (\o)   $.
   Adding   $\int_0^{t_i}  \n  g_r(\o) dr$ to both sides, one can deduce from \eqref{eq:cc601} that
  \if{0}
    Adding   $\int_0^{t_i}  \n  g_r(\o) dr$ to both sides yields that
   \bea  \label{eq:db011}
   \U^n_{\nu^{n,\d}_k \land t_i}   (\o)  =  \U^n_{  t_i}   (\o)  =
   V^n_{t_i} (\o) + \int_0^{t_i} g_r (\o) dr
  \le     \ul{\sE}_{t_i}
   \bigg[ V^n_{t_{i+1}}
    \n + \n \int_{t_i}^{t_{i+1}} g_r dr  \bigg] (\o) + \int_0^{t_i} g_r (\o) dr .
    \eea
    For any $\hP \in \cP(t_i,\o)$, since \eqref{eq:bb421} implies that
    $g_r \big(\o \otimes_{t_i} \O^{t_i} \big) = g_r (\o)$, $\fa r \in [0, t_i]$, we can deduce from
    \eqref{eq:cc601} that
    \bea
     \big( \U^n_{\nu^{n,\d}_k \land t_{i+1}} \big)^{t_i, \o} (\wt{\o})
  & \tn =  & \tn   \U^n \big( \nu^{n,\d}_k ( \o  \n \otimes_{t_i} \n  \wt{\o} )   \n \land \n  t_{i+1} ,
   \o  \n \otimes_{t_i} \n  \wt{\o} \big)
    \n = \n  \U^n \big(   t_{i+1} ,  \o  \n \otimes_{t_i} \n  \wt{\o} \big)
     \n = \n  V^n \big(   t_{i+1} ,  \o  \n \otimes_{t_i} \n  \wt{\o} \big)
    \n + \n \int_0^{t_{i+1}} \n  g_r   \big(     \o  \n \otimes_{t_i} \n  \wt{\o} \big)   dr  \nonumber \\
   & \tn  =  & \tn   V^n \big(   t_{i+1} ,  \o  \n \otimes_{t_i} \n  \wt{\o} \big)
    \n + \n \int_{t_i}^{t_{i+1}}  \n g_r   \big(     \o  \n \otimes_{t_i} \n  \wt{\o} \big)   dr
    \n + \n \int_0^{t_i}  \n g_r   \big(     \o   \big)   dr   \nonumber  \\
      & \tn  =  & \tn       \bigg( V^n_{t_{i+1}}
    \n + \n \int_{t_i}^{t_{i+1}}  \n g_r dr \bigg)^{t_i, \o} (\wt{\o})
    \n + \n \int_0^{t_i}  \n g_r   \big(     \o   \big)   dr     ,
    ~ \q \fa \wt{\o}  \n \in \n  \O^{t_i} .   \label{eq:ek173}
    \eea
    Taking expectation $\hE_\hP [~]$ and   taking infimum over $\hP \in \cP(t_i,\o)$,
    we see from \eqref{eq:db011}    that
    \fi
    \beas
     \ul{\sE}_{t_i}   \Big[ \U^n_{\nu^{n,\d}_k \land t_{i+1}}  \Big] (\o)
  \=  \ul{\sE}_{t_i}   \Big[ \U^n_{  t_{i+1}}  \Big] (\o)
  \=     \ul{\sE}_{t_i}   \bigg[ V^n_{t_{i+1}}
    \n + \n \int_{t_i}^{t_{i+1}} g_r dr  \bigg] (\o) + \int_0^{t_i} g_r (\o) dr
    \ge  \U^n_{  t_i}   (\o)
    \= \U^n_{\nu^{n,\d}_k \land t_i}   (\o)   , 
    \eeas
    which together with \eqref{eq:db017} proves \eqref{eq:db014} for $   i = i_k + 1 , \cds , 2^k - 1 $.

 \no {\bf 3c)} {\it  As a consequence of \eqref{eq:db014}, one then has  }
  \bea   \label{eq:db019}
   \ul{\sE}_t \Big[ \U^n_{\nu^{n,\d}_k \land t_i} \Big] (\o)
    \le  \ul{\sE}_t \Big[ \U^n_{\nu^{n,\d}_k \land t_{i+1}} \Big] (\o) , \q
     i = i_k + 1    , \cds , 2^k - 1 .
   \eea

   Let   $i \= i_k  \+ 1   , \cds \n , 2^k \- 1  $ and $\hP  \n \in \n  \cP(t,\o)$.
   As 
   $\xi_i \df \U^n_{\nu^{n,\d}_k \land t^k_{i+1}} $  is $   \cF_T  -$measurable by Remark \ref{rem_V_adapted},
   Proposition   \ref{prop_shift0} (1) shows that $\eta_i : = \xi^{t,\o}_i $ is $ \cF^t_T  -$measurable.
   Since  \eqref{eq:uxu170} and the first equality in \eqref{eq:xax141}  show  that for any $\wt{\o} \n \in \n \O^t$
   \bea
   |\eta_i (\wt{\o})| 
    \n \le \n   \Psi \big( \nu^{n,\d}_k (\o  \n \otimes_t \n  \wt{\o} )  \n \land \n  t^k_{i+1}  ,
    \o  \n \otimes_t \n  \wt{\o} \big)
    \n + \dn  \int_0^T \n  \big| g_r (\o  \n \otimes_t \n  \wt{\o}) \big| dr
   \ls  \underset{r \in [t,T]}{\sup}  \Psi_r (\o \oti_t \wt{\o} ) \+
   \int_0^t \n | g_r (\o  )  | dr \+ \int_t^T \n \big| g^{t,\o}_r (  \wt{\o}) \big| dr , \q
   \label{eq:db131}
   \eea
       an analogy to \eqref{eq:uxu183} and     \eqref{eq:f475} imply   that
      for all $\wt{\o}  \n \in \n  \O^t $ except on a $\hP-$null set $\cN_i$,
      \bea   \label{eq:db211}
       \hE_{\hP^{{t_i}, \wt{\o}}} \Big[ \eta^{{t_i}, \wt{\o}}_i  \Big]
         =    \hE_\hP \big[ \eta_i \big| \cF^t_{t_i} \big] (\wt{\o})   \in   \hR     .
       \eea
   By (P2),    there exists an   extension $ \big(\O^t,\cF',\hP' \big)$ of $(\O^t,\cF^t_T,\hP)$ and
   $\O'  \n \in \n  \cF'$ with $\hP' \big(\O'\big)  \n = \n  1$ such that
    $\hP^{{t_i},   \wt{\o}}  \n \in \n    \cP (s ,  \o  \otimes_t \wt{\o} ) $
    for any $\wt{\o}  \n \in \n  \O'$. Given $ \wt{\o} \in \O' \cap  \cN^c_i  $, since
 \beas
  \eta^{{t_i}, \wt{\o}}_i (\wh{\o}) = \eta_i (\wt{\o} \otimes_{t_i} \wh{\o})
  = \xi^{t,\o}_i  (\wt{\o} \otimes_{t_i} \wh{\o})
  = \xi_i \big( \o \otimes_t (\wt{\o} \otimes_{t_i} \wh{\o}) \big)
  = \xi_i \big( (\o \otimes_t  \wt{\o} ) \otimes_{t_i} \wh{\o}  \big)
  = \xi^{{t_i}, \o \otimes_t  \wt{\o} }_i (\wh{\o}), \q \fa \wh{\o} \in \O^s ,
 \eeas
 we can deduce from \eqref{eq:db014} and \eqref{eq:db211} that
 \beas
  \big( \U^n_{\nu^{n,\d}_k \land t_i} \big)^{t,\o}  (\wt{\o}) & \tn =& \tn
  \big( \U^n_{\nu^{n,\d}_k \land t_i} \big) (\o \otimes_t \wt{\o})
  \le  \ul{\sE}_{t_i}
   \Big[ \U^n_{\nu^{n,\d}_k \land t_{i+1}}  \Big] (\o \otimes_t \wt{\o})
   = \ul{\sE}_{t_i}  [ \xi_i   ] (\o \otimes_t \wt{\o})
=      \underset{\wt{\hP} \in \cP({t_i}, \o \otimes_t \wt{\o} )}{\inf}
\hE_{\wt{\hP} } \Big[\xi^{{t_i}, \o \otimes_t \wt{\o}}_i \Big] \\
 & \tn \le & \tn  \hE_{\hP^{{t_i}, \wt{\o}}} \Big[\xi^{{t_i}, \o \otimes_t \wt{\o}}_i \Big]
= \hE_{\hP^{{t_i}, \wt{\o}}} \Big[\eta^{{t_i},\wt{\o}}_i \Big]
     =    \hE_\hP \big[\eta_i |\cF^t_{t_i}\big] (\wt{\o})
= \hE_\hP \Big[\xi^{t,\o}_i  \Big|\cF^t_{t_i}\Big]  (\wt{\o})   ,
\eeas
 which shows that   $\O'  \n \cap \n   \cN^c_i  \n \subset \n  \cA'  \n := \n
  \Big\{  \big( \U^n_{\nu^{n,\d}_k \land t_i} \big)^{t,\o}
   \ls \hE_\hP \big[\xi^{t,\o}_i  \big|\cF^t_{t_i} \big]  \Big\} \ins \cF^t_T    $.
   It   follows that
   $   \hP (\cA') \n = \n  \hP'(\cA')  \n \ge \n  \hP' \big(\O' \cap   \cN^c_i \big)  \n = \n 1 $.
   Hence,  $ \big( \U^n_{\nu^{n,\d}_k \land t_i} \big)^{t,\o}
    \ls \hE_\hP \Big[\xi^{t,\o}_i  \Big|\cF^t_{t_i} \Big] $,   \pas ~
   Taking the expectation $\hE_\hP [\cd]$  yields that
 $  \hE_\hP \Big[ ( \U^n_{\nu^{n,\d}_k \land t_i} )^{t,\o} \Big] \ls
\hE_\hP \big[ \xi^{t,\o}_i  \big]
=  \hE_\hP \Big[ \big( \U^n_{\nu^{n,\d}_k \land t_{i+1}} \big)^{t,\o} \Big]  $.
 Then  taking infimum over $\hP \in \cP(t, \o )$, we obtain \eqref{eq:db019}.
 \if{0}
\beas
    \ul{\sE}_t \big[  \U^n_{\nu^{n,\d}_k \land t_i}  \big] (\o)
 = \underset{\hP \in \cP(t, \o )}{\inf} \hE_\hP \Big[ ( \U^n_{\nu^{n,\d}_k \land t_i} )^{t,\o} \Big]
 \le \underset{\hP \in \cP(t, \o )}{\inf} \hE_\hP \big[ ( \U^n_{\nu^{n,\d}_k \land t_{i+1}} )^{t,\o} \big]
 =  \ul{\sE}_t [\U^n_{\nu^{n,\d}_k \land t_{i+1}}] (\o) .
\eeas
 \fi

  \no {\bf 3d)} {\it Finally, we will use \eqref{eq:db019} as well as the continuity of process $\ol{V}$
   to reach \eqref{eq:dc251} for the case $t \< T$ with $\ol{V}_t(\o)  \n > \n  L_t(\o) $. }

  Taking $i \= i_k$ in \eqref{eq:db014} shows that
  $    \U^n_t   (\o) \n =  \n   \U^n_{\nu^{n,\d}_k \land t}   (\o)
  \n \le \n
       \ul{\sE}_t
   \Big[ \U^n_{\nu^{n,\d}_k \land t_{i_k+1}}  \Big] (\o)    $,
    which together with   \eqref{eq:db019} and \eqref{eq:da371} yields that
 \bea
 \U^n_t   (\o) \n \le  \n    \ul{\sE}_t  \Big[ \U^n_{\nu^{n,\d}_k \land t_{i_k+1}}  \Big] (\o)
    \n \le  \n    \ul{\sE}_t   \Big[ \U^n_{\nu^{n,\d}_k \land t_{i_k+2}}  \Big] (\o)
    \n \le \n  \cds  \n \le  \n       \ul{\sE}_t   \Big[ \U^n_{\nu^{n,\d}_k \land t_{2^k}}  \Big] (\o)
    \n = \n   \ul{\sE}_t   \big[ \U^n_{\nu^{n,\d}_k  }  \big] (\o)
    \n \le \n   \ul{\sE}_t   \big[ \ol{\U}_{\nu^{n,\d}_k  }  \big] (\o)  . \q   \label{eq:dc021}
 \eea

 Since $\lmtd{k \to \infty} \nu^{n,\d}_k  \= \nu^{n,\d}  $,
 the continuity of   $\ol{V}$ by Proposition \ref{prop_conti_V}   implies that
 $ \lmt{k \to \infty}   \ol{\U}_{\nu^{n,\d}_k  }   \=   \ol{\U}_{ \nu^{n,\d}  } $.
    Also,  an analogy to \eqref{eq:db131} that for any $\wt{\o} \n \in \n \O^t$
   \bea
   \Big| \big( \ol{\U}_{\nu^{n,\d}_k}  \big)^{t,\o} (\wt{\o}) \Big|
   \ls \Psi^{t,\o}_* (\wt{\o}) \+
   \int_0^t \n | g_r (\o  )  | dr \+ \int_t^T \n \big| g^{t,\o}_r (  \wt{\o}) \big| dr . \qq
   \label{eq:dc023}
   \eea
 Then for any $\hP \n  \in  \n  \cP(t,\o)$,  the  dominated convergence theorem and
  an analogy to \eqref{eq:uxu183} imply that
 $ \lmt{k \to \infty} \hE_\hP \Big[ \big( \ol{\U}_{\nu^{n,\d}_k  } \big)^{t, \o} \Big]
 \= \hE_\hP \Big[ \big( \ol{\U}_{ \nu^{n,\d}  } \big)^{t, \o} \Big] $. Taking infimum over $\hP \n  \in  \n  \cP(t,\o)$
 and letting $ k \to \infty $ in \eqref{eq:dc021},  we obtain
 \beas
  \U^n_t   (\o) \ls 
  \lsup{k \to \infty} \, \underset{\hP \in \cP(t , \o)}{\inf}
    \hE_\hP    \Big[ \big( \ol{\U}_{\nu^{n,\d}_k  } \big)^{t, \o} \Big]
    \ls  \underset{\hP \in \cP(t , \o)}{\inf} \, \lsup{k \to \infty}
    \hE_\hP    \Big[ \big( \ol{\U}_{\nu^{n,\d}_k  } \big)^{t, \o} \Big]
    \=  \underset{\hP \in \cP(t , \o)}{\inf} \,  \hE_\hP \Big[ \big( \ol{\U}_{ \nu^{n,\d}  } \big)^{t, \o} \Big] .
 \eeas
 As $\|\o\|_{0,t} \le \|\o\|_{0, T} < \a $, we further see from \eqref{eq:da365}    that
 \bea
    \ol{\U}_t (\o) 
 \n \le \n  \U^n_t (\o)  \n + \n  \rho_\a (2^{-n})     \n + \n  2^{-n} \big( |g_t(\o)|  \n + \n  \rho_\a (T \n - \n t) \big)
  \n \le  \n  \underset{\hP \in \cP(t , \o)}{\inf} \,  \hE_\hP \Big[ \big( \ol{\U}_{ \nu^{n,\d}  } \big)^{t, \o} \Big]
   \n + \n  \rho_\a (2^{-n})    \n  + \n  2^{-n} \big( |g_t(\o)|  \n + \n  \rho_\a (T \n - \n t) \big) .
   \q \; \;    \label{eq:dc027}
 \eea

 The path regularity of $V^n$ in Proposition \ref{prop_conti_V} implies   that
  \bhe  \bea \label{eq:uxu201}
   \lmtu{\d \to 0}  \lmtu{n \to \infty} \tau^{n,\d}_{(t,\o)} (\wt{\o})  \n = \n    \tau^*_{(t,\o)} (\wt{\o}) ,
   \q   \fa   \wt{\o} \ins \O^t .
   \eea \ehe
 The continuity of   $\ol{V}$ thus   shows that
  $   \lmt{\d \to 0} \, \lmt{n \to \infty}     \ol{\U}_{ \nu^{n,\d}  }
 =  \lmt{\d \to 0} \, \lmt{n \to \infty}  \ol{\U}_{   \big( \tau^{n,\d}_{(t',\o)} (\Pi^0_{t'})  \land  \z  \big) \vee t }
 =  \ol{\U}_{ \big( \tau^*_{(t',\o)} (\Pi^0_{t'})     \land    \z  \big) \vee t  }  $.
 Also, letting $k \to \infty$ in \eqref{eq:dc023} yields that
 $  \big|  ( \ol{\U}_{\nu^{n,\d}}   )^{t,\o}  \big|
 \ls  \Psi^{t,\o}_*   \+  \int_0^t  | g_r (\o  )  | dr \+ \int_t^T \n \big| g^{t,\o}_r   \big| dr $.
 Then for any $\hP \n  \in  \n  \cP(t,\o)$,    applying the   dominated convergence theorem
 and  an analogy to \eqref{eq:uxu183} again, we obtain    that
 $  \lmt{\d \to 0} \,  \lmt{n \to \infty}  \hE_\hP \Big[ \big( \ol{\U}_{ \nu^{n,\d}  } \big)^{t, \o} \Big]
 =  \hE_\hP \bigg[ \Big( \,
 \ol{\U}_{   \big( \tau^*_{(t',\o)} (\Pi^0_{t'})     \land    \z  \big) \vee t  }  \Big)^{t, \o} \bigg]  $.
 Eventually,  letting $n \to \infty$ and $\d \to 0$ in \eqref{eq:dc027}  yields that
 \beas
 \q   \ol{\U}_t (\o) & \tn \le & \tn   \lsup{\d \to 0} \, \lsup{n \to \infty} \,
    \underset{\hP \in \cP(t , \o)}{\inf} \,  \hE_\hP \Big[ \big( \ol{\U}_{ \nu^{n,\d}  } \big)^{t, \o} \Big]
    \ls   \lsup{\d \to 0} \,   \underset{\hP \in \cP(t , \o)}{\inf} \,
    \lsup{n \to \infty} \,   \hE_\hP \Big[ \big( \ol{\U}_{ \nu^{n,\d}  } \big)^{t, \o} \Big]
    \ls    \underset{\hP \in \cP(t , \o)}{\inf} \,  \lsup{\d \to 0} \,
    \lsup{n \to \infty} \,   \hE_\hP \Big[ \big( \ol{\U}_{ \nu^{n,\d}  } \big)^{t, \o} \Big] \\
     & \tn = & \tn  \underset{\hP \in \cP(t , \o)}{\inf} \,
   \hE_\hP \bigg[ \Big( \,
 \ol{\U}_{   \big( \tau^*_{(t',\o)} (\Pi^0_{t'})     \land    \z  \big) \vee t  }  \Big)^{t, \o} \bigg]  .
 \hspace{10cm} \hb{\qed}
 \eeas

 \subsection{Proofs of Starred Statements in Section \ref{sec:proofs}}

 \no {\bf Proof of \eqref{eq:da041}:}
   When  $n = \infty$,  applying \eqref{eq:xxx617c}
   with $A = \{\tau  \n \land \n  \ga   \n \ge \n   s \}  \n \in \n  \cF^t_s $
   and $\tau \= \tau  \n \vee \n  s \ins  \cT^t_s   $ shows that
     \beas
     \q   \sum^\l_{j=1}  \hE_{\wh{\hP}_\l}  \big[ \b1_{\{\t  \land \ga  \ge  s\} \cap \cA_j}
         R^{t,\o}  (\tau   , \wp^n_j )   \big]
         \=  \sum^\l_{j=1}   \hE_{\wh{\hP}_\l}  \big[ \b1_{\{\tau  \land \ga  \ge  s\} \cap \cA_j}
         R^{t,\o}  (\tau  \n \vee \n  s  , \wp^n_j )   \big]
     \ls  \hE_{\hP}  \bigg[ \b1_{  \{\tau  \land \ga  \ge  s\} \cap \cA_0^c }
     \Big(  V^{t,\o}_s     \n + \n  \int_t^s \n g^{t,\o}_r   dr  \Big)          \bigg] \+   \e .
     \eeas
     On the other hand,  if  $n < \infty$, let $ i_s $ be the smallest integer such that
      $ i_s 2^{-n} T \ge s    $. Clearly, $ \tau \n \vee \n  (i_s 2^{-n} T) \ins \cT^t_s  (n) $.
      Since $ \{\tau  \land \ga  \ge  s\}
      \subset \{\tau     \ge  s\} = \{\tau \ge i_s 2^{-n} T \} $,  applying \eqref{eq:xxx617c} again
      with $ A = \{\tau  \n \land \n  \ga   \n \ge \n   s \} $
      and $\tau \= \tau  \n \vee \n  (i_s 2^{-n} T)    $ yields that
          \beas
          \hspace{-5mm}
   \sum^\l_{j=1}  \hE_{\wh{\hP}_\l}  \big[ \b1_{\{\t  \land \ga  \ge  s\} \cap \cA_j}
         R^{t,\o}  (\tau   , \wp^n_j )   \big]
         \=  \sum^\l_{j=1}   \hE_{\wh{\hP}_\l}  \big[ \b1_{\{\tau  \land \ga  \ge  s\} \cap \cA_j}
         R^{t,\o}  (\tau \n \vee \n (i_s 2^{-n} T)  , \wp^n_j )   \big]
     \ls  \hE_{\hP}  \bigg[ \b1_{  \{\tau  \land \ga  \ge  s\} \cap \cA_0^c }
     \Big( ( V^n_s )^{t,\o}     \n + \dn  \int_t^s \n g^{t,\o}_r   dr  \Big)          \bigg] \+   \e .
     \eeas

 \no {\bf Proof of \eqref{eq:uxu151}:}
 We set $ A^s_0 \n : = \n  \{  \ga   \n   < \n  s\}
   \n  \cup  \n   \big( \{  \ga  \n \ge \n  s\}  \n \cap \n  \cA_0 \big)  \n \in \n  \cF^t_s $
  and $A^s_j  \n := \n  \{  \ga  \n \ge \n  s\}  \n \cap \n  \cA_j   \n \in \n  \cF^t_s   $.
  Given $r  \n \in \n  [t,T]$,
  \bea  \label{eq:da031}
  \{\wh{\ga}_\l \le r\} = \big( A^s_0 \cap \{\ga \le r\} \big) \cup
   \Big(   \underset{j=1}{\overset{\l}{\cup}} \n \big( A^s_j \cap \{ \wp^n_j \le r \} \big) \Big) .
   \eea
 If $r  \n < \n   s$, since  $ \{\ga   \n \le \n  r\}
   \n \subset \n  \{   \ga  \n < \n  s \}    $
  and since each $\wp^n_j  \n \in \n  \cT^t_s$,   one has
  $  \{\wh{\ga}_\l  \n \le \n  r\}    \n  = \n \{   \ga  \n < \n  s \} \n  \cap \n \{\ga  \n \le \n  r\}
   \n = \n   \{\ga  \n \le \n  r\}  \n \in \n  \cF^t_r $.
 Otherwise,   if $ r  \n \ge \n   s  $,
 as $ A^s_j  \n \in \n  \cF^t_s  \n \subset \n  \cF^t_r $
 for $j \n = \n 0,1,\cds \n , \l$,   \eqref{eq:da031} also implies that
   $\{\wh{\ga}_\l  \n \le \n  r\}  \n \in \n  \cF^t_r$.
  Hence, $\wh{\ga}_\l  \n \in \n  \cT^t $. \qed

 \no {\bf Proof of \eqref{eq:uxu181}:}
   Since $\z_{\wt{\o}}  \n = \n  \lmtd{k \to \infty} \z^k_{\wt{\o}}$, we see that
      $ \{\z_{\wt{\o}}  \n < \n  \ga_{\wt{\o}}  \}  \n \subset \n
      A_{\wt{\o}} \n := \n \underset{k \in \hN}{\cup} \big\{  \z^k_{\wt{\o}}  \n \le \n  \ga_{\wt{\o}} \big\}
       \n \subset \n  \{\z_{\wt{\o}}  \n \le \n  \ga_{\wt{\o}}  \} $ and thus
       $\{\z_{\wt{\o}}  \n \le \n  \ga_{\wt{\o}}  \} \backslash A_{\wt{\o}}  \n \subset \n  \{\z_{\wt{\o}}  \n = \n  \ga_{\wt{\o}}  \} $.
       Then the continuity of process $L$ implies that
\beas
&& \hspace{-1.2cm} R^{s,\o  \otimes_t  \wt{\o}} (\z_{\wt{\o}},\ga_{\wt{\o}} )
\n  = \n   \int_s^{\z_{\wt{\o}} \land \ga_{\wt{\o}}} g^{s,\o  \otimes_t  \wt{\o}}_r dr
+ \b1_{\{\z_{\wt{\o}} \le \ga_{\wt{\o}}\}} L^{s,\o  \otimes_t  \wt{\o}}_{\z_{\wt{\o}}}
+ \b1_{\{\ga_{\wt{\o}} < \z_{\wt{\o}}   \}} U^{s,\o  \otimes_t  \wt{\o}}_{\ga_{\wt{\o}}}  \nonumber \\
& &  \hspace{-0.7cm}  =    \int_s^{\z_{\wt{\o}} \land \ga_{\wt{\o}}} \n g^{s,\o  \otimes_t  \wt{\o}}_r dr
\+ \b1_{A_{\wt{\o}}} L^{s,\o  \otimes_t  \wt{\o}}_{\z_{\wt{\o}}}
\+ \b1_{\{\z_{\wt{\o}}  \n \le \n  \ga_{\wt{\o}}  \} \backslash A_{\wt{\o}}} L^{s,\o  \otimes_t  \wt{\o}}_{\ga_{\wt{\o}}}
\+ \b1_{\{\ga_{\wt{\o}} < \z_{\wt{\o}}   \}} U^{s,\o  \otimes_t  \wt{\o}}_{\ga_{\wt{\o}}}
   \ls    \int_s^{\z_{\wt{\o}} \land \ga_{\wt{\o}}} \n g^{s,\o  \otimes_t  \wt{\o}}_r dr
\+ \b1_{A_{\wt{\o}}} L^{s,\o  \otimes_t  \wt{\o}}_{\z_{\wt{\o}}}
\+ \b1_{A^c_{\wt{\o}}} U^{s,\o  \otimes_t  \wt{\o}}_{\ga_{\wt{\o}}}    \nonumber  \\
 & &  \hspace{-0.7cm}  =
   \lmt{k \to \infty} \bigg( \int_s^{\z^k_{\wt{\o}} \land \ga_{\wt{\o}}} \n  g^{s,\o  \otimes_t  \wt{\o}}_r dr
 \n + \n  \b1_{\{\z^k_{\wt{\o}} \le \ga_{\wt{\o}}\}} L^{s,\o  \otimes_t  \wt{\o}}_{\z^k_{\wt{\o}}}
 \n + \n  \b1_{\{\ga_{\wt{\o}} < \z^k_{\wt{\o}}   \}} U^{s,\o  \otimes_t  \wt{\o}}_{\ga_{\wt{\o}}} \bigg)
 \n = \n   \lmt{k \to \infty} R^{s,\o  \otimes_t  \wt{\o}} (\z^k_{\wt{\o}},\ga_{\wt{\o}} ) .
\eeas

   \no {\bf Proof of \eqref{eq:uxu175}:}
 For any $\tau_1, \tau_2  \n \in \n  \cT^t_s $, letting
$A  \n : = \n  \{ \hE_\hP \big[R^{t,\o} (\tau_1,\wh{\ga}')\big|\cF^t_s\big]
 \n \ge  \n  \hE_\hP \big[R^{t,\o} (\tau_2,\wh{\ga}')\big|\cF^t_s\big]  \}  \n \in \n  \cF^t_s $
and   $ \ol{\tau}  \n : = \n  \b1_A \tau_1  \n + \n  \b1_{A^c} \tau_2  \n \in \n  \cT^t_s $, we can deduce that
\beas
\hE_\hP \big[R^{t,\o} (\ol{\tau},\wh{\ga}')\big|\cF^t_s\big]
 & \tn = &  \tn  \hE_\hP \big[ \b1_A R^{t,\o} ( \tau_1 ,\wh{\ga}')
 + \b1_{A^c} R^{t,\o} ( \tau_2 ,\wh{\ga}') \big|\cF^t_s\big]
 =  \b1_A \hE_\hP \big[ R^{t,\o} ( \tau_1 ,\wh{\ga}') \big|\cF^t_s \big]
 + \b1_{A^c} \hE_\hP \big[  R^{t,\o} ( \tau_2 ,\wh{\ga}') \big|\cF^t_s \big] \\
& \tn = & \tn  \hE_\hP \big[R^{t,\o} (\tau_1,\wh{\ga}')\big|\cF^t_s\big]
 \vee \hE_\hP \big[R^{t,\o} (\tau_2,\wh{\ga}')\big|\cF^t_s\big] .
\eeas
 So the family
 $ \Big\{  \hE_\hP \big[R^{t,\o} (\t,\wh{\ga}')\big|\cF^t_s\big] \Big\}_{\tau \in \cT^t_s}  $
 is directed upwards. Appealing to the basic properties of the
essential infimum (e.g., \cite[Proposition VI-\b1-1]{Neveu_1975}),
  we can find a sequence  $\{\tau_n\}_{n \in \hN} $ in $ \cT^t_s $ such that  \eqref{eq:uxu175} holds. \qed

\no {\bf Proof of \eqref{eq:uxu177}:}   For any $r \ins [t,s)$, since $\tau_n  \ins \cT^t_s$
  and since $ \{\tau  \n \le \n  r \}  \n \subset \n  \{\tau  \n < \n  s \}
  \n \subset \n  \{\tau \land \wh{\ga}  \n < \n  s \} $, one can deduce that
  $\{ \ol{\t}_n   \n \le \n  r \}  \n = \n  \{\t \land \wh{\ga}  \n < \n  s \} \cap \{\tau  \n \le \n  r \}
   \n = \n  \{ \tau  \n \le \n  r \}  \n \in \n  \cF^t_r   $.
  On the other hand,  for any $r  \n \in \n  [s,T]$, $\{ \ol{\t}_n   \n \le \n  r \}
     \n = \n  \big( \{\t \land \wh{\ga}  \n < \n s \} \cap \{\tau  \n \le \n  r \} \big) \cup
  \big(  \{\t \land \wh{\ga}  \n \ge \n  s \} \cap \{\tau_n   \n \le \n  r\} \big)    \n \in \n  \cF^t_r $.
  Hence, $ \ol{\t}_n \ins \cT^t $.

 \no {\bf Proof of \eqref{eq:uxu214}:}
  Given $r  \n \in \n  [s,T]$,  as $A_r  \n := \n  \{\wp  \n \le \n  r \}  \n \in \n  \cF^s_r $,
   \eqref{eq:xxx439} shows that
\bea \label{eq:xax119}
 \{\wp_j   \le r\} = \big\{\wh{\o} \in \O^s : \wp \big( \cX^j (\wh{\o}) \big) \le r \big\}
  =  \{\wh{\o} \in \O^s : \cX^j (\wh{\o}) \in A_r \}
  = (\cX^j)^{-1} (A_r) \in \ol{\cF}^s_r   .
\eea
 Also,  Lemma A.3 in the ArXiv version of \cite{ROSVU} implies that
 $ \big\{  \nu_j   \n \le \n  r \big\}  \n  = \n  \big\{\wt{\o}  \n \in \n  \O^t  \n :
  \Pi^t_s  (\wt{\o})  \n  \in \n  \{ \wp_j   \n \le \n  r \} \big\}
   \n = \n   (\Pi^t_s)^{-1} \big(  \{\wp_j    \n \le \n  r\}  \big) \n \in \n \ol{\cF}^t_r $, then
 one can deduce from \eqref{eq:xax047} that
 $  \{\wh{\ga}_j   \n  \le \n  r\}  \n = \n  \big\{\wt{\o}  \n \in \n  \O^t  \n :
    \wh{\cW} (\wt{\o})   \n \in \n  \{ \nu_j  \n  \le \n  r \} \big\}
   \n = \n   \wh{\cW}^{-1} \big(   \{ \nu_j  \n  \le \n  r \}   \big)
   \n \in \n  \wh{\fF}_r $.
 Hence,  $   \wp_j \n \in  \n  \ol{\cT}^s  $, $\nu_j   \n \in \n  \ol{\cT}^t_s $ while
 $\wh{\ga}_j$ is a $\wh{\fF}-$stopping time  that takes values in $[s,T]$. \qed

\no {\bf Proof of \eqref{eq:xax149}:}
   When $n \n < \n  \infty$, as induced by $\tau  \n \in \n  \cT^t_s (n)$, $\z_{\wt{\o}}$   takes   values in
$\{t^n_i\}^{2^n}_{i=i_s}$, where $i_s$ be the smallest integer such that $i_s 2^{-n} T  \n \ge \n  s$.
Similar to \eqref{eq:xax143}, there exists   $\z'_{\wt{\o}}  \n \in \n  \cT^s (n)$ such that
$\z'_{\wt{\o}}  \n = \n  \z_{\wt{\o}} $, $\fp_j-$a.s. So we have
   \beas
     \hE_{\fp^j} \Big[ R^{s,\o  \otimes_t  \cX  (\wt{\o}  )}
    (\z_{\wt{\o}},\wp )  \Big]
  \n  =   \n  \hE_{\fp^j} \Big[ R^{s,\o  \otimes_t  \cX  (\wt{\o}  )}
    (\z'_{\wt{\o}},\wp )  \Big]
     \n  =   \n  \hE_{\hP^j} \Big[ R^{s,\o  \otimes_t  \cX  (\wt{\o}  )}    (\z'_{\wt{\o}},\wp )  \Big]
   \n  \le  \n    \underset{\vs \in \cT^s (n)}{\sup} \hE_{\hP^j}
   \Big[ R^{s,\o  \otimes_t  \cX  (\wt{\o}  )}  (\vs ,\wp )  \Big] . \q
   \eeas

     Suppose  $n  \n = \n  \infty$ now. Let $k  \n \in \n  \hN$ and
    set    $s^k_i    \n  :=  \n   s    \vee    (i 2^{-k} T)   $, $i  \n   =   \n   0 , \cds   \n   , 2^k  $.
    With  $s^k_{-1}  \n := \n  -1 $,
     $  \z^k_{\wt{\o}}  \n := \n  \sum^{2^k}_{i=0}  \n  \b1_{\{ s^k_{i-1}  <  \z_{\wt{\o}}
    \le    s^k_i  \}}  s^k_i$ defines a $ \fF^j-$stopping time.
    By similar arguments to  those that lead to \eqref{eq:xax143}, one can construct a $\cT^s-$stopping time $\vs^k_{\wt{\o}}$
     valued in $\{s^k_i\}^{2^k}_{i=0}$ such that $\z^k_{\wt{\o}}  \n = \n  \vs^k_{\wt{\o}}$, $\fp^j-$a.s.
     Since $\z_{\wt{\o}}  \n = \n  \lmtd{k \to \infty} \z^k_{\wt{\o}}$,
     an analogy to \eqref{eq:uxu181} shows that
     \bea
  R^{s,\o  \otimes_t  \cX  (\wt{\o}  )} (\z_{\wt{\o}},\wp )
\ls  \lmt{k \to \infty} R^{s,\o  \otimes_t  \cX  (\wt{\o}  )} (\z^k_{\wt{\o}},\wp ) .
 \label{eq:xax147}
 \eea
  \if{0}
     we see that
      $ \{\z_{\wt{\o}}  \n < \n  \wp  \}  \n \subset \n
      A_{\wt{\o}} \n := \n \underset{k \in \hN}{\cup} \big\{  \z^k_{\wt{\o}}  \n \le \n  \wp \big\}
       \n \subset \n  \{\z_{\wt{\o}}  \n \le \n  \wp  \} $ and thus
       $\{\z_{\wt{\o}}  \n \le \n  \wp  \} \backslash A_{\wt{\o}}  \n \subset \n  \{\z_{\wt{\o}}  \n = \n  \wp  \} $.
       Then the continuity of process $L$ implies that
\bea
&& \hspace{-1cm} R^{s,\o  \otimes_t  \cX  (\wt{\o}  )} (\z_{\wt{\o}},\wp )
\n  = \n   \int_s^{\z_{\wt{\o}} \land \wp} g^{s,\o  \otimes_t  \cX  (\wt{\o}  )}_r dr
+ \b1_{\{\z_{\wt{\o}} \le \wp\}} L^{s,\o  \otimes_t  \cX  (\wt{\o}  )}_{\z_{\wt{\o}}}
+ \b1_{\{\wp < \z_{\wt{\o}}   \}} U^{s,\o  \otimes_t  \cX  (\wt{\o}  )}_\wp \nonumber \\
& &  =    \int_s^{\z_{\wt{\o}} \land \wp} g^{s,\o  \otimes_t  \cX  (\wt{\o}  )}_r dr
+ \b1_{A_{\wt{\o}}} L^{s,\o  \otimes_t  \cX  (\wt{\o}  )}_{\z_{\wt{\o}}}
+ \b1_{\{\z_{\wt{\o}}  \n \le \n  \wp  \} \backslash A_{\wt{\o}}} L^{s,\o  \otimes_t  \cX  (\wt{\o}  )}_\wp
+ \b1_{\{\wp < \z_{\wt{\o}}   \}} U^{s,\o  \otimes_t  \cX  (\wt{\o}  )}_\wp \nonumber  \\
& &  \le    \int_s^{\z_{\wt{\o}} \land \wp} g^{s,\o  \otimes_t  \cX  (\wt{\o}  )}_r dr
+ \b1_{A_{\wt{\o}}} L^{s,\o  \otimes_t  \cX  (\wt{\o}  )}_{\z_{\wt{\o}}}
+ \b1_{A^c_{\wt{\o}}} U^{s,\o  \otimes_t  \cX  (\wt{\o}  )}_\wp \nonumber  \\
 & &  =
   \lmt{k \to \infty} \bigg( \int_s^{\z^k_{\wt{\o}} \land \wp} \n  g^{s,\o  \otimes_t  \cX  (\wt{\o}  )}_r dr
 \n + \n  \b1_{\{\z^k_{\wt{\o}} \le \wp\}} L^{s,\o  \otimes_t  \cX  (\wt{\o}  )}_{\z^k_{\wt{\o}}}
 \n + \n  \b1_{\{\wp < \z^k_{\wt{\o}}   \}} U^{s,\o  \otimes_t  \cX  (\wt{\o}  )}_\wp \bigg)
 \n = \n   \lmt{k \to \infty} R^{s,\o  \otimes_t  \cX  (\wt{\o}  )} (\z^k_{\wt{\o}},\wp ) .   \qq
 \label{eq:xax147}
\eea
\fi
  By \eqref{eq:ab015},  $ \big| R^{s,\o  \otimes_t  \cX  (\wt{\o}  )} (\z^k_{\wt{\o}},\wp ) \big|
 \le \int_s^T \Big| g^{s,\o  \otimes_t  \cX  (\wt{\o}  )}_r \Big| dr + \Psi^{s,\o  \otimes_t  \cX  (\wt{\o}  )}_* $,
 $\fa k \in \hN$. Since $\hP_j \ins \wh{\fP}_s$ by \eqref{eq:uxu217}, \eqref{eq:xxx111}   shows that
 \beas
  \hE_{\fp^j} \bigg[ \int_s^T  \Big| g^{s,\o  \otimes_t  \cX  (\wt{\o}  )}_r \Big|  dr
 +   \Psi^{s,\o  \otimes_t  \cX  (\wt{\o}  )}_* \bigg] =
  \hE_{\hP^j} \bigg[ \int_s^T  \Big| g^{s,\o  \otimes_t  \cX  (\wt{\o}  )}_r \Big|  dr
 +   \Psi^{s,\o  \otimes_t  \cX  (\wt{\o}  )}_* \bigg] < \infty .
 \eeas
 Taking expectation $\hE_{\fp^j} [~] $ in  \eqref{eq:xax147},
 one can deduce from  the  dominated convergence theorem  that
 \beas
 \hspace{2cm}
  \hE_{\fp^j} \Big[ R^{s,\o  \otimes_t  \cX  (\wt{\o}  )} (\z_{\wt{\o}},\wp ) \Big]
 & \tn \le & \tn  \lmt{k \to \infty} \hE_{\fp^j} \Big[ R^{s,\o  \otimes_t  \cX  (\wt{\o}  )} (\z^k_{\wt{\o}},\wp ) \Big]
  = \lmt{k \to \infty} \hE_{\fp^j} \Big[ R^{s,\o  \otimes_t  \cX  (\wt{\o}  )} (\vs^k_{\wt{\o}},\wp ) \Big] \\
  & \tn  = & \tn  \lmt{k \to \infty} \hE_{\hP^j} \Big[ R^{s,\o  \otimes_t  \cX  (\wt{\o}  )} (\vs^k_{\wt{\o}},\wp ) \Big]
    \le \underset{\vs \in \cT^s}{\sup}
    \hE_{\hP^j} \Big[ R^{s,\o  \otimes_t  \cX  (\wt{\o}  )} (\vs ,\wp ) \Big] .  \hspace{2.3cm} \hb{\qed}
\eeas

      \no {\bf Proof of \eqref{eq:uxu211}:}
      If $t_i \< \tau^{n,\d}_{(t',\o)} \big(\Pi^0_{t'} (\o)\big)$,
      the definition of  $\tau^{n,\d}_{(t',\o)}$   shows that
 $ (V^n   \- L ) ( t_i , \o)  \n = \n \big( (V^n)^{t',\o} \- L^{t',\o} \big) \big( t_i , \Pi^0_{t'} (\o)\big)
   \n   \ge  \n   \d  \gs 0  $. 
    On the other hand, if  $t_i \= \tau^{n,\d}_{(t',\o)} \big(\Pi^0_{t'} (\o)\big)$
    the left-upper-semicontinuity of $ (V^n)^{t',\o}     \n - \n  L^{t',\o}  $
    implies that
    \beas
    \hspace{2.4cm}
    (V^n   \-  L ) ( t_i , \o)   \= \big( (V^n)^{t',\o}      \- L^{t',\o} \big) \big( t_i , \Pi^0_{t'} (\o)\big)
    \gs \lsup{s \nearrow t_i} \big( (V^n)^{t',\o}      \- L^{t',\o} \big) \big( s , \Pi^0_{t'} (\o)\big)
    \gs \d . \hspace{2.4cm} \hb{\qed}
    \eeas

  \no {\bf Proof of \eqref{eq:uxu201}:}
   Fix   $\wt{\o}  \n \in \n  \O^t$
   and set $\wt{\a}  \n : = \n  1  \n + \n  \|\o  \n \otimes_t \n  \wt{\o} \|_{0 , T}$.
 \if{0}
   For any $\d \n > \n 0$, similar to $\tau^*_{(t,\o)}$ in \eqref{eq:uxu163},
  $   \tau^{*,\d}_{(t,\o)} \n  : = \n
      \inf \big\{ s  \n \in \n  [t,T] \n :  \ol{V}^{t,\o}_s   \n \le  \n  L^{t,\o}_s \n + \n \d  \big\}   $
defines an  $\bF^t-$stopping time.
 \fi

 We Let $\d \> 0$, $n \ins \hN$ and simply denote $ t_{n,\d} \df \tau^{n,\d}_{(t,\o)}(\wt{\o})$,
   $ t_* \df \tau^*_{(t,\o)}(\wt{\o}) $.
  Let us first show that
 \bea   \label{eq:dc011}
   (V^n)^{t,\o} \big(t_{n,\d},\wt{\o} \big)   \le L^{t,\o} \big(t_{n,\d},\wt{\o} \big)  + \d .
 \eea
 If $ t_{n,\d} = T$, \eqref{eq:ek011} shows that
 \bea
   (V^n)^{t,\o} \big(t_{n,\d},\wt{\o} \big)
  \n = \n  (V^n)^{t,\o}  (T,\wt{\o}  )
  \n = \n  L^{t,\o} (T, \wt{\o})
    \n = \n  L^{t,\o} \big(t_{n,\d}, \wt{\o} \big)  .    \label{eq:ej014}
 \eea
 On the other hand, if $ t_{n,\d} < T$,
 let $\{\ft_i \= \ft_i (t, \o, \wt{\o},n ,\d)\}_{i \in \hN}$ be a sequence in $ \big[ t_{n,\d} , T \big]$
 such that $\lmtd{i \to \infty} \ft_i = t_{n,\d} $
 and that $ (V^n)^{t,\o} (\ft_i, \wt{\o})   \n < \n  L^{t,\o} (\ft_i, \wt{\o})  \n +  \n \d $, $\fa i \in \hN$
 by the definition of $t_{n,\d} \= \tau^{n,\d}_{(t,\o)}(\wt{\o})$.
 The right-lower-semicontinuity of path $ V^n_\cd (\o \otimes_t \wt{\o}) $ by Proposition \ref{prop_conti_V}
 and the continuity of path $L_\cd (\o \otimes_t \wt{\o}) $ then imply that
 \beas
   (V^n)^{t,\o} \big(t_{n,\d}, \wt{\o} \big)
    \=   V^n  \big(t_{n,\d}, \o \otimes_t \wt{\o} \big)
 \ls \linf{s \searrow t_{n,\d}} V^n  (s, \o \otimes_t \wt{\o}) \ls
    \linf{i \to \infty} V^n (\ft_i, \o \otimes_t  \wt{\o})
      \ls  L \big(t_{n,\d}, \o \otimes_t  \wt{\o} \big)  \n +  \n \d
     \=  L^{t,\o} \big(t_{n,\d},\wt{\o} \big)  \n +  \n \d ,
     \eeas
     which together with \eqref{eq:ej014}  proves  \eqref{eq:dc011}.

     As  $ \|\o \otimes_t \wt{\o}\|_{0 , t_{n,\d} }
    \le  \|\o \otimes_t \wt{\o} \|_{0 , T} < \wt{\a}  $,
  we see from \eqref{eq:dc011} and  \eqref{eq:da365}     that
  \bea   \hspace{-3mm}
  \ol{V}^{t,\o} \n \big(t_{n,\d},\wt{\o} \big)
  \-   L^{t,\o} \big(t_{n,\d},\wt{\o} \big)
    \ls \ol{V}  \n  \big(t_{n,\d}, \o \oti_t \wt{\o} \big)
  \-  V^n   \big(t_{n,\d}, \o \oti_t \wt{\o} \big) \+ \d
    \ls    \rho_{\wt{\a}} (2^{-n})
   \+ 2^{-n} \Big( \,  \underset{r \in [0,T]}{\sup}|g_r (\o \oti_t \wt{\o})| \+ \rho_{\wt{\a}} (T  ) \Big)  \+ \d   . \q
   \label{eq:dc014}
  \eea
  For any $s \n \in \n  [t,T]$, since $ \cT^s (n)  \n \subset \n  \cT^s (n+1)
   \n \subset \n  \cT^s $, an analogy to \eqref{eq:da371} shows that
  $  V^n_s   (\o \otimes_t \wt{\o})  \n \le \n   V^{n+1}_s (\o \otimes_t \wt{\o})
   \n \le \n   \ol{V}_s (\o \otimes_t \wt{\o})    $.
  It follows that    $\wh{t}_\d  \n : = \n  \lmtu{n \to \infty} t_{n,\d}
   \n \le \n  t_* $.
  As  $n \n \to \n  \infty$ in \eqref{eq:dc014}, the continuity of the path
  $\ol{V}^{t,\o} (\wt{\o}) \n - \n  L^{t,\o} (\wt{\o}) $ by Proposition \ref{prop_conti_V} yields  that
  $  \ol{V}^{t,\o} \big(\wh{t}_\d,\wt{\o} \big)
   \n - \n   L^{t,\o} \big(\wh{t}_\d,\wt{\o} \big)  \n \le \n  \d     $, and thus
  \bea \label{eq:ej021}
   t_* \gs  \lmtu{n \to \infty} t_{n,\d} \= \wh{t}_\d \ge t_{*,\d} \df
   \inf\big\{s  \n \in \n  [t,T] \n :  \ol{V}^{t,\o}_s (\wt{\o})  \n \le  \n  L^{t,\o}_s (\wt{\o})  \n + \n \d \big\} .
    \eea
    The continuity of the path
  $\ol{V}^{t,\o} (\wt{\o}) \n - \n  L^{t,\o} (\wt{\o}) $ also implies that
  $ t_* =   \lmtu{\d \to 0} t_{*,\d}  $
  $ $
 \if{0}
 Clearly, $ t_{*,\d} $ is decreasing in $\d$, so we set
 $\nu (\wt{\o}) : = \lmtu{\d \to 0} t_{*,\d}
  \n \le \n  t_*  $. By the continuity of  $\ol{V}^{t,\o}  \n - \n  L^{t,\o} $ again,
  $ \ol{V}^{t,\o} \big( t_{*,\d} , \wt{\o} \big)
   \n \le  \n  L^{t,\o} \big( t_{*,\d} , \wt{\o} \big)  \n + \n \d $, $\fa \d > 0$.
   Letting $\d \to 0$, we can deduce from the continuity of  $\ol{V}^{t,\o}  \n - \n  L^{t,\o} $
   and   \eqref{eq:uxu170} that
   $ \ol{V}^{t,\o} \big( \nu (\wt{\o}) , \wt{\o} \big)
   \n =  \n  L^{t,\o} \big( \nu (\wt{\o}) , \wt{\o} \big)   $.
   Therefore, $ t_* = \nu (\wt{\o}) = \lmtu{\d \to 0} t_{*,\d}  $,
 \fi
   which together with \eqref{eq:ej021} leads to that
   $  \lmtu{\d \to 0}  \lmtu{n \to \infty} t_{n,\d} =  t_*  $, i.e.,
   $  \lmtu{\d \to 0}  \lmtu{n \to \infty} \tau^{n,\d}_{(t,\o)} (\wt{\o}) =  \tau^*_{(t,\o)} (\wt{\o})  $.  \qed

{\small

\bibliographystyle{siam}
\bibliography{RDG_bib}

\begin{thebibliography}{10}

\bibitem{ALM_1982}
{\sc M.~Alario-Nazaret, J.-P. Lepeltier, and B.~Marchal}, {\em Dynkin games},
  in Stochastic differential systems ({B}ad {H}onnef, 1982), vol.~43 of Lecture
  Notes in Control and Inform. Sci., Springer, Berlin, 1982, pp.~23--32.

\bibitem{Alvarez_2008}
{\sc L.~H.~R. Alvarez}, {\em A class of solvable stopping games}, Appl. Math.
  Optim., 58 (2008), pp.~291--314.

\bibitem{OS_CRM}
{\sc E.~Bayraktar, I.~Karatzas, and S.~Yao}, {\em Optimal stopping for dynamic
  convex risk measures}, Illinois J. Math., 54 (2010), pp.~1025--1067.

\bibitem{MR3162260}
{\sc E.~Bayraktar and M.~S{\^{\i}}rbu}, {\em Stochastic {P}erron's method and
  verification without smoothness using viscosity comparison: obstacle problems
  and {D}ynkin games}, Proc. Amer. Math. Soc., 142 (2014), pp.~1399--1412.

\bibitem{ROSVU}
{\sc E.~Bayraktar and S.~Yao}, {\em On the robust optimal stopping problem},
  SIAM J. Control Optim., 52 (2014), pp.~3135--3175.

\bibitem{DRBSDEUI}
\leavevmode\vrule height 2pt depth -1.6pt width 23pt, {\em Doubly reflected
  {BSDE}s with integrable parameters and related {D}ynkin games}, Stochastic
  Process. Appl., 125 (2015), pp.~4489--4542.

\bibitem{RDOSRT}
\leavevmode\vrule height 2pt depth -1.6pt width 23pt, {\em Optimal stopping
  with random maturity under nonlinear expectations},  (2015).
\newblock Available on \url{http://arxiv.org/abs/1505.07533}.

\bibitem{Bensoussan_Friedman_1974}
{\sc A.~Bensoussan and A.~Friedman}, {\em Nonlinear variational inequalities
  and differential games with stopping times}, J. Functional Analysis, 16
  (1974), pp.~305--352.

\bibitem{Bensoussan_Friedman_1977}
\leavevmode\vrule height 2pt depth -1.6pt width 23pt, {\em Nonzero-sum
  stochastic differential games with stopping times and free boundary
  problems}, Trans. Amer. Math. Soc., 231 (1977), pp.~275--327.

\bibitem{Bismut_DG_1977}
{\sc J.-M. Bismut}, {\em Sur un probl\`eme de {D}ynkin}, Z.
  Wahrscheinlichkeitstheorie und Verw. Gebiete, 39 (1977), pp.~31--53.

\bibitem{Bismut_DG_1979}
\leavevmode\vrule height 2pt depth -1.6pt width 23pt, {\em Contr\^ole de
  processus alternants et applications}, Z. Wahrsch. Verw. Gebiete, 47 (1979),
  pp.~241--288.

\bibitem{Boetius_2005}
{\sc F.~Boetius}, {\em Bounded variation singular stochastic control and
  {D}ynkin game}, SIAM J. Control Optim., 44 (2005), pp.~1289--1321
  (electronic).

\bibitem{Buckdahn_Li_3}
{\sc R.~Buckdahn and J.~Li}, {\em Probabilistic interpretation for systems of
  {I}saacs equations with two reflecting barriers}, NoDEA Nonlinear
  Differential Equations Appl., 16 (2009), pp.~381--420.

\bibitem{Cattiaux_Lepeltier_1990}
{\sc P.~Cattiaux and J.-P. Lepeltier}, {\em Existence of a quasi-{M}arkov
  {N}ash equilibrium for nonzero sum {M}arkov stopping games}, Stochastics
  Stochastics Rep., 30 (1990), pp.~85--103.

\bibitem{Cosso_2013}
{\sc A.~Cosso}, {\em Stochastic differential games involving impulse controls
  and double-obstacle quasi-variational inequalities}, SIAM J. Control Optim.,
  51 (2013), pp.~2102--2131.

\bibitem{Cvitanic_Karatzas_1996}
{\sc J.~Cvitani{\'c} and I.~Karatzas}, {\em Backward stochastic differential
  equations with reflection and {D}ynkin games}, Ann. Probab., 24 (1996),
  pp.~2024--2056.

\bibitem{Dolinsky_2014}
{\sc Y.~Dolinsky}, {\em Hedging of game options under model uncertainty in
  discrete time}, Electron. Commun. Probab., 19 (2014), pp.~no. 19, 11.

\bibitem{Dynkin_1967}
{\sc E.~Dynkin}, {\em Game variant of a problem on optimal stopping}, Soviet
  Math. Dokl., 10 (1967), pp.~270–--274.

\bibitem{ETZ_2012}
{\sc I.~Ekren, N.~Touzi, and J.~Zhang}, {\em Optimal stopping under nonlinear
  expectation}, Stochastic Process. Appl., 124 (2014), pp.~3277--3311.

\bibitem{ETZ_2014}
\leavevmode\vrule height 2pt depth -1.6pt width 23pt, {\em Optimal stopping
  under nonlinear expectation}, Stochastic Process. Appl., 124 (2014),
  pp.~3277--3311.

\bibitem{EkrenZhang1016}
{\sc I.~Ekren and J.~Zhang}, {\em Pseudo markovian viscosity solutions of fully
  nonlinear degenerate ppdes},  (2016).
\newblock Available on \url{https://arxiv.org/abs/1604.02239}.

\bibitem{Ekstrom_2006}
{\sc E.~Ekstr{\"o}m}, {\em Properties of game options}, Math. Methods Oper.
  Res., 63 (2006), pp.~221--238.

\bibitem{Lp_DRBSDE}
{\sc B.~El~Asri, S.~Hamad{\`e}ne, and H.~Wang}, {\em {$L^p$}-solutions for
  doubly reflected backward stochastic differential equations}, Stoch. Anal.
  Appl., 29 (2011), pp.~907--932.

\bibitem{Friedman_1973}
{\sc A.~Friedman}, {\em Stochastic games and variational inequalities}, Arch.
  Rational Mech. Anal., 51 (1973), pp.~321--346.

\bibitem{Fukushima_Taksar_2002}
{\sc M.~Fukushima and M.~Taksar}, {\em Dynkin games via {D}irichlet forms and
  singular control of one-dimensional diffusions}, SIAM J. Control Optim., 41
  (2002), pp.~682--699 (electronic).

\bibitem{Hamadene_2006}
{\sc S.~Hamad{\`e}ne}, {\em Mixed zero-sum stochastic differential game and
  {A}merican game options}, SIAM J. Control Optim., 45 (2006), pp.~496--518.

\bibitem{Hamadene_Hassani_2005}
{\sc S.~Hamad\`ene and M.~Hassani}, {\em {BSDE}s with two reflecting barriers:
  the general result.}, Probab. Theory Relat. Fields, 132 (2005), pp.~237--264.

\bibitem{Hamadene_Hassani_2014a}
{\sc S.~Hamad{\`e}ne and M.~Hassani}, {\em The multi-player nonzero-sum
  {D}ynkin game in discrete time}, Math. Methods Oper. Res., 79 (2014),
  pp.~179--194.

\bibitem{Hamadene_Hdhiri_2006}
{\sc S.~Hamad{\`e}ne and I.~Hdhiri}, {\em Backward stochastic differential
  equations with two distinct reflecting barriers and quadratic growth
  generator}, J. Appl. Math. Stoch. Anal.,  (2006), pp.~Art. ID 95818, 28.

\bibitem{Hamad_Lepeltier_2000}
{\sc S.~Hamad{\`e}ne and J.-P. Lepeltier}, {\em Reflected {BSDE}s and mixed
  game problem}, Stochastic Process. Appl., 85 (2000), pp.~177--188.

\bibitem{Hamadene_Lepeltier_Wu_1999}
{\sc S.~Hamad{\`e}ne, J.-P. Lepeltier, and Z.~Wu}, {\em Infinite horizon
  reflected backward stochastic differential equations and applications in
  mixed control and game problems}, Probab. Math. Statist., 19 (1999),
  pp.~211--234.

\bibitem{Hamadene_Hassani_2014b}
{\sc S.~Hamad{\`e}ne and H.~Mohammed}, {\em The multiplayer nonzero-sum
  {D}ynkin game in continuous time}, SIAM J. Control Optim., 52 (2014),
  pp.~821--835.

\bibitem{HRZ_2009}
{\sc S.~Hamad{\`e}ne, E.~Rotenstein, and A.~Z{\u{a}}linescu}, {\em A
  generalized mixed zero-sum stochastic differential game and double barrier
  reflected {BSDE}s with quadratic growth coefficient}, An. \c Stiin\c t. Univ.
  Al. I. Cuza Ia\c si. Mat. (N.S.), 55 (2009), pp.~419--444.

\bibitem{Hamadene_Zhang_2010}
{\sc S.~Hamad{\`e}ne and J.~Zhang}, {\em The continuous time nonzero-sum
  {D}ynkin game problem and application in game options}, SIAM J. Control
  Optim., 48 (2009/10), pp.~3659--3669.

\bibitem{Kallsen_Kuhn_2004}
{\sc J.~Kallsen and C.~K{\"u}hn}, {\em Pricing derivatives of {A}merican and
  game type in incomplete markets}, Finance Stoch., 8 (2004), pp.~261--284.

\bibitem{Kara_Shr_BMSC}
{\sc I.~Karatzas and S.~E. Shreve}, {\em Brownian motion and stochastic
  calculus}, vol.~113 of Graduate Texts in Mathematics, Springer-Verlag, New
  York, second~ed., 1991.

\bibitem{Karatzas_Wang_2001}
{\sc I.~Karatzas and H.~Wang}, {\em Connections between bounded variation
  control and {D}ynkin games}, in Optimal Control and Partial Differential
  Equations (Volume in honor of A. Bensoussan), J.L. Menaldi, E. Rofman, and A.
  Sulem, eds., IOS Press, Amsterdam (2001), pp.~363--373.

\bibitem{Kara_Zam_2008}
{\sc I.~Karatzas and I.~M. Zamfirescu}, {\em Martingale approach to stochastic
  differential games of control and stopping}, Ann. Probab., 36 (2008),
  pp.~1495--1527.

\bibitem{Kifer_2000}
{\sc Y.~Kifer}, {\em Game options}, Finance Stoch., 4 (2000), pp.~443--463.

\bibitem{Kifer_2013_survey}
\leavevmode\vrule height 2pt depth -1.6pt width 23pt, {\em Dynkin games and
  {I}sraeli options}, ISRN Probability and Statistics, Volume 2013,  (2013).
\newblock Article ID 856458, 17 pages.

\bibitem{KQC_2014}
{\sc M.~Kobylanski, M.-C. Quenez, and M.~R. de~Campagnolle}, {\em Dynkin games
  in a general framework}, Stochastics, 86 (2014), pp.~304--329.

\bibitem{Laraki_Solan_2005}
{\sc R.~Laraki and E.~Solan}, {\em The value of zero-sum stopping games in
  continuous time}, SIAM J. Control Optim., 43 (2005), pp.~1913--1922
  (electronic).

\bibitem{Lepeltier_Maingueneau_1984}
{\sc J.-P. Lepeltier and M.~A. Maingueneau}, {\em Le jeu de {D}ynkin en
  th\'eorie g\'en\'erale sans l'hypoth\`ese de {M}okobodski}, Stochastics, 13
  (1984), pp.~25--44.

\bibitem{Ma_Cvitanic_2001}
{\sc J.~Ma and J.~Cvitani{\'c}}, {\em Reflected forward-backward {SDE}s and
  obstacle problems with boundary conditions}, J. Appl. Math. Stochastic Anal.,
  14 (2001), pp.~113--138.

\bibitem{MPP_2014}
{\sc A.~Matoussi, L.~Piozin, and D.~Possama{\"{\i}}}, {\em Second-order {BSDE}s
  with general reflection and game options under uncertainty}, Stochastic
  Process. Appl., 124 (2014), pp.~2281--2321.

\bibitem{Morimoto_1984}
{\sc H.~Morimoto}, {\em Dynkin games and martingale methods}, Stochastics, 13
  (1984), pp.~213--228.

\bibitem{Nagai_1987}
{\sc H.~Nagai}, {\em Non-zero-sum stopping games of symmetric {M}arkov
  processes}, Probab. Theory Related Fields, 75 (1987), pp.~487--497.

\bibitem{Neveu_1975}
{\sc J.~Neveu}, {\em Discrete-parameter martingales}, North-Holland Publishing
  Co., Amsterdam, revised~ed., 1975.
\newblock Translated from the French by T. P. Speed, North-Holland Mathematical
  Library, Vol. 10.

\bibitem{Nutz_2012a}
{\sc M.~Nutz}, {\em A quasi-sure approach to the control of non-{M}arkovian
  stochastic differential equations}, Electron. J. Probab., 17 (2012),
  pp.~1--23.

\bibitem{NZ_2012}
{\sc M.~Nutz and J.~Zhang}, {\em Optimal stopping under adverse nonlinear
  expectation and related games}, Ann. Appl. Probab., 25 (2015),
  pp.~2503--2534.

\bibitem{Ohtsubo_1987}
{\sc Y.~Ohtsubo}, {\em A nonzero-sum extension of {D}ynkin's stopping problem},
  Math. Oper. Res., 12 (1987), pp.~277--296.

\bibitem{Peng_G_2007b}
{\sc S.~Peng}, {\em {G}-{B}rownian motion and dynamic risk measure under
  volatility uncertainty}, 2007.
\newblock Available at \url{http://lanl.arxiv.org/abs/0711.2834}.

\bibitem{Rogers_Williams_2}
{\sc L.~C.~G. Rogers and D.~Williams}, {\em Diffusions, {M}arkov processes, and
  martingales. {V}ol. 2}, Cambridge Mathematical Library, Cambridge University
  Press, Cambridge, 2000.
\newblock It{\^o} calculus, Reprint of the second (1994) edition.

\bibitem{RSV_2001}
{\sc D.~Rosenberg, E.~Solan, and N.~Vieille}, {\em Stopping games with
  randomized strategies}, Probab. Theory Related Fields, 119 (2001),
  pp.~433--451.

\bibitem{Snell_1952}
{\sc J.~L. Snell}, {\em Applications of martingale system theorems}, Trans.
  Amer. Math. Soc., 73 (1952), pp.~293--312.

\bibitem{STZ_2011b}
{\sc H.~Soner, N.~Touzi, and J.~Zhang}, {\em Dual formulation of second order
  target problems}, Ann. Appl. Probab., 23 (2013), pp.~308--347.

\bibitem{Stettner_1982}
{\sc L.~Stettner}, {\em Zero-sum {M}arkov games with stopping and impulsive
  strategies}, Appl. Math. Optim., 9 (1982/83), pp.~1--24.

\bibitem{Stroock_Varadhan}
{\sc D.~W. Stroock and S.~R.~S. Varadhan}, {\em Multidimensional diffusion
  processes}, Classics in Mathematics, Springer-Verlag, Berlin, 2006.
\newblock Reprint of the 1997 edition.

\bibitem{Taksar_1985}
{\sc M.~I. Taksar}, {\em Average optimal singular control and a related
  stopping problem}, Math. Oper. Res., 10 (1985), pp.~63--81.

\bibitem{Touzi_Vieille_2002}
{\sc N.~Touzi and N.~Vieille}, {\em Continuous-time {D}ynkin games with mixed
  strategies}, SIAM J. Control Optim., 41 (2002), pp.~1073--1088 (electronic).

\bibitem{Xu_2007}
{\sc M.~Xu}, {\em Reflected backward {SDE}s with two barriers under
  monotonicity and general increasing conditions}, J. Theoret. Probab., 20
  (2007), pp.~1005--1039.

\bibitem{Yasuda_1985}
{\sc M.~Yasuda}, {\em On a randomized strategy in {N}eveu's stopping problem},
  Stochastic Process. Appl., 21 (1985), pp.~159--166.

\bibitem{Yin_JL_2012}
{\sc J.~Yin}, {\em Reflected backward stochastic differential equations with
  two barriers and {D}ynkin games under {K}nightian uncertainty}, Bull. Sci.
  Math., 136 (2012), pp.~709--729.

\end{thebibliography}
}

\end{document}